\documentclass{article}
\usepackage{arxiv}
\usepackage{amsmath,amssymb}
\usepackage{amsthm}

\usepackage{graphicx}
\graphicspath{ {./images/} }
\usepackage{float}
\usepackage[ruled,vlined]{algorithm2e}
\usepackage[titles]{tocloft}

\usepackage{hyperref}       

\usepackage{booktabs}       
\usepackage{amsfonts}       
\usepackage{nicefrac}       
\usepackage{microtype}      
\usepackage{siunitx}
\usepackage[]{natbib}
\usepackage{appendix}
\usepackage{chngcntr}  
\usepackage{caption}   

\newtheorem{Thm}{Theorem}
\newtheorem{Lemma}{Lemma}
\newtheorem{Coro}{Corollary}
\newtheorem{Pro}{Proposition}
\theoremstyle{definition}
\newtheorem{Def}{Definition}

\newtheorem{remark}{Remark}

\usepackage{geometry}
\usepackage[pagewise]{lineno}
\usepackage{setspace}
\allowdisplaybreaks

\newcommand{\PP}{\mathbb{P}}           
\newcommand{\EE}{\mathbb{E}}            
\newcommand{\RR}{\mathbb{R}}            

\newcommand{\FDP}{\mathrm{FDP}}          
\newcommand{\FDR}{\mathrm{FDR}}          
\newcommand{\PRDS}{\mathrm{PRDS}}        
\newcommand{\BH}{\mathrm{BH}}            
\newcommand{\cp}{\mathrm{cp}}            
\newcommand{\ind}{\mathbf{1}}            
\newcommand{\re}{\mathrm{e}}             
\newcommand{\rd}{\mathrm{d}}             
\newcommand{\test}{\mathrm{test}}        
\newcommand{\rcal}{\mathrm{cal}}         

\newcommand{\bW}{\mathbf{W}}
\newcommand{\bL}{\mathbf{L}}
\newcommand{\bP}{\mathbf{P}}
\newcommand{\bX}{\mathbf{X}}
\newcommand{\bY}{\mathbf{Y}}
\newcommand{\bZ}{\mathbf{Z}}
\newcommand{\bE}{\mathbf{E}}
\newcommand{\bJ}{\mathbf{J}}
\newcommand{\bQ}{\mathbf{Q}}
\newcommand{\bU}{\mathbf{U}}

\newcommand{\cR}{\mathcal{R}}           
\newcommand{\cH}{\mathcal{H}}           
\newcommand{\cN}{\mathcal{N}}           
\newcommand{\cQ}{\mathcal{Q}}           
\newcommand{\cD}{\mathcal{D}}           
\newcommand{\cC}{\mathcal{C}}           
\newcommand{\cP}{\mathcal{P}}           
\newcommand{\cA}{\mathcal{A}}           
\newcommand{\cB}{\mathcal{B}}           
\newcommand{\cV}{\mathcal{V}}           

\title{The Randomized BH: A Generalized Framework Encompassing Conformal and Competition Tests}
\date{}

\author{
  Mingzhou Deng \\
  State Key Laboratory of Mathematical Sciences, Academy of Mathematics and Systems Science\\ Chinese Academy of Sciences, Beijing 100190, China\\
  \texttt{dengmingzhou@amss.ac.cn} \\
  \And
  Yan Fu \\
  State Key Laboratory of Mathematical Sciences, Academy of Mathematics and Systems Science\\ Chinese Academy of Sciences, Beijing 100190, China\\
  \texttt{yfu@amss.ac.cn} \\
}

\begin{document}
\raggedbottom

\maketitle

\begin{abstract}
This paper introduces the randomized Benjamini--Hochberg (BH) procedure, a generalized framework for false discovery rate control that unifies and extends conformal and competition methods for multiple hypothesis testing. The conformal BH and competition procedures are shown to be mutually analogous reformulations, revealing their common core: testing hypotheses through mutual references across hypotheses. This parallelism yields generalizations of the classical BH procedure in three directions---weighted, dominated, and randomized testing problems---each addressing specific gaps between the two families, and leads to the randomized BH framework and embedding theorems showing that classical BH, conformal, and competition tests are special instances of it. The framework relaxes exchangeability for conformal methods and provides competition methods with access to established p-value-based theory. Applied to compound p-values under negative dependence, it explains why exact false discovery rate control fails and yields procedures that integrate heterogeneous tests with global false discovery rate control. Its practical value is further demonstrated on a real single-cell gene-expression data set.
\end{abstract}

\keywords{BH procedure, Competition test, Conformal test, False discovery rate, Randomization}

\section{Introduction}
\label{sec:introduction}

Multiple hypothesis testing is a fundamental problem in statistics: it concerns controlling the cumulative Type I error across many simultaneous tests. The FDR \citep{benjamini1995controlling} was introduced as a criterion for error control, together with the corresponding BH procedure. As a p-value-based method, the BH procedure controls the FDR when the null p-values are super-uniform and satisfy suitable dependence conditions. It is nonetheless limited by its reliance on p-values, which require knowledge of the null distribution or a reliable approximation to it. For complex data or statistics, valid p-values are often infeasible to obtain, restricting the applicability of classical p-value-based methods. Several null-distribution-free approaches have therefore emerged, notably competition methods and conformal methods.

Competition methods include the Target--Decoy Competition (TDC) procedure \citep{elias2007,he2015theoretical,he2022null}, the Knockoff filter \citep{barber2015controlling,weinstein2017power} with its Model-X extension \citep{candes2018panning}, and the Data Splitting and Symmetrized Data Aggregation approaches \citep{dai2023false,du2023false}. They construct pseudo-variables or mirror statistics that preserve the dependence structure of the original variables under the null. Comparing the original test statistics against their constructed counterparts then controls the FDR without specifying the null distribution, and the resulting procedure adapts to complex dependence structures.

Conformal methods originated from the conformal prediction framework \citep{vovk2005} and were adapted for multiple hypothesis testing \citep{mary2022,bates2023}. They generate valid p-values, called conformal p-values, by evaluating the rank of a test statistic within a calibration set. Under exchangeability, they rely solely on the relative ordering of test statistics across the observed and calibration data, requiring no knowledge of the null distribution. They are thus interpretable without knowing the null distribution, which makes them particularly meaningful for black-box models in machine learning.

Competition methods and conformal methods offer distinct advantages yet face notable limitations. Competition methods discard p-values entirely, sacrificing their natural benefits: p-values lie in $[0,1]$, provide a unified scale, and directly quantify evidence against the null hypothesis. Moreover, the raw values of the test statistics carry no intrinsic meaning beyond their relative ordering, and the rejection threshold is a complex random variable without a closed form. Conformal methods, by contrast, require exchangeability to generate valid conformal p-values, which restricts the applicable data and, with the calibration set, reduces the effective sample size for testing.

A closer inspection suggests a structural parallelism between the two. Neither paradigm evaluates test statistics against a null distribution; both evaluate them against other statistics. Conformal methods compare test scores to a calibration set, whereas competition methods compare original variables against constructed counterparts. Since each hypothesis is thus judged by reference to the others, a natural question arises: despite their different formal structures and independent origins, do they share a common theoretical core? Their complementary strengths and weaknesses motivate a unified framework that encompasses both, overcomes their limitations, and lets each benefit from the theoretical tools developed for the other.

The main contributions of this work are summarized as follows.

\textit{Generalizations of the BH procedure.} We establish that the conformal BH and competition procedures are mutually analogous reformulations sharing a common core: testing hypotheses through mutual references across hypotheses. This insight yields generalizations of the classical BH procedure in three directions---weighted BH, dominated BH, and BH for randomized testing problems---each addressing specific gaps between the two families.

\textit{The Randomized BH framework.} We develop the randomized BH framework, with randomized BH-validity at its core, and prove embedding theorems showing that the classical BH procedure, conformal methods, and competition methods are all special instances of it.

\textit{Theoretical and practical applications.} We provide two applications of the framework. First, we analyse compound p-values \citep{barber2025falsediscoveryratecontrol} through the lens of randomized BH, explaining why exact FDR control fails under heterogeneous marginal distributions. Second, we develop a general approach for integrating heterogeneous multiple testing procedures under global FDR control.

\section{Preliminaries}
\label{sec:preliminaries}

\subsection{Notation}

For a positive integer $n$, let $[n]=\{1,\dots,n\}$. For a set $A$, let $|A|$ denote its cardinality and $\cP(A)$ its power set. For $a,b\in\RR$, write $a\vee b=\max\{a,b\}$ and $a\wedge b=\min\{a,b\}$. Let $\RR_+=(0,\infty)$ be the positive half-line, $\RR_*=\RR_+\cup\{0\}=[0,\infty)$ the nonnegative half-line, and $\RR_*^m$ its $m$-dimensional product space. The indicator function is $\ind\{\cdot\}$. For a vector $\bX=(X_1,\dots,X_m)$, let $\bX_{-j}$ denote the subvector with the $j$th component removed, $\bX_S=(X_i)_{i\in S}$ the restriction to $S\subseteq[m]$, and $\bX_{(\cdot)}$ the components in non-decreasing order. For an index $j$, let $\bX^{j}_{(\cdot)}$ denote $\bX$ with the $j$th component removed and the remaining components sorted in non-decreasing order.

The set of hypotheses is $\cH=[m]$, with $\cH_0\subseteq\cH$ the true null hypotheses and $\cH_1=\cH\setminus\cH_0$ the alternatives. Let $m_0=|\cH_0|$, $\pi_0=m_0/m$, and $m_1=m-m_0$. For a rejection set $\cR\subseteq\cH$, the false discovery proportion, the FDR, and the power (average power) are defined as
\[
\FDP=\frac{|\cR\cap\cH_0|}{|\cR|\vee1},\qquad
\FDR=\EE[\FDP],\qquad
\mathrm{Power}=\EE\!\left[\frac{|\cR\cap\cH_1|}{m_1\vee1}\right].
\]

\subsection{The classical BH procedure}

The BH procedure is the most widely used method for FDR control \citep{benjamini1995controlling}. We recall its definition and the conditions under which it guarantees FDR control.

\begin{Def}[BH procedure]
\label{def:BH-procedure}
Let $\bX\in\RR_*^m$ be test statistics associated with hypotheses $\cH=[m]$.
The BH procedure at level $\alpha\in(0,1)$ rejects
$\cR^{\BH}_{\alpha}(\bX)=\{j\in\cH: X_j\leq q^{m,\alpha}_K\}$,
where
\[
K=\sup\{k\in[m]: X_{(k)}\leq q^{m,\alpha}_k\},
\qquad
q^{m,\alpha}_k=\frac{k}{m}\alpha,
\]
with the conventions $\sup\emptyset=0$ and $q^{m,\alpha}_0=0$.
\end{Def}

FDR control of the BH procedure relies on two properties of the test statistics \citep{benjamini2001}: super-uniformity (the p-variable condition) and positive regression dependence on a subset (PRDS).

\begin{Def}[BH-validity]
\label{def:BH-valid}
A random variable $X$ is a \emph{p-variable} if
$\PP\{X\leq t\}\leq t$ for all $t\in[0,1]$.
A random vector $\bX\in\RR_*^m$ is \emph{$\PRDS_{\cN}$} if for any
non-decreasing set $D\subseteq[0,1]^m$, the conditional probability
$\PP\{\bX\in D\mid X_j\leq t\}$ is non-decreasing in $t$ for all
$j\in\cN$.
We say $\bX\in\RR_*^m$ is \emph{BH-valid} with respect to
$\cN\subseteq[m]$ if each $X_j$ ($j\in\cN$) is a p-variable and $\bX$ is $\PRDS_{\cN}$.
\end{Def}

In practice, the BH procedure evaluates p-values only at the discrete thresholds $q^{m,\alpha}_k=\alpha k/m$, so PRDS need only hold on the finite set $\cQ_{m,\alpha}=\{\alpha k/m\}_{k=1}^m$. We denote this restricted condition by $\PRDS_{\cN,\cQ_{m,\alpha}}$. Under BH-validity, the BH procedure controls the FDR at the nominal level.

\begin{Thm}[FDR control for the BH procedure]
\label{thm:BH-controls-FDR}
If $\bP\in\RR_*^m$ are BH-valid with respect to $\cH_0$, then
for any $\alpha\in(0,1)$, the BH procedure $\cR^{\BH}_{\alpha}(\bP)$ controls $\FDR \leq \pi_0\alpha$.
\end{Thm}

\subsection{Conformal methods}

Conformal methods circumvent the need for an explicit null distribution by constructing conformal p-values directly from the data.

\begin{Def}[Conformal p-variables]
\label{def:conformal-p-variables}
Let $\bX\in\RR^m$ be a test set and $\bY\in\RR^n$ a calibration
set. The conformal p-variables $\bP=\bP(\bX,\bY)$ are defined by
\[
P_j = \frac{1+\sum_{i=1}^n\ind\{Y_i\geq X_j\}}{1+n},
\qquad j=1,\ldots,m.
\]
\end{Def}

The conformal p-values form an empirical distribution function with the test statistics as query points and the calibration set as reference. Their validity relies on the exchangeability of the test and calibration sets.

\begin{Def}[Exchangeability]
\label{def:exchangeability}
A sequence of random variables $\bX=(X_1,\ldots,X_m)$ is
exchangeable if for any permutation $\pi$ of $[m]$, the joint
distribution of $(X_{\pi(1)},\ldots,X_{\pi(m)})$ is the same as
that of $(X_1,\ldots,X_m)$.
\end{Def}

\begin{Thm}
\label{thm:conformal-is-BH-valid}
If $\bX_{\cH_0}$ and $\bY$ are exchangeable conditional on
$\bX_{\cH_1}$, and $\bP=\bP(\bX,\bY)$ are conformal p-variables,
then $\bP$ are BH-valid with respect to $\cH_0$.
\end{Thm}

By BH-validity, applying the BH procedure to conformal p-values constructed from exchangeable data guarantees FDR control at the nominal level.

\subsection{Competition methods}

Competition methods construct a synthetic counterpart for each test statistic via permutations, data splitting, or knockoff constructions. Competing the original variables against their counterparts yields, for each hypothesis, a competition label $L_j\in\{0,1\}$ and a competition score $W_j$. The competition procedure is then defined as follows.

\begin{Def}[Competition procedure]
\label{def:competition-procedure}
Let $(\bW,\bL)\in\RR_*^m\times\{0,1\}^m$ be test statistics.
The competition procedure with parameter $r>0$ at level
$\alpha\in(0,1)$ rejects
$\cR^{\cp}_{r,\alpha}(\bW,\bL)=\{j\in\cH: W_j\geq T,\;L_j=1\}$,
where
\[
T=\inf\left\{t\in\{W_j:L_j=1\}:
r\,\frac{1+\sum_{j=1}^m\ind\{W_j\geq t,\;L_j=0\}}
{1\vee\sum_{j=1}^m\ind\{W_j\geq t,\;L_j=1\}}\leq\alpha\right\},
\]
with the convention $\inf\emptyset=+\infty$.
\end{Def}

The additive constant $1$ in the numerator of the threshold is essential for exact FDR control \citep{uri2023}. Without it, the procedure guarantees only the modified FDR \citep{barber2015controlling}.

\begin{Def}[(Sub-)competition statistics]
\label{def:competition-statistics}
The pair $(\bW,\bL)\in\RR_*^m\times\{0,1\}^m$ are
\emph{sub-competition statistics} with respect to parameter $r>0$
and a subset $\cN\subseteq[m]$ if for every $j\in\cN$,
\[
\PP\{L_j=1\mid \bW,\bL_{-j}\}\leq\frac{r}{1+r}.
\]
If the inequality holds with equality for all $j\in\cN$, then
$(\bW,\bL)$ are \emph{competition statistics} with respect to
$r$ and $\cN$.
\end{Def}

Several well-known methods instantiate this condition: the Knockoff filter and Data Splitting approaches typically take $r=1$, while the Target--Decoy Competition procedure relies on independence among test statistics.

\begin{Thm}[FDR control for the competition procedure]
\label{thm:competition-controls-FDR}
If $(\bW,\bL)\in\RR_*^m\times\{0,1\}^m$ are sub-competition
statistics with respect to parameter $r>0$ and $\cH_0$, then
for any $\alpha\in(0,1)$, the competition procedure $\cR^{\cp}_{r,\alpha}(\bW,\bL)$ controls $\FDR \leq \alpha$.
\end{Thm}

The Supplementary Material provides three proofs of Theorem~\ref{thm:competition-controls-FDR}, each offering a distinct perspective on the mechanism of FDR control in competition methods.

\section{The randomized BH framework}
\label{sec:randomized-BH-framework}

\subsection{Mutually analogous reformulations between conformal and competition methods}

First, we reformulate the conformal BH procedure to resemble the competition procedure. Let $\bX\in\RR^m$ be the test statistics, $\bY\in\RR^n$ the calibration data, and $\cD^{\test}$ and $\cD^{\rcal}$ the test and calibration sets. We construct a competition-like procedure that, at level $\alpha\in(0,1)$, rejects
\[
\cR^{*\cp}_{m(1+n)^{-1},\alpha}\bigl((\bX,\bY),\;\ind_{\cD^{\test}}\bigr)
=\{j: X_j\geq T,\;\ind_{\cD^{\test}}(j)=1\},
\]
where the threshold $T$ is defined by
\[
T=\inf\left\{t\in\{X_j\}_{j=1}^m:
\frac{m}{1+n}\,
\frac{1+\sum_{i=1}^n\ind\{Y_i\geq t,\;\ind_{\cD^{\rcal}}(i)=1\}}
{1\vee\sum_{j=1}^m\ind\{X_j\geq t,\;\ind_{\cD^{\test}}(j)=1\}}\leq\alpha\right\}.
\]
The procedure $\cR^{*\cp}_{m(1+n)^{-1},\alpha}$ resembles the competition procedure $\cR^{\cp}_{r,\alpha}$ with $r=m/(1+n)$, in which the labels $L_j$ are replaced by the test-set indicator $\ind_{\cD^{\test}}$. The key difference is that the competition procedure generates the labels $L_j$ at random, whereas the conformal procedure uses the fixed labels $\ind_{\cD^{\test}}$ determined by the data split. Indeed,
\[
\frac{m}{1+n} = \frac{1}{\pi_0}
\,\frac{\sum_{j\in\cH_0}\ind_{\cD^{\test}}(j)}
{1+\sum_{i=1}^n\ind_{\cD^{\rcal}}(i)}=\frac{\widehat{r}}{\pi_0},
\]
where $\widehat{r}$ is the empirical competition parameter, the ratio of the number of nulls in the test set to one plus the size of the calibration set. This also explains in part why competition methods control the FDR at level $\alpha$ without the $\pi_0$ factor.

Second, we reformulate the competition procedure to resemble the conformal BH procedure. For $(\bW,\bL)\in\RR_*^m\times\{0,1\}^m$ and threshold $t\geq 0$, define
$\cR_+(t)=\{j: W_j\geq t,\;L_j=1\}$ and
$\cR_-(t)=\{j: W_j\geq t,\;L_j=0\}$,
and let $\cV_\pm(t)=\cH_0\cap\cR_\pm(t)$. Their cardinalities are $R_\pm(t),V_\pm(t)$. Rewriting the competition threshold $T$ in Definition \ref{def:competition-procedure} gives
\begin{align*}
T&=\inf\left\{t\in\{W_j:L_j=1\}:
r\,\frac{1+R_-(t)}{1\vee R_+(t)}\leq\alpha\right\}\\
&=\inf\left\{t\in\{W_j:L_j=1\}:
\frac{1+R_-(t)}{1+R_-(0)}
\bigg/ \frac{1\vee R_+(0)}{r+rR_-(0)}
\leq\frac{\alpha R_+(t)}{1\vee R_+(0)}\right\}.
\end{align*}
Given $\bP\in\RR_*^m$ and weights $\bE\in\RR_*^m$, the conformal-like
(weighted BH) procedure at level $\alpha\in(0,1)$ rejects
$\cR^{*\BH}_{\alpha}(\bP;\bE)=\{j: P_j/E_j\leq q^{R_+(0),\alpha}_K\}$,
where $K$ follows Definition~\ref{def:BH-procedure} with $m$ replaced by $R_+(0)$:
\[
K=\sup\{k\in[R_+(0)]: P_{(k)}/E_{(k:\bP)}\leq q^{R_+(0),\alpha}_k\},
\]
and $E_{(k:\bP)}$ is the weight corresponding to $P_{(k)}$. Now set
\begin{equation}
  \label{eq:competition-as-conformal1}
  P_j =
  \begin{cases}
    \dfrac{1+R_-(W_j)}{1 + R_-(0)}, & L_j = 1,\\[6pt]
    +\infty,                           & L_j = 0,
  \end{cases}
  \qquad
  E_j = \frac{1\vee R_+(0)}{r+rR_-(0)}.
\end{equation}
When $L_j=0$, hypothesis $j$ is excluded from rejection as its p-value is infinite. An alternative formulation replaces $R_-(0)$ with $V_-(0)$:
\begin{equation}
  \label{eq:competition-as-conformal2}
  P_j =
  \begin{cases}
    \dfrac{1+R_-(W_j)}{1+ V_-(0)}, & L_j = 1,\\[6pt]
    +\infty,                           & L_j = 0,
  \end{cases}
  \qquad
  E_j = \frac{1\vee R_+(0)}{r+rV_-(0)}.
\end{equation}
Since $P_j/E_j$ is the same in (\ref{eq:competition-as-conformal1}) and (\ref{eq:competition-as-conformal2}), the rejection set is unchanged; although $V_-(0)$ is unobservable, the procedure remains operational. Under exchangeability, for any null $j$ with $L_j=1$, $P_j$ in (\ref{eq:competition-as-conformal2}) is super-uniform, but not in (\ref{eq:competition-as-conformal1}), where $R_-(0)$ replaces $V_-(0)$.

\begin{Thm}[Mutually analogous reformulations]
\label{thm:mutual-reformulations}
The conformal BH procedure and the competition procedure admit
the following mutually analogous reformulations:
\[
\cR^{\BH}_{\alpha}
= \cR^{*\cp}_{m(1+n)^{-1},\alpha},\qquad
\cR^{\cp}_{r,\alpha}
= \cR^{*\BH}_{\alpha}.
\]
\end{Thm}

\begin{remark}
The two identities in Theorem~\ref{thm:mutual-reformulations} are exact set equalities. The qualifier ``analogous'' refers not to the identities but to the reformulations' membership in the other family: each identity realizes the reformulated procedure as a member of an extended class, rather than as the standard procedure of that family. The resulting discrepancies are listed below.
\end{remark}

Several structural discrepancies persist between the two
reformulations: (i)~the conformal-like reformulation uses random weights $E_j$ that depend on the data and on the p-values themselves, unlike most existing weighted BH methods \citep{genovese2006,ignatiadis2016,li2019}; (ii)~the contamination of $R_-(t)$ by non-null terms breaks the distributional tractability of $P_j$ in (\ref{eq:competition-as-conformal2}): the ratio $(1+V_-(t))/(1+V_-(0))$ would behave like a conformal p-value on the known finite set $\{k/(1+V_-(0)): k=1,\dots,1+V_-(0)\}$, which would greatly simplify the analysis. However, the structure is lost under the observable $R_-(t)$ version, where non-null contributions distort the numerator (the ranks); (iii)~the competition procedure restricts testing to the random subset $\cR_+(0)$, whereas the conformal framework uses a fixed test set; (iv)~the conformal BH procedure relies on exchangeability, requiring the null test statistics and calibration data to be identically distributed, whereas competition methods impose no such distributional assumption, and although extensions to weighted exchangeability exist \citep{barber2024}, they remain difficult to apply in practice; (v)~conformal p-values satisfy PRDS, whereas competition methods require only a conditional coin-flip property on the labels. These gaps motivate the generalizations introduced in the next subsection and summarized in Table~\ref{tab:gaps-generalizations}.

\begin{table}[htbp]
\centering
\caption{Correspondence between the gaps and the generalizations of the BH procedure.}
\label{tab:gaps-generalizations}
\begin{tabular}{cl}
\toprule
Gap & Generalization \\
\midrule
(i)~Random weights & Weighted BH \\
(ii)~Non-null contamination of calibration & Dominated BH \\
(iii)~Random hypothesis subset & BH for randomized testing problems \\
(iv)~Exchangeability & Randomized super-uniformity \\
(v)~Dependence structure & Randomized $\PRDS$ \\
\bottomrule
\end{tabular}
\end{table}

\subsection{Generalizations of the BH procedure}
\label{subsubsec:dominated-BH}

The first three gaps are addressed by three generalizations of the BH procedure: weighted BH, dominated BH, and BH for randomized testing problems; the remaining two by randomized BH-validity in Section~\ref{sec:rbh}. 

\textbf{Weighted BH.} \citet{ignatiadis2023} examined the weighted BH procedure, treating the random weights $\bE$ as e-values. Given BH-valid p-variables $\bP\in\RR_*^m$ and random weights $\bE\in\RR_+^m$, the procedure at level $\alpha\in(0,1)$ rejects $\cR^{\BH}_{\alpha}(\bP')$, where $\bP'=(P_1/E_1,\dots,P_m/E_m)$. They established that
\[
\FDR\leq \frac{\alpha}{m}\sum_{j\in\cH_0}\EE[E_j]
\]
under any of the following conditions: (i)~$\bE$ is degenerate (non-random); (ii)~$\bP_{\cH_0}$ and $\bE_{\cH_0}$ are mutually independent; or (iii)~$P_j$ and $E_j$ are mutually independent for each $j$, and $\{(P_j,E_j)\}_{j=1}^m$ satisfy the basic independence condition \citep{blanchard2008}. Each of these conditions alone suffices for FDR control, and they share a common mathematical core: knowledge of the weights does not break the BH-validity of the p-values.

\begin{Thm}[Weighted BH procedure and FDR control]
\label{thm:weighted-BH}
Let $\bP\in\RR_*^m$, $\bE$ be random weights. Suppose that, conditional on $\bE$, the p-values $\bP$ are BH-valid with respect to $\cH_0$: for any $j\in\cH_0$, $\PP(P_j\leq t\mid \bE)\leq t$ for all $t\in[0,1]$, and for any non-decreasing set $D\subseteq[0,1]^{m-1}$, $\PP\{\bP_{-j}\in D\mid \bE,\;P_j\leq t\}$ is non-decreasing in $t$. Then for any $\alpha\in(0,1)$,
\[
\FDR = \EE\left[\frac{|\cR^{\BH}_{\alpha}(\bP')\cap\cH_0|}
{1\vee|\cR^{\BH}_{\alpha}(\bP')|}\right]
\leq \frac{\alpha}{m}\sum_{j\in\cH_0}\EE[E_j],
\]
where $\bP'=(P_1/E_1,\dots,P_m/E_m)$.
\end{Thm}

To apply the weighted BH perspective to the competition procedure, we first examine the random weights $E_j$ in (\ref{eq:competition-as-conformal1}) and (\ref{eq:competition-as-conformal2}). In both formulations, hypotheses with $L_j=0$ are assigned $P_j=+\infty$ and never rejected, so their weights are negligible; one may simply set $E_j=0$ when $L_j=0$. The following lemma confirms that the expected total weight on the null hypotheses is appropriately bounded.

\begin{Lemma}[Weights of competition]
\label{lem:weights-of-competition}
If $(\bW,\bL)\in\RR_*^m\times\{0,1\}^m$ are competition statistics
with respect to parameter $r>0$ and $\cH_0$, then
\[
\sum_{j\in\cH_0}\EE\!\left[\frac{1\vee R_+(0)}{r+rR_-(0)}\,
\ind\{L_j=1\}\right]
\leq \sum_{j\in\cH_0}\EE\!\left[\frac{1\vee R_+(0)}{r+rV_-(0)}\,
\ind\{L_j=1\}\right]
\leq m.
\]
\end{Lemma}

\textbf{Dominated BH.} Except in specific experiments such as the TDC procedure in proteomics peptide-spectrum matching, there is no guarantee that all hypotheses assigned label $0$ are truly null. The count $R_-(t)$ in the numerator of (\ref{eq:competition-as-conformal1}) and (\ref{eq:competition-as-conformal2}) then departs from the ideal $V_-(t)$. The super-uniformity of the resulting $P_j$ can still be established, but the PRDS property enjoyed by the $V_-(t)$ version no longer holds. The following theorem addresses this gap, showing that domination of the observable p-values by the ideal ones suffices for FDR control.

\begin{Thm}[Dominated BH procedure and FDR control]
\label{thm:dominated-BH}
Let $\bP,\bQ\in\RR_*^m$ be statistics. Suppose that: (i)~$\PP(\bP\geq \bQ)=1$; (ii)~$\PP(Q_j\leq t)\leq t$ for all $t\in[0,1]$ and $j\in\cH_0$; (iii)~for any non-decreasing set $D\subseteq[0,1]^{m-1}$, $\PP\{\bP_{-j}\in D\mid Q_j\leq t\}$ is non-decreasing in $t$ for all $j\in\cH_0$. Then for any $\alpha\in(0,1)$, the BH procedure $\cR^{\BH}_{\alpha}(\bP)$ controls $\FDR \leq \pi_0\alpha$.
We refer to this condition as \emph{dominated BH-validity} of $\bP$ with respect to $\cH_0$ and $\bQ$.
\end{Thm}

The optimal domination is achieved when $\bQ$ determines $\bP$ uniquely and monotonically, as in the competition construction where $Q_j$ is the ideal p-value generated from $V_-(t)$ rather than $R_-(t)$.

The theorem also suggests a broader perspective on conformal testing. In practice, the calibration set may contain samples that are not exchangeable with the test data, owing to experimental contamination, batch effects, or other sources of heterogeneity. Write $\cD^{\rcal}=\cD^{\rcal}_0+\cD^{\rcal}_1$, where $\cD^{\rcal}_0$ is exchangeable with the test data and $\cD^{\rcal}_1$ is arbitrary (typically unidentifiable), and consider the statistics and weights
\[
  P_j = \frac{1+\sum_{i\in\cD^{\rcal}}\ind\{Y_i\geq X_j\}}{1+|\cD^{\rcal}_0|},\quad Q_j = \frac{1+\sum_{i\in\cD^{\rcal}_0}\ind\{Y_i\geq X_j\}}{1+|\cD^{\rcal}_0|},\quad
  E_j = \frac{1+|\cD^{\rcal}|}{1+|\cD^{\rcal}_0|}.
\]
Here $\PP(\bP\geq\bQ)=1$, and $\bP$ is dominated BH-valid with respect to $\cH_0$ and $\bQ$: conditional on $Y_{\cD^{\rcal}_1}$ and the order statistics of the exchangeable part, $\bP$ is determined by $\bQ$, and $Y_{\cD^{\rcal}_1}$ is independent of $\bQ$. The constant weight $E_j$ absorbs the discrepancy in denominators. Crucially,
\[
P_j' = \frac{P_j}{E_j} = \frac{1+\sum_{i\in\cD^{\rcal}}\ind\{Y_i\geq X_j\}}{1+|\cD^{\rcal}|}
\]
is observable and exactly the standard conformal p-value computed from the full calibration set. Theorems~\ref{thm:weighted-BH} and~\ref{thm:dominated-BH} then give
\[
\FDR=\EE\left[\frac{|\cR^{\BH}_{\alpha}(\bP')\cap\cH_0|}
{1\vee|\cR^{\BH}_{\alpha}(\bP')|}\right]\leq \alpha\pi_0\,\EE\!\left[\frac{1+|\cD^{\rcal}|}{1+|\cD^{\rcal}_0|}\right].
\]
FDR control is thus preserved whenever the exchangeable subset $\cD^{\rcal}_0$ is sufficiently large. Practitioners can apply BH to the standard conformal p-values from the full calibration set, without identifying which samples are exchangeable. If instead the contaminated portion is relatively large, $\bP'$ may be inflated relative to the contamination-free case, reducing power or even eliminating rejections.

\textbf{BH for randomized testing problems.} The competition procedure restricts testing to the random subset $\cR_+(0)$, whereas the standard BH procedure assumes a fixed hypothesis set. The following theorem extends BH to a random hypothesis set; the proof follows from the law of total expectation by conditioning on the random mask $\eta$.

\begin{Thm}[BH for the randomized testing problems]
\label{thm:BH-randomized-testing}
Let $\bP\in\RR_*^m$ and let $\eta\in\cP([m])$ be a random mask. If, conditional on $\eta$,
the subvector $\bP_{\eta}$ is BH-valid with respect to $\eta\cap\cH_0$, then
for any $\alpha\in(0,1)$,
\[
\FDR = \EE\left[\frac{|\cR^{\BH}_{\alpha}(\bP_{\eta})\cap\cH_0|}
{1\vee|\cR^{\BH}_{\alpha}(\bP_{\eta})|}\right]
\leq \alpha\;\EE\!\left[\frac{|\cH_0\cap\eta|}{1\vee|\eta|}\right],
\]
where $\cR^{\BH}_{\alpha}(\bP_{\eta})$ is the BH procedure applied
to the sub-problem $\eta$ with critical values
$q^{|\eta|,\alpha}_k = k\alpha/|\eta|$,
$k=1,\ldots,|\eta|$.
\end{Thm}

The theorem also admits an intuitive interpretation via dominated BH. Extend $\bP_\eta$ to the full $m$-dimensional vector $\widetilde{\bP}$ by setting $\widetilde P_j = P_j$ for $j\in\eta$ and $\widetilde P_j = +\infty$ for $j\notin\eta$, so that $\PP(\bP\leq \widetilde{\bP})=1$. In the notation of Theorem~\ref{thm:dominated-BH}, the untruncated vector $\bP$ plays the role of $\bQ$ (the valid bounding statistics) and $\widetilde{\bP}$ the role of $\bP$. Taking expectation over $\eta$, we have
\[
\PP(\widetilde P_j \leq t) = \EE_\eta\bigl[\PP(\widetilde P_j \leq t\mid\eta)\bigr]
= \EE_\eta\bigl[\PP(P_j \leq t\mid\eta)\,\ind\{j\in\eta\}\bigr]
\leq t\,\PP(j\in\eta),
\]
so the tail of $\widetilde P_j$ is deflated by the selection probability $\PP(j\in\eta)$. This stronger super-uniformity compensates for the more aggressive critical values $q^{|\eta|,\alpha}_k = k\alpha/|\eta|$ used in place of the standard $q^{m,\alpha}_k = k\alpha/m$. Conditional on $\eta$, the BH-validity of $\bP_\eta$ and the domination $\widetilde{\bP}\geq\bP$ yield the FDR bound by Theorem~\ref{thm:dominated-BH}.

\subsection{The randomized BH-validity and FDR control}
\label{sec:rbh}

The remaining two gaps, exchangeability and the PRDS condition, are more fundamental. Conformal BH requires the null test statistics and the calibration data to be identically distributed. Statistics satisfying the same null hypothesis may nonetheless follow different distributions: for example, when testing equality of means across groups, the second moments may differ even under the null, so test statistics based on different sample sizes or variances need not be identically distributed. Moreover, competition-type procedures can involve intricate dependence structures that are difficult to characterize. In the knockoff setting, the design matrix and the $\ell_1$ penalty induce complex finite-sample correlations, making PRDS difficult to establish.

A Bayesian perspective offers a unifying intuition \citep{storey2003,sun2007,efron2008}: even when the hypotheses concern distinct parameters or distributions, a prior renders their marginal distributions identical, and conformal inference requires only that this equality holds under the null. Averaging over the hypotheses, with each null treated as equally likely to be selected, may thus relax exchangeability and PRDS simultaneously. The following definition formalizes this idea.

\begin{Def}[Randomized BH-valid]
\label{def:randomized-BH-valid}
$\bP\in\RR_*^m$ are randomized BH-valid with respect to $\cN$ if the
following hold for an index $J$ uniformly random on $\cN$ and
independent of $\bP$: (i)~$\PP(P_J\leq t)\leq t$ for all $t\in[0,1]$
(marginal super-uniformity under random indexing); and (ii)~for any non-decreasing set $D\subseteq[0,1]^{m-1}$,
$\PP\{\bP^{J}_{(\cdot)}\in D\mid P_J\leq t\}$
is non-decreasing in $t$.
\end{Def}

We refer to condition (i) as \emph{randomized super-uniformity} and condition (ii) as \emph{randomized PRDS}, which can likewise be relaxed to hold only on the grid $\cQ_{m,\alpha}$. Randomized BH-validity replaces the per-index requirements of classical BH-validity with average conditions over $\cH_0$: rather than requiring each null $P_j$ to be a p-variable satisfying PRDS, it requires only that a random null index $J$ satisfy these properties. The effective marginal distribution of the null p-value is then the average $|\cH_0|^{-1}\sum_{j\in\cH_0}\PP(P_j\leq t)$, which can be super-uniform even when some individual $P_j$ are not. Under this formulation, the BH procedure retains FDR control.

\begin{Thm}[FDR control for the randomized BH procedure]
\label{thm:RBH-controls-FDR}
If $\bP\in\RR_*^m$ are randomized BH-valid with respect to
$\cH_0$, then for any $\alpha\in(0,1)$, the BH procedure $\cR^{\BH}_{\alpha}(\bP)$ controls $\FDR \leq\pi_0 \alpha$.
\end{Thm}

The validity of this theorem rests on a simple fact: FDR, power, and the total number of rejections depend only on how many null and how many non-null hypotheses are rejected, not on which ones, so each hypothesis carries equal weight. Randomization can thus be understood as reassigning rejections uniformly among the null hypotheses, which preserves the FDR. Equivalently, randomized BH-validity amounts to requiring that after a uniformly random permutation of the null hypotheses, the permuted vector $\widetilde{\bP}$ with $\widetilde P_i = P_{J_i}$ be classically BH-valid.

In a conformal procedure, this permutation can be realized by permuting the underlying data. Beyond implementing randomization, this confers an additional advantage: it induces exchangeability where none was present. The following lemma makes this precise.

\begin{Lemma}[Random permutation induces exchangeability]
\label{lem:random-permutation-to-exchangeable}
Let $Z_i = X_{J_i}$ for $i=1,\dots,p$ and $Z_i = X_i$ for
$i=p+1,\dots,p+q$, where $\mathbf{J}$ is a uniformly random
permutation of $[p]$. Then $(Z_1,\dots,Z_p)$ are exchangeable
conditional on $(Z_{p+1},\dots,Z_{p+q})$.
\end{Lemma}

This lemma shows that, without any assumptions on the original data, a random permutation of the null hypotheses yields an exchangeable sequence. For a conformal procedure with fixed test and calibration sets, however, permuting data across them would break the separation between test and calibration units. The required randomization can instead be achieved by randomizing the partition itself: a uniformly random split of all units into two subsets induces the same distribution as does a uniformly random permutation of all units. Here, for a generic set $\cN$, a uniformly random split $\cN=\cA+\cB$ means that
\[
\PP(\cA=A\mid |\cA|=k)=\binom{|\cN|}{k}^{-1},
\qquad\forall\;A\in\cP(\cN),\;|A|=k.
\]
That is, conditionally on the size of $\cA$, every subset of that size is equally likely to be selected. This yields the following lemma.

\begin{Lemma}[Random split induces random permutation]
\label{lem:random-split-to-random-permutation}
Let $[m]$ be split uniformly at random into two disjoint subsets
$[m]=\cA+\cB$. Define $Z_i = X_{J_i}$ for $i\in\cA$ and
$Z_i = X_{K_i}$ for $i\in\cB$, where $\mathbf{J}$ is a random
permutation of $\cA$ and $\mathbf{K}$ is a random permutation of
$\cB$. Then $\bZ \overset{d}{=} (X_{S_1},\dots,X_{S_m})$,
where $\mathbf{S}$ is a random permutation of $[m]$.
\end{Lemma}

Lemmas~\ref{lem:random-permutation-to-exchangeable} and~\ref{lem:random-split-to-random-permutation} yield the following corollary of Theorem~\ref{thm:RBH-controls-FDR}.

\begin{Coro}[Random split induces randomized BH-validity]
\label{cor:random-split-to-RBH-valid}
Let $\cH_0$ be split by a uniformly random partition into
$\cH_0=\cH_0^{\test} + \cH_0^{\rcal}$, and set
$\cH^{\test} = \cH_0^{\test} + \cH_1$,
$\cH^{\rcal} = \cH_0^{\rcal}$.
Then the BH procedure applied to the conformal p-values controls the FDR in the sense that
\[
\FDR=\EE\left[\frac{|\cR^{\BH}_{\alpha}(\bP(\bX_{\cH^{\test}},\bX_{\cH^{\rcal}}))\cap\cH_0|}
{1\vee|\cR^{\BH}_{\alpha}(\bP(\bX_{\cH^{\test}},\bX_{\cH^{\rcal}}))|}\right]\leq\alpha \EE\left[\frac{|\cH_0^{\test}|}{1\vee|\cH^{\test}|}\right].
\]
\end{Coro}

We sketch a heuristic proof. The conformal p-values are constructed in the standard manner on the randomized sub-problem $\cH^{\test}$, using $\cH^{\rcal}$ as the calibration set. Those for indices outside $\cH^{\test}$ are set to $+\infty$ to avoid the indexing ambiguity of a random hypothesis set. This is the domination-based perspective through which Theorem~\ref{thm:BH-randomized-testing} treats randomized testing problems. Conditional on the split, the conformal p-values are not BH-valid, but conditional on its size, they are randomized BH-valid.

\begin{remark}[randomized BH as a guideline for constructing multiple testing procedures]
\label{rmk:RBH-guideline}
The preceding argument illustrates a general principle: an independent random split, constructed through any effective mechanism, enables multiple hypothesis testing. For example, one may run an independent test and treat its outcome as the partition. One then constructs conformal p-values on the induced sub-problem and applies the BH procedure, thereby controlling the FDR. Randomized BH thus provides a guideline for constructing new multiple testing procedures: any method for which randomized BH-validity can be established yields a valid procedure. Our analysis of how random permutations and random splits induce exchangeability already provides one such route.
\end{remark}

The three generalizations of Section~\ref{subsubsec:dominated-BH} carry over directly to the randomized BH setting, since they can be carried out in the permuted form $\widetilde{\bP}$ introduced above.

\subsection{Embedding theorems}
\label{subsec:embedding-theorems}

We now show that the randomized BH framework encompasses the main paradigms considered in this paper: the following theorems establish that conformal tests, classical BH procedures, and competition tests are each special instances of it.

We first show that conformal p-values are randomized BH-valid. This follows most directly, as they have the most desirable properties, including identical marginal distributions.

\begin{Thm}[Embedding conformal tests into randomized BH]
\label{thm:embedding-conformal-into-RBH}
(i)~If $\bP\in\RR_*^m$ are BH-valid with respect to $\cH_0$
and have identical marginal distributions, then $\bP$ are randomized
BH-valid with respect to $\cH_0$.
(ii)~Consequently, if $\bX_{\cH_0}$ and $\bY$ are exchangeable
conditional on $\bX_{\cH_1}$, then
conformal p-values $\bP=\bP(\bX,\bY)$ are randomized BH-valid
with respect to $\cH_0$.
\end{Thm}

Under identical marginal distributions, randomized PRDS reduces to PRDS directly. In most p-value-based tests, such as point-null parametric tests, the null p-values are exactly uniform and hence identically distributed, so BH-validity automatically implies randomized BH-validity. If the condition is relaxed to super-uniformity without identical distributions, the implication no longer holds directly. Classical BH nevertheless remains a special case of the framework via the dominated formulation, as the following theorem shows.

\begin{Thm}[Embedding classical BH into randomized BH]
\label{thm:embedding-classical-BH-into-RBH}
If $\bP\in\RR_*^m$ are BH-valid with respect to $\cH_0$, then
there exists a random vector $\bQ$ such that
$\bP$ are dominated randomized BH-valid with respect to $\cH_0$
and $\bQ$. That is, $\PP(\bP\geq\bQ)=1$, and for $J$ uniform on $\cH_0$ and independent of $\bP$:
(i)~$\PP(Q_J\leq t)\leq t$ for all $t\in[0,1]$; and
(ii)~for any non-decreasing set $D\subseteq[0,1]^{m-1}$,
$\PP\{\bP^{J}_{(\cdot)}\in D\mid Q_J\leq t\}$ is
non-decreasing in $t$.
\end{Thm}

Recall from the previous subsection that a uniformly random split of the null hypotheses provides the randomization needed for randomized BH. For competition methods this split is not an additional design choice but an intrinsic feature of the procedure: the competition labels already partition the null hypotheses uniformly at random. The following lemma makes this precise and establishes the key link allowing competition tests to be embedded into the framework.

\begin{Lemma}[Competition labels induce random split]
\label{lem:competition-to-random-split}
If $(\bW,\bL)\in\RR_*^m\times\{0,1\}^m$ are competition statistics
with respect to parameter $r>0$ and $\cH_0$, then the competition labels
$\bL$ induce a uniformly random split of $\cH_0$ into
$\cH_0 = \cV_+(0) + \cV_-(0)$.
\end{Lemma}

This lemma and Corollary~\ref{cor:random-split-to-RBH-valid} now show that competition tests are randomized BH-valid.

\begin{Thm}[Embedding competition tests into randomized BH]
\label{thm:embedding-competition-into-RBH}
Let $(\bW,\bL)\in\RR_*^m\times\{0,1\}^m$ be competition statistics
with respect to parameter $r>0$ and $\cH_0$. Let $\bP$ and $\bE$ be as defined in (\ref{eq:competition-as-conformal2}), and set $\bP' = (P_1/E_1,\ldots,P_m/E_m)$. Define
\[
Q_j =
\begin{cases}
\dfrac{1+V_-(W_j)}{1+V_-(0)}, & j\in\cR_+(0),\\[10pt]
+\infty, & j\in\cR_-(0).
\end{cases}
\]
Conditional on $V_+(0)$, $\bP$ are dominated randomized BH-valid with respect to $\cH_0$ and $\bQ$. Consequently, the BH procedure $\cR^{\BH}_{\alpha}(\bP')$ controls $\FDR \leq \alpha$.
\end{Thm}

FDR control for the competition procedure was already established in Theorem~\ref{thm:competition-controls-FDR}, with the mutual reformulations in Theorem~\ref{thm:mutual-reformulations}. The derivation from the randomized BH framework nonetheless offers a novel proof and, more importantly, a fresh perspective on the structural reasons underlying the result. The present theorem covers only the standard competition procedure; the sub-competition case is discussed in the Supplementary Material.

The three embedding theorems thus show that conformal, classical BH, and competition tests are all special instances of the randomized BH framework, answering the question that motivated this work.

\section{Theoretical applications of the randomized BH framework}
\label{sec:applications}

\subsection{Compound p-values and randomized BH framework}

The first application concerns compound p-values. The randomized BH perspective explains why exact FDR control fails under heterogeneous marginal distributions and yields an explicit FDR bound under a specified model.

\citet{wang2022} proposed the eBH method based on e-values and established its FDR control, requiring only that the input statistics be compound e-values \citep{ramdas2025}. This raises the natural question of whether an analogous theory holds for compound p-values, a notion introduced by \citet{armstrong2025}, who termed them \emph{average p-values} in the original work.

\begin{Def}[Compound p-variables]
\label{def:compound-p-variables}
$\bP\in\RR_*^m$ are compound p-variables with respect to $\cN\subseteq[m]$ if $\sum_{j\in\cN}\PP(P_j\leq t)\leq mt$ for all $t\in[0,1]$.
\end{Def}

\citet{barber2025falsediscoveryratecontrol} established finite-sample FDR control theory for compound p-variables.

\begin{Thm}[FDR control for compound p-variables]
\label{thm:ren_compound}
Let $\bP\in\RR_*^m$ be compound p-variables with respect to $\cH_0$. Then for any $\alpha\in(0,1)$, the BH procedure $\cR^{\BH}_{\alpha}(\bP)$ controls the FDR as follows: (i)~under independence, $\FDR\leq 1.93\alpha$; (ii)~under arbitrary dependence, $\FDR\leq \alpha h_m$, where $h_m=\sum_{k=1}^m 1/k$; (iii)~under independence and $\cH_0=\cH$, $\FDR\leq \alpha+2\alpha^2$.
\end{Thm}

This inflation is not an artefact of the proof but unavoidable: \citet{barber2025falsediscoveryratecontrol} established matching lower bounds showing that the enlarged constants are necessary under the compound p-value condition. We restate the lower bounds in Proposition~\ref{prop:counterexample-compound} for completeness.

\begin{Pro}[Lower bounds for compound p-variables]
\label{prop:counterexample-compound}
Let $\bP\in\RR_*^m$ be compound p-variables with respect to $\cH_0$.
The following lower bounds can be attained: (i)~there exist independent $\bP$ with $\FDR \geq 7\alpha/6$; (ii)~there exist $\bP$ with $\PRDS_{\cH_0}$ and $\FDR \geq (3/8)\alpha h_m$, where $h_m = \sum_{i=1}^m i^{-1}$; (iii)~there exist independent $\bP$ under the global null $\cH_0=\cH$ with $\FDR \geq \alpha + \alpha^2/4$.
\end{Pro}

These bounds were obtained by constructive proofs. Theorem~\ref{thm:embedding-conformal-into-RBH} explains this through its proof structure. Compound p-variables differ from randomized super-uniformity only by a factor of $\pi_0$, which does not account for the inflated bounds; the gap must therefore originate from the failure of randomized PRDS, and the mechanism is a selection bias. The compound condition $\sum_{j\in\cH_0}\PP(P_j\leq t)\leq mt$ neither requires identical marginal distributions nor is dominated by any p-variable construction with identical distributions. When the marginal distributions differ, the uniformly selected index $J$ is no longer uniform after conditioning on $P_J\leq t$, since
\begin{align*}
  \PP\{\bP^{J}_{(\cdot)}\in D\mid P_J\leq t\}=\sum_{j\in\cH_0}\PP\{\bP^{j}_{(\cdot)}\in D\mid P_j\leq t\}\frac{\PP(J=j)\PP(P_j\leq t)}{\sum_{k\in\cH_0}\PP(J=k)\PP(P_k\leq t)}.
\end{align*}
Because $\PP(J=j)=1/|\cH_0|$, the weight for each null index $j$ is proportional to $\PP(P_j\leq t)$, and these weights vary with $t$, which can break the monotonicity required for randomized PRDS. In the Supplementary Material we show that the counterexample fails randomized PRDS \citep{barber2025falsediscoveryratecontrol}. Identical marginal distributions are sufficient but not necessary: what matters is that the weights be independent of $t$, which holds precisely when $\PP(P_j\leq t)=c_j f(t)$ for some function $f$ and constants $c_j>0$, with identical marginals corresponding to $c_j=c$.

Random selection of a null index is equivalent to a random permutation of the hypotheses, which tends to induce negative dependence.

\begin{Lemma}[Weak negative dependence is preserved under permutation]
\label{lem:permutation-to-WND}
A random vector $\bX\in\RR^m$ is \emph{weakly negatively dependent}
(WND) if
\[
\PP\Bigl(\bigcap_{k\in A}\{X_k\leq x\}\Bigr)
\leq \prod_{k\in A}\PP(X_k\leq x)
\]
for all $A\subseteq[m]$ and $x\in\RR$. If $\bX$ are WND, and
$Z_i = X_{J_i}$ for $i\in[m]$ where $\mathbf{J}$ is a random
permutation of $[m]$, then $\bZ$ are WND.
\end{Lemma}

PRDS is a positive dependence property, and randomly permuting null p-values with heterogeneous marginal distributions induces strict negative dependence, which destroys PRDS. FDR control under negative dependence is nevertheless well studied \citep{chi2024}, which motivates the following result.

\begin{Thm}[FDR control under WND for compound p-variables]
\label{thm:WND-compound-FDR}
Let $\bP\in\RR_*^m$ be randomized super-uniform and WND with respect to $\cH_0$. Then for any $\alpha\in(0,1)$, the BH procedure $\cR^{\BH}_{\alpha}(\bP)$ controls the FDR as follows:
\[
\FDR \leq \alpha + \alpha\int_{\alpha}^{1}\frac{1}{x^2}\,
\min\!\Bigl\{1,\; x+2x^2+\frac{9}{2}x^3
+ \sum_{k=4}^{\infty}\frac{(\re x)^k}{\sqrt{2\pi k}}\Bigr\}\rd x.
\]
\end{Thm}

The theorem requires only randomized super-uniformity rather than individual super-uniformity, thereby admitting compound p-values. More importantly, the dependence condition is confined to the null set: weak negative dependence is required only in the marginal sense, among the p-values in $\cH_0$, and no condition is imposed on their dependence with those in $\cH_1$. This relaxation is possible because the proof relies on the FDR-Linking theorem \citep{su2018}. How the strength of the negative dependence affects the FDR bound is discussed in the Supplementary Material.

\subsection{Integration of heterogeneous multiple testing procedures}

The second application combines heterogeneous multiple testing procedures under global FDR control, reconciling disparate sources of evidence within a single procedure.

Combining p-value-free procedures naively can severely erode power and may even invalidate FDR control: unlike p-values, their statistics need not share a common scale and may lack intrinsic interpretability (Section~\ref{sec:introduction}). The following counterexample illustrates this failure for the simple strategy of taking the union of rejection sets from two competition procedures; the explicit construction is deferred to the Supplementary Material.

Specifically, there exist two independent competition statistics $(\bW^1,\bL^1)$ and $(\bW^2,\bL^2)$ with respect to parameters $r_1,\cH^1_0$ and $r_2,\cH^2_0$, such that the FDR of the pooled rejection set exceeds the nominal level:
\[
\FDR = \EE\left[
\frac{\bigl|(\cR^{\cp}_{r_1,\alpha}(\bW^1,\bL^1)\cap\cH^1_0)
\;\cup\;(\cR^{\cp}_{r_2,\alpha}(\bW^2,\bL^2)\cap\cH^2_0)\bigr|}
{1\vee\bigl|(\cR^{\cp}_{r_1,\alpha}(\bW^1,\bL^1))
\;\cup\;(\cR^{\cp}_{r_2,\alpha}(\bW^2,\bL^2))\bigr|}
\right] > \alpha.
\]

The randomized BH framework provides a natural remedy: embedding each procedure into it yields randomized BH-valid p-values that share a common scale and can be pooled under a single BH procedure while retaining FDR control. The following theorem formalizes this construction, the direct analogue of pooling classical p-values.

\begin{Thm}[Integration of two independent randomized BH procedures]
\label{thm:integration-RBH-controls-FDR}
Let $\bP^1\in\RR_*^{m^1}$ and $\bP^2\in\RR_*^{m^2}$ be independent
and randomized BH-valid with respect to $\cH^1_0$ and $\cH^2_0$,
respectively. Define
$\bP^{\mathrm{full}} = (\bP^1,\bP^2)\in\RR_*^{m^1+m^2}$ and
$\cH^{\mathrm{full}}_0 = \cH^1_0 + \cH^2_0$.
Then for any $\alpha\in(0,1)$, the BH procedure $\cR^{\BH}_{\alpha}(\bP^{\mathrm{full}})$ controls $\FDR \leq \alpha$.
\end{Thm}

Owing to the selection bias, the concatenation of two independent randomized BH-valid statistics is not in general randomized BH-valid; the following remark shows that this global property is not required for FDR control.

\begin{remark}[Integration does not require global randomized BH-validity]
\label{rmk:integration-RBH-decomposition}
In fact, suppose $\cH$ is partitioned into blocks $\cH=\cup_{i=1}^d\cH^i$ such that, for each $i\in[d]$, $\bP_{\cH^i}$ is randomized BH-valid with respect to $\cH_0\cap\cH^i$ conditional on $\bP_{\cup_{j\neq i}\cH^j}$. Then the BH procedure applied to $\bP$ controls the FDR at any level $\alpha\in(0,1)$.
\end{remark}

The embedding theorems of Section~\ref{subsec:embedding-theorems} provide a principled way to bring p-value-free procedures back into the p-value framework: each procedure is embedded to produce randomized BH-valid p-values, which are pooled and subjected to a single BH procedure. The following two algorithms illustrate this approach.

\begin{algorithm}[htbp]
\caption{CP-CP integration}
\LinesNotNumbered
\label{alg:cp-cp-integration}
\emph{Input:} $(\bW^1,\bL^1)\in\RR_*^{m^1}\times\{0,1\}^{m^1}$, $(\bW^2,\bL^2)\in\RR_*^{m^2}\times\{0,1\}^{m^2}$, $r_1,r_2>0$, $\alpha\in(0,1)$. \emph{Output:} Rejection sets $\cR^1,\cR^2$.
\begin{enumerate}
  \item For $g=1,2$, compute for each $j=1,\ldots,m^g$:
        $X^g_j = (r_g+r_gR^g_-(W^g_j))/m^g$ if $j\in\cR^g_+(0)$, and $X^g_j=+\infty$ if $j\in\cR^g_-(0)$,
        where $R^g_-(t)=\sum_{j=1}^{m^g}\ind\{W^g_j\geq t,\,L^g_j=0\}$ and $\cR^g_+(t)=\{j: W^g_j\geq t,\,L^g_j=1\}$.
  \item Find the largest $K\in\{1,\ldots,m^1+m^2\}$ such that $\sum_{g=1}^2\sum_{j=1}^{m^g} \ind\{X^g_j \leq q^{m^1+m^2,\alpha}_K\} \geq K$.
  \item Return $\cR^g = \{j: X^g_j \leq q^{m^1+m^2,\alpha}_K\}$, $g=1,2$.
\end{enumerate}
\end{algorithm}

Algorithm~\ref{alg:cp-cp-integration} integrates two procedures of the same competition type. The embedding approach applies equally when one procedure is parametric and the other non-parametric, as the following algorithm, which merges a competition procedure with a conventional p-value-based one, illustrates.

\begin{algorithm}[htbp]
\caption{CP-P integration}
\LinesNotNumbered
\label{alg:cp-p-integration}
\emph{Input:} $(\bW,\bL)\in\RR_*^{m^1}\times\{0,1\}^{m^1}$, $\bX^2\in\RR_*^{m^2}$, $r>0$, $\alpha\in(0,1)$. \emph{Output:} Rejection sets $\cR^1,\cR^2$.
\begin{enumerate}
  \item Compute for each $j=1,\ldots,m^1$:
        $X^1_j = (r+rR_-(W_j))/m^1$ if $j\in\cR_+(0)$, and $X^1_j=+\infty$ if $j\in\cR_-(0)$,
        where $R_-(t)=\sum_{j=1}^{m^1}\ind\{W_j\geq t,\,L_j=0\}$ and $\cR_+(t)=\{j: W_j\geq t,\,L_j=1\}$.
  \item Find the largest $K\in\{1,\ldots,m^1+m^2\}$ such that $\sum_{g=1}^2\sum_{j=1}^{m^g}\ind\{X^g_j \leq q^{m^1+m^2,\alpha}_K\} \geq K$.
  \item Return $\cR^g = \{j: X^g_j \leq q^{m^1+m^2,\alpha}_K\}$, $g=1,2$.
\end{enumerate}
\end{algorithm}

The following theorem guarantees FDR control for both procedures.

\begin{Thm}[Integration controls global FDR]
\label{thm:integration-controls-FDR}
The integrated procedures described in
Algorithms~\ref{alg:cp-cp-integration}--\ref{alg:cp-p-integration}
control the global FDR under the following conditions:
\begin{enumerate}
\renewcommand{\labelenumi}{(\roman{enumi})}
  \item For Algorithm~\ref{alg:cp-cp-integration}:
        $(\bW^1,\bL^1)$ and $(\bW^2,\bL^2)$ are independent and
        are competition statistics with respect to $r_1,\cH^1_0$
        and $r_2,\cH^2_0$, respectively;
  \item For Algorithm~\ref{alg:cp-p-integration}:
        $(\bW,\bL)$ and $\bX^2$ are independent, $(\bW,\bL)$ are
        competition statistics with respect to $r,\cH^1_0$, and
        $\bX^2$ are randomized BH-valid with respect to $\cH^2_0$.
\end{enumerate}
Under either condition, for any $\alpha\in(0,1)$,
\[
\FDR = \EE\left[
\frac{|(\cR^1\cap\cH^1_0)\cup(\cR^2\cap\cH^2_0)|}
{1\vee|\cR^1\cup\cR^2|}
\right] \leq \alpha.
\]
\end{Thm}

In the integration algorithms we replace the denominator $R_+(0)$ of the original competition procedure (\ref{eq:competition-as-conformal2}) by the total number of hypotheses $m$. This is immaterial in the single-problem setting, since the constant factor $m/R_+(0)$ is absorbed into the BH threshold through the equivalence $cX_j \leq q^{m,\alpha}_k \iff X_j \leq q^{cm,\alpha}_k$, but necessary in the integrated setting: $R_+(0)$ is usually positively correlated with $V_+(0)$, and by the rearrangement inequality this association between a group's random size and its error count can push the global FDR above $\alpha$, whereas the fixed total $m$ weights each group deterministically. Details are given in the Supplementary Material.

The replacement also suggests a flexible form of group weighting. To reject more hypotheses from $\cH^1$ and fewer from $\cH^2$, multiply the transformed statistics of $\cH^2$ by $c>1$ and run BH with the effective size $m^1+c^{-1}m^2$. Since $cX^2_j \leq q^{m^1+c^{-1}m^2,\alpha}_k$ is equivalent to $X^2_j \leq q^{cm^1+m^2,\alpha}_k < q^{m^1+m^2,\alpha}_k$, the effective threshold for $\cH^2$ is tighter than the baseline $q^{m^1+m^2,\alpha}_k$. Consequently $\cH^1$ receives relatively more rejections and $\cH^2$ relatively fewer, as desired. In the Supplementary Material we analyse this scheme further, showing that FDR control remains in force under various weightings, with numerical results on power.

When the multiplier is random rather than constant, FDR control may still be attainable: the Supplementary Material shows that retaining $R_-(0)$ as the denominator preserves it.

Finally, the randomized BH procedures derived from competition procedures can accommodate arbitrary inter-group dependence, a property that generic randomized p-variables do not possess. The following theorem makes this precise.

\begin{Thm}[FDR control under arbitrary inter-group dependence]
\label{thm:cp_cp_any_dep}
Let $g=1,2$, $(\bW^g,\bL^g)$ be competition statistics for $\cH^g_0\subseteq\cH^g$ with parameters $r_g$. Then there exists a function $c(r)$ such that for any $\alpha\in(0,1)$, the integrated test in Algorithm~\ref{alg:cp-cp-integration} with $X^g_j$ replaced by $X^g_j/c(r_g)$ controls FDR even when $(\bW^1,\bL^1)$ and $(\bW^2,\bL^2)$ have arbitrary inter-group dependence. In particular, one may take $c(1)=0.5201$.
\end{Thm}

The two-problem setting is for simplicity only: all the results above extend directly to any number of testing problems, and the function $c(r)$ in Theorem~\ref{thm:cp_cp_any_dep} does not depend on their number, so it applies equally when integrating any number of competition methods with arbitrary inter-group dependence. The Supplementary Material identifies an additional property of the transformed statistics that ensures FDR control under arbitrary dependence, and provides further details on $c(r)$.

\section{Simulations}
\label{sec:simulations}

\subsection{Conformal tests with prescreened calibration sets}

We give a simple example of the prescreening construction of Section~\ref{subsubsec:dominated-BH}. Using data splitting, we generate two independent sets of random variables for the prescreening and testing stages. For $i=1,\dots,m$, $X^{\mathrm{pre}}_i, X_i \sim N(\mu_i, 1)$, with $m=2000$, $m_0=\pi_0 m$, and $\mu_i=\beta\ind\{i>m_0\}$. The test and calibration sets are obtained by prescreening the first-stage variables at threshold $\tau=0.5$:
\[
\cD^{\test} = \bigl\{i\in[m]: |X^{\mathrm{pre}}_i|\geq 0.5\bigr\}, \qquad
\cD^{\rcal} = \bigl\{i\in[m]: |X^{\mathrm{pre}}_i| < 0.5\bigr\}.
\]
We then construct the conformal p-values $\bP=\bP(\bX_{\cD^{\test}},\bX_{\cD^{\rcal}})$ and apply the BH procedure $\cR^{\BH}_{\alpha}(\bP)$ on the test set. Over $300$ replications, we estimate the FDR and power by the average false discovery proportion and the true discovery proportion (TDP), respectively.
\begin{figure}[htbp]
\centering
\includegraphics[width=368.97pt]{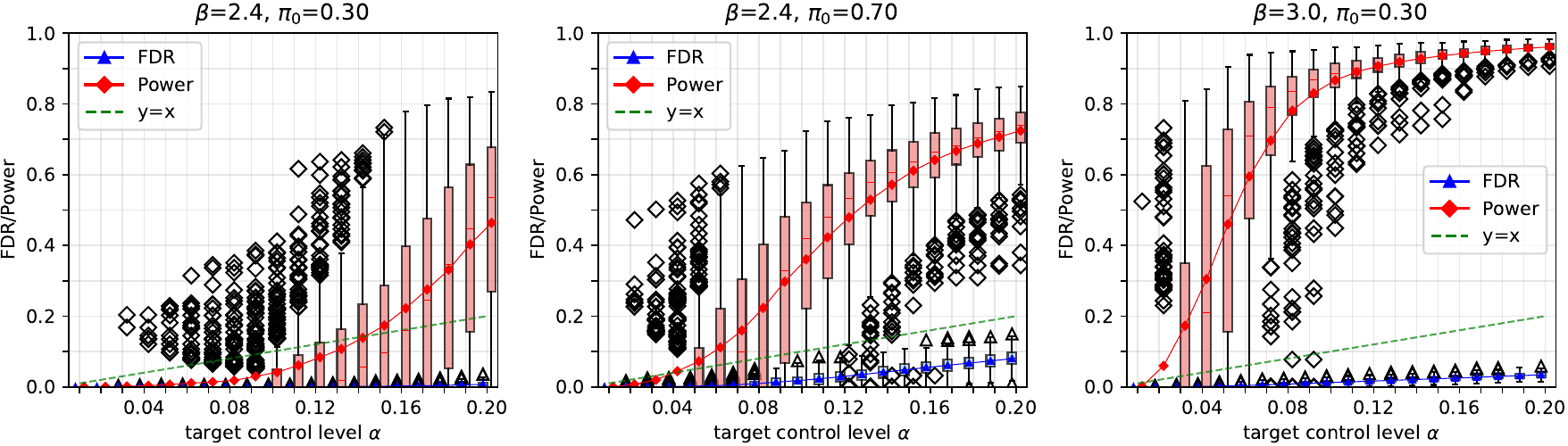}
\caption{Boxplots of FDR and power under different target levels $\alpha$ and parameter settings $(\beta,\pi_0)$.}
\label{fig:prescreen_cal_box3}
\end{figure}

The results, summarized in Fig.~\ref{fig:prescreen_cal_box3}, confirm that prescreening the calibration set, followed by conformal p-values and the BH procedure, achieves FDR control without a $\pi_0$ correction \citep{storey2002}, provided that the prescreening is accurate enough to keep the calibration set from being over-contaminated with alternatives. Further results are deferred to the Supplementary Material.

\subsection{Numerical results for the integration procedures}

Both the conformal p-values and the competition statistics are built from normal random variables with a common structure. For alternatives ($i>m_0$), $\mu_i=\beta$ and $\sigma_i=1$, whereas for nulls ($i\leq m_0$), $\mu_i=0$ and the $\sigma_i$ follow a periodic weight vector.

For randomized p-values, we generate $X_i \sim N(\mu_i, \sigma_i^2)$ with weights $\boldsymbol{w}\in\RR_+^b$, where $\sigma_i = w_k$ if $i\equiv k-1 \pmod{b}$ and $i\leq m_0$. Let $\bY$ permute the first $m_0$ components of $\bX$ uniformly at random and leave the remaining $m-m_0$ unchanged, and set $\cD^{\rcal}=[m^{\rcal}]$, $\cD^{\test}=[m]\setminus\cD^{\rcal}$ with $m^{\rcal}<m_0$. The conformal p-values $\bP=\bP(\bY_{\cD^{\test}},\bY_{\cD^{\rcal}})$ are then randomized BH-valid.

For competition statistics, we generate original variables $X^{\mathrm{org}}_i\sim N(\mu_i,\sigma_i^2)$ and pseudo-variables $X^{\mathrm{pse}}_i\sim N(0,\sigma_i^2)$ for $i=1,\dots,m$, with null variances following the periodic weight vector $\boldsymbol{w}^{\mathrm{cp}}\in\RR_+^{b^{\mathrm{cp}}}$ as above. We then set $W_i=\max(|X^{\mathrm{org}}_i|,|X^{\mathrm{pse}}_i|)$ and $L_i=\ind\{|X^{\mathrm{org}}_i|\geq |X^{\mathrm{pse}}_i|\}$, so that $(\bW,\bL)$ are competition statistics with $r=1$.

\begin{figure}[htbp]
\centering
\includegraphics[width=333.67pt]{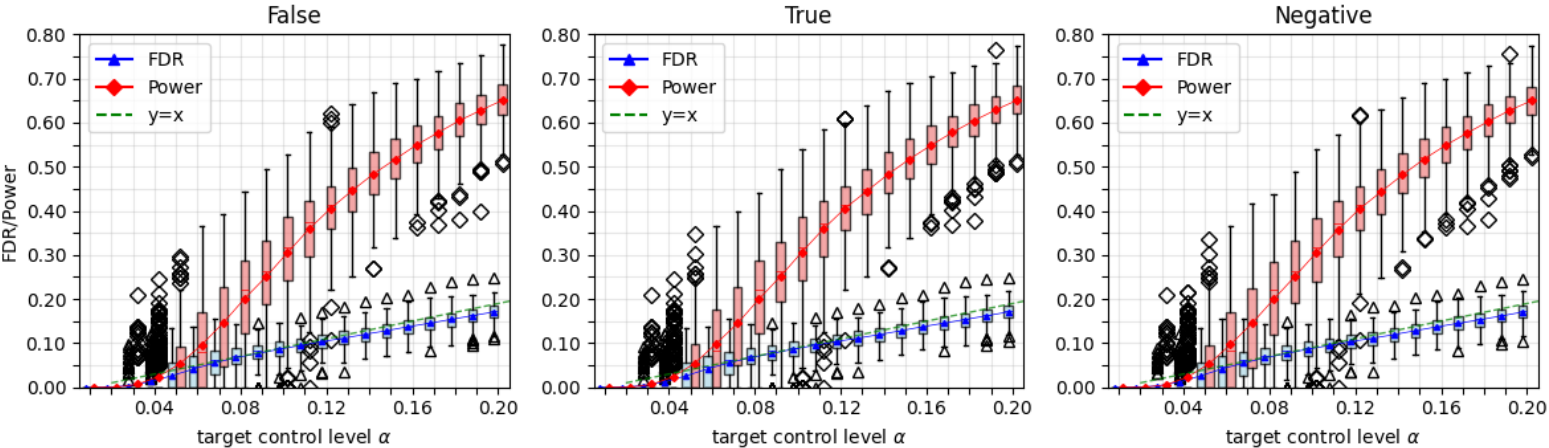}
\caption{Integration test of two competition procedures.}
\label{fig:cp_cp_integration_box}
\end{figure}

\begin{figure}[htbp]
\centering
\includegraphics[width=333.67pt]{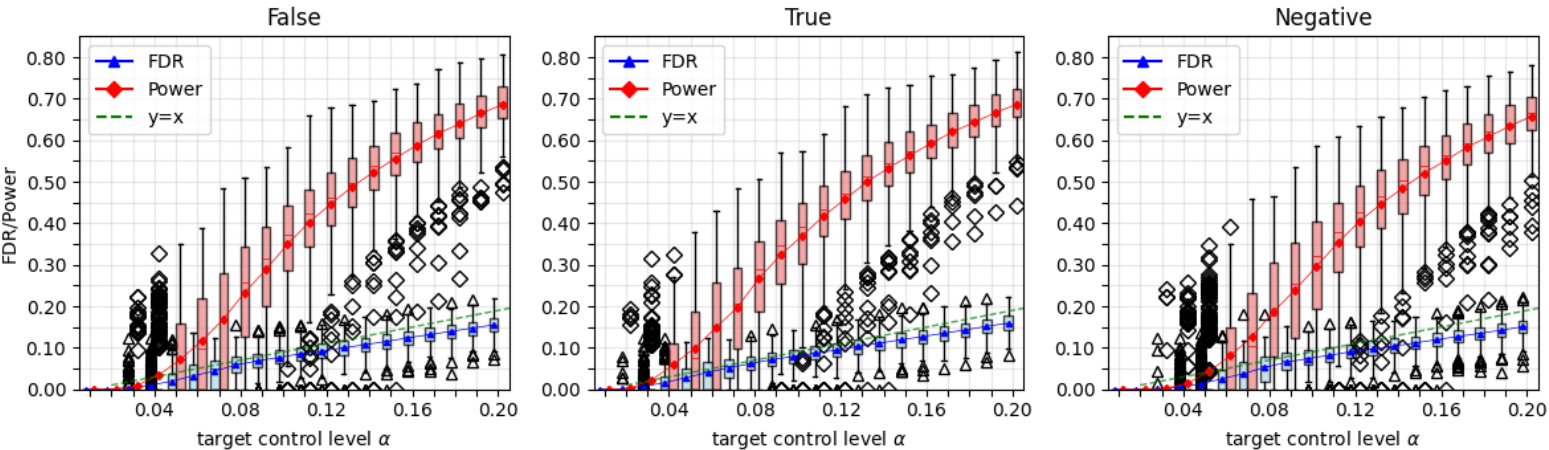}
\caption{Integration test of randomized BH procedure and competition procedure.}
\label{fig:rp_cp_integration_box}
\end{figure}

For Algorithm~\ref{alg:cp-cp-integration}, we construct $(\bW^1,\bL^1)$ with $m^1=1500$, $m^1_0=1000$, $\beta^1=2.4$, $\boldsymbol{w}^1=(0.6,1,1.4)$ and $(\bW^2,\bL^2)$ with $m^2=1000$, $m^2_0=700$, $\beta^2=2.4$, $\boldsymbol{w}^2=(0.6,1,1.4)$. Over $300$ replications, the FDR and power are estimated by the average FDP and TDP (Fig.~\ref{fig:cp_cp_integration_box}).
For Algorithm~\ref{alg:cp-p-integration}, we construct $\bX^2\in\RR_*^{m^2}$ with $m^2=2000$, $m^2_0=1700$ and $\beta^2=2.0$, together with a calibration set $\bX^{2,\rcal}\in\RR_*^{m^{2,\rcal}}$ of size $m^{2,\rcal}=1000$ and weights $\boldsymbol{w}^2=(0.6,1,1.4)$. The competition statistics $(\bW,\bL)$ have $m^1=1000$, $m^1_0=700$, $\beta^1=2.8$ and $\boldsymbol{w}^1=(0.6,1,1.4)$. Over $300$ replications, FDR and power are again estimated by the average FDP and TDP (Fig.~\ref{fig:rp_cp_integration_box}).
In Figs.~\ref{fig:cp_cp_integration_box} and~\ref{fig:rp_cp_integration_box}, `False' indicates the use of $R_+(0)$ as the weight denominator, for which FDR control is not guaranteed; `True' corresponds to executing Algorithms~\ref{alg:cp-cp-integration} and~\ref{alg:cp-p-integration} as written; and `Negative' indicates using $R_-(0)$ as the denominator. Further details are provided in the Supplementary Material.

\subsection{Comparison with e-value-based integration}

Transforming test statistics into randomized p-values provides a unified perspective on various testing methods and converting statistics into e-values offers an alternative approach. Two recent constructions turn multiple testing procedures into e-values: one converts any procedure that controls the FDR \citep[Theorem 4.2]{ignatiadis2025asymptoticcompoundevaluesmultiple}, and the other any competition procedure \citep[Theorem 1]{ren2023}. We recall them in the Supplementary Material.

Let $\cD$ be a multiple testing procedure. When $\cD$ is a competition procedure with parameter $r$, its threshold is $T_{\cD}$. Consider the e-values
\[
E^{\mathrm{mht}}_j = \frac{|\cH|\ind\{j\in\cD\}}{\alpha_{\cD}\sum_{k\in\cH}\ind\{k\in\cD\}},\quad E^{\cp}_j = \frac{|\cH|}{r}\frac{\ind\{W_j\geq T_{\cD},L_j=1\}}{1+R_-(T_{\cD})},\quad j\in\cH,
\]
where $\alpha_{\cD}$ is the level that $\cD$ controls the FDR. Applying the eBH procedure \citep{wang2022} to these two sets of e-values yields the methods MHTE (Multiple Hypothesis Testing E-values) and CPE (ComPetition E-values), respectively, whereas our randomized p-value approach is denoted by CPCF (ComPetition-ConFormal).
Using the same data generation as before, we compute the competition statistics $(\bW^1,\bL^1)\in\RR_*^{m^1}\times\{0,1\}^{m^1}$ and $(\bW^2,\bL^2)\in\RR_*^{m^2}\times\{0,1\}^{m^2}$, with $m^1=1500$, $m^1_0=1000$, $m^2=1000$, $m^2_0=800$, $\beta^1=\beta^2=2.6$ and $\boldsymbol{w}^1=\boldsymbol{w}^2=(0.6,1,1.4)$. Following the level recommended in these works, we set $\alpha_{\cD}=\alpha/(1+\alpha)$. Over $300$ replications, the FDR and power are estimated by the average FDP and TDP.

\begin{figure}[htbp]
\centering
\includegraphics[width=333.67pt]{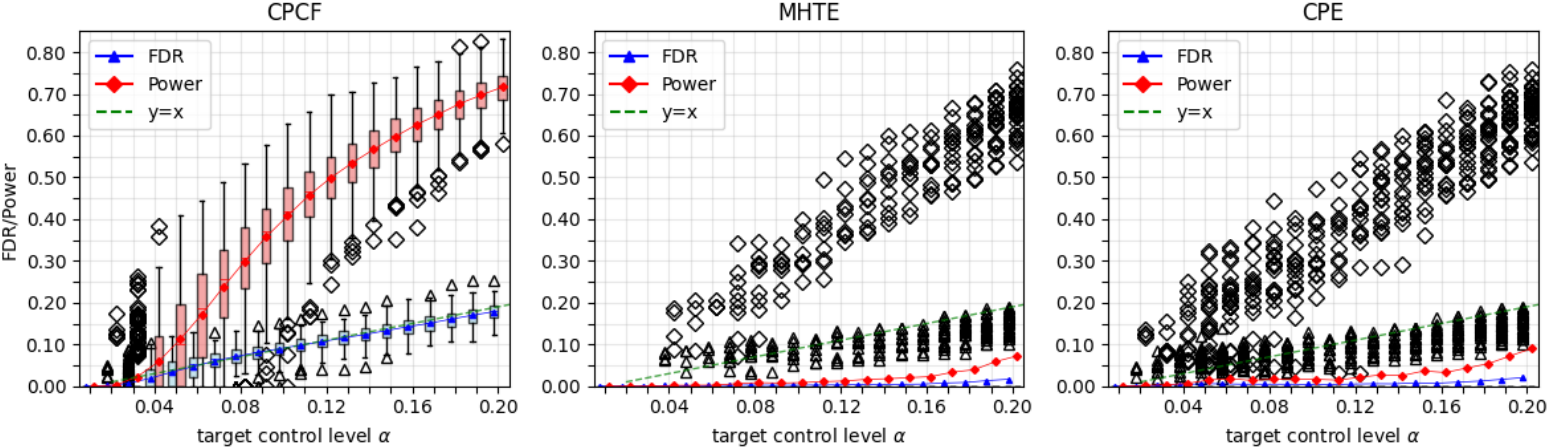}
\caption{Integration test and e-value-based integration.}
\label{fig:cp_integration_vs_e_values_box}
\end{figure}

The boxplots in Fig.~\ref{fig:cp_integration_vs_e_values_box} show that the e-value transformation sometimes matches the randomized p-value approach, but frequently makes few rejections at all, suggesting that e-value transformations may lack practical effectiveness in such settings.

The issue is that e-values constructed from multiple testing procedures are inherently tight, being either zero or only marginally above the threshold of the original procedure. When procedures are integrated, statistics from different problems perturb the combined threshold, so hypotheses whose e-values are barely above their individual thresholds may no longer be rejected. A potential remedy is to reduce the target FDR level of the procedure generating the e-values, which forces more e-values to zero while enlarging the non-zero ones \citep{ren2023}. However, reducing the level too much makes most e-values vanish and no rejections possible, while an insufficient reduction leaves the tightness problem. The stability analysis of \citet{ren2023} guides the choice of the reduced level for a single testing problem, but no universally good choice exists for heterogeneous merging problems.

\section{An application to single-cell gene-expression data}
\label{sec:data-application}

We illustrate the randomized BH framework on a real single-cell RNA-sequencing data set, reconstructing the conditional-dependence graph among a set of genes. We consider the CD4$^+$/CD45RO$^+$ memory T-cell data set of \citet{zheng2017}, distributed by 10x Genomics (Cell Ranger 1.1.0, hg19 reference). It has $32{,}738$ genes and $10{,}224$ cells. We follow a standard pre-processing pipeline: retain the $2{,}000$ genes with the largest total counts, transform the counts by $\log_2(\text{count}+1)$, keep the $50$ genes with the largest transformed variance and standardize each column. Denote by $\bX\in\RR^{10224\times 50}$ the resulting expression matrix. The object of inference is the conditional-dependence graph over $p=50$ genes, that is, the set of $1{,}225$ undirected edges $\{i,j\}$ representing a nonzero partial association between genes $i$ and $j$.

We compare the CPCF method with six existing methods for Gaussian graphical model inference, each run at nominal FDR level $q=0.2$:
\begin{enumerate}
  \item \emph{BH.} We obtain p-values from two-sided tests for zero partial correlations $\rho_{ij\mid[p]\setminus\{i,j\}}$, $i=1,\ldots,p$, $j<i$, and apply the BH procedure.
  \item \emph{BY.} The p-values are obtained as above, and we apply the BY procedure.
  \item \emph{GFC-L and GFC-SL.} Two methods for high-dimensional Gaussian graphical models proposed by \citet{liu2013}; we use the R package SILGGM \citep{zhang2018silggm} with the default values of the tuning parameters.
  \item \emph{KO2.} A knockoff-based method proposed by \citet{yu2019}; we use the authors' R implementation (\texttt{https://github.com/LedererLab/GGM-FDR}).
  \item \emph{GKF-Re+.} The GGM knockoff filter of \citet{10.1111/rssb.12430}.
\end{enumerate}

For CPCF, we first model the transformed matrix $\bX$ as a Gaussian graphical model, which considers $50$ regressions, each involving $49$ variables. Applying the knockoff construction \citep{barber2015controlling} to each regression yields the competition statistics. Because $(i,j)$ and $(j,i)$ denote the same edge, we assign each unordered pair to exactly one orientation using a balanced round-robin (tournament) scheme fixed independently of $\bX$, partitioning the $1{,}225$ pairs into $50$ groups of sizes $24$ or $25$. Embedding the competition statistics of each group into randomized p-values at level $\alpha=0.2$, as in Algorithm~\ref{alg:cp-cp-integration} extended from $2$ to $50$ groups, yields the CPCF result.

\begin{figure}[htbp]
\centering
\includegraphics[width=321.1pt]{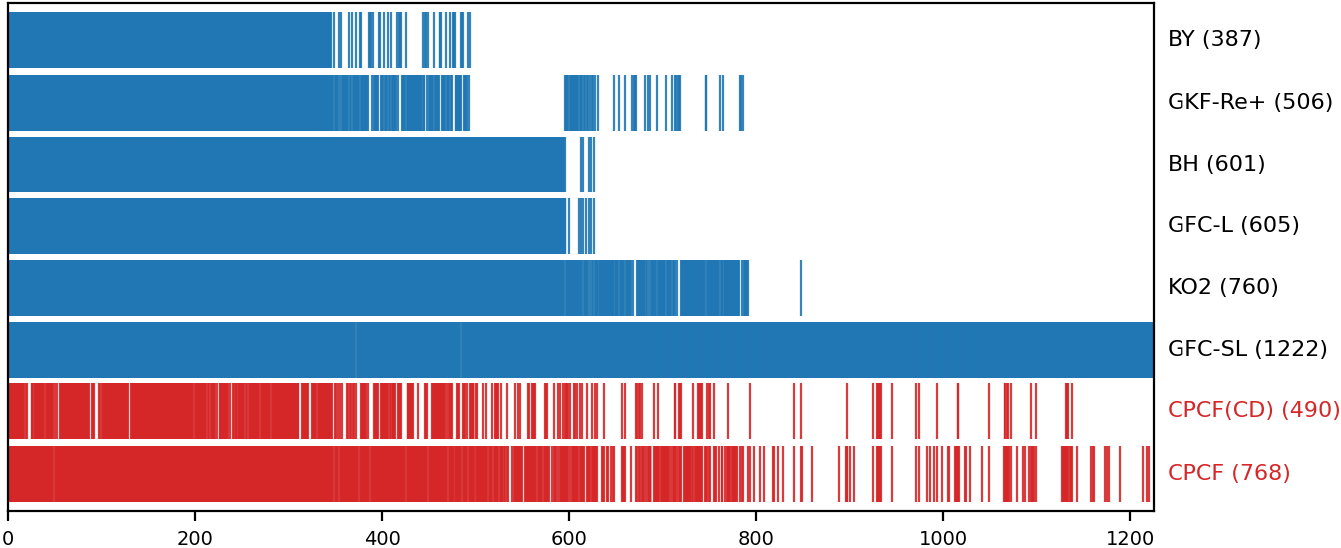}
\caption{Rejected edges of the six comparison methods, CPCF(CD), and CPCF on the memory T-cell data. Only BY, GKF-Re+, and CPCF(CD) control the FDR.}
\label{fig:memory_t}
\end{figure}

Because all groups are built from the same expression matrix, the randomized p-values of different groups are dependent. The variant CPCF(CD) (ComPetition-ConFormal with Correction for Dependence) uses the corrected level $0.105$, and Theorem~\ref{thm:cp_cp_any_dep} then guarantees FDR control at level $\alpha$ under arbitrary inter-group dependence, once each test is corrected by $c(r_g)$. Figure~\ref{fig:memory_t} displays the rejected edges of all eight procedures. 

\section{Discussion}
\label{sec:discussion}

Our randomized framework resolves the complementary limitations of the two families, freeing conformal tests from strict exchangeability and giving competition methods the scale and analytical tractability of p-values, with both resting on a common foundation in drawing evidence from the remaining hypotheses. The framework also explains why FDR control fails for compound p-values under heterogeneous marginal distributions and yields procedures that integrate heterogeneous tests with global FDR control. The simulations and the real-data application demonstrate the practical value of the framework.

Several directions for future work remain open. First, the framework may be used to construct new nonparametric multiple testing procedures, as indicated by Remark~\ref{rmk:RBH-guideline}. Second, the two applications presented here illustrate its theoretical potential; a natural next step is to investigate which deeper results from classical p-value and BH theory carry over to randomized p-values, such as optimal weighting and the structure-adaptive BH procedure \citep{li2019}. Third, the randomization introduces negative dependence, which for independent compound p-values bounds the FDR by a factor of $1.93$ and weakens as the marginal distributions become more similar; how this effect interacts with a pre-existing dependence structure, and how to derive bounds accounting for both sources, remain unresolved. Fourth, although e-values currently lie outside the framework, incorporating them would unify all null-distribution-free approaches under a single structure; whether a randomization scheme adapted to e-values suffices, or a different construction is needed, is not yet clear.




\clearpage

\part*{Supplementary Material}

\setcounter{section}{0}
\setcounter{equation}{0}
\setcounter{Thm}{0}
\setcounter{Lemma}{0}
\setcounter{Coro}{0}
\setcounter{Pro}{0}
\setcounter{Def}{0}
\setcounter{example}{0}
\setcounter{remark}{0}
\setcounter{figure}{0}
\setcounter{table}{0}
\setcounter{algocf}{0}
\def\theequation{S\arabic{section}.\arabic{equation}}
\def\thesection{S\arabic{section}}
\renewcommand{\theThm}{S\arabic{Thm}}
\renewcommand{\theLemma}{S\arabic{Lemma}}
\renewcommand{\theDef}{S\arabic{Def}}
\renewcommand{\theCoro}{S\arabic{Coro}}
\renewcommand{\thefigure}{S\arabic{figure}}
\renewcommand{\thetable}{S\arabic{table}}
\renewcommand{\thealgocf}{S\arabic{algocf}}
\renewcommand{\theHsection}{S\arabic{section}}
\def\theHThm{S\arabic{Thm}}
\def\theHLemma{S\arabic{Lemma}}
\def\theHDef{S\arabic{Def}}
\def\theHCoro{S\arabic{Coro}}
\def\theHPro{S\arabic{Pro}}
\def\theHremark{S\arabic{remark}}
\def\theHexample{S\arabic{example}}

This supplementary material is organized as follows. Section~\ref{asec:bh} proves the FDR control of the BH, eBH, and conformal BH procedures. Section~\ref{asec:competition} provides three distinct proofs for the FDR control of the competition procedure and bridges sub-competition statistics with competition statistics. Section~\ref{app:ef} proves the equivalent forms of the procedures. Section~\ref{app:gbh} establishes the FDR control of the generalized BH and randomized BH procedures. Section~\ref{asec:embedding} proves the embedding theorems and the embedding of sub-competition statistics. Section~\ref{app:compound} analyses compound p-values under the randomized BH framework, covering the violation of PRDS, negative dependence, and marginal heterogeneity. Section~\ref{asec:integration} analyses integrated tests, including a counterexample to FDR control under integration, the proofs of the integration tests, the inflation of FDR, and integrated tests with different weights. Section~\ref{asec:lemmas} collects the proofs of the lemmas used in the supplement. Section~\ref{app:morefig} reports additional numerical results, including further analysis of the prescreen test, a comparison between the conformal and competition methods, the integration of more groups of hypotheses, and the influence of weights in integration.


\section{Proofs for FDR Control for BH, eBH and Conformal BH Procedures}
\label{asec:bh}

We first recall the definitions of e-values and the eBH procedure.

\begin{Def}[eBH procedure]
    Let $\bX\in\RR_*^m$ be test statistics associated with hypotheses $\cH=[m]$. The eBH procedure at level $\alpha$ rejects
    \[
    \cR^{\mathrm{eBH}}_{\alpha}(\bX) = \{i\in\cH: X_i\geq f^{m,\alpha}_K\},
    \]
    where
    \[
    K=\sup\{k\in[m]:X_{(m-k+1)}\geq f^{m,\alpha}_k\},\qquad f^{m,\alpha}_k = \frac{m}{k\alpha},
    \]
    with the convention $\sup\emptyset=0$ and $f^{m,\alpha}_0 = +\infty$.
\end{Def}

\begin{Def}[e-value]
    A random variable $X$ is an e-value if $\EE[X]\leq 1$. Random variables $\bX\in\RR_*^m$ are compound e-values with respect to $\cN\subseteq[m]$ if $\sum_{j\in\cN}\EE[X_j]\leq m$.
\end{Def}

\begin{Thm}[FDR control for eBH procedure]
    \label{athm:ebh}
    If $\bX\in\RR_*^m$ are compound e-values with respect to $\cH_0$, then for any $\alpha\in(0,1)$,
    \[
    \FDR=\EE\left[\frac{|\cR^{\mathrm{eBH}}_{\alpha}(\bX)\cap\cH_0|}{|\cR^{\mathrm{eBH}}_{\alpha}(\bX)|\vee 1}\right]\leq \alpha.
    \]
\end{Thm}

Two constructions turn multiple testing procedures into e-values; we recall them here.

\begin{Thm}[E-value construction from a multiple testing procedure]
    \label{athm:evalue-construction}
    \citet[Theorem 4.2]{ignatiadis2025asymptoticcompoundevaluesmultiple}. Let $\cD$ be any procedure that controls the FDR at a given level $\alpha_{\cD}\in(0,1]$. Define
    \[
    E_j=\frac{|\cH|\ind\{j\in\cD\}}{\alpha_{\cD}\sum_{k\in\cH}\ind\{k\in\cD\}},\quad j\in\cH.
    \]
    Then $E_j$ are compound e-values satisfying $\sum_{j\in\cH_0}\EE(E_j)\leq |\cH|$, and applying the eBH procedure at level $\alpha_{\cD}$ using this set of e-values yields the same rejection set as directly using $\cD$.
\end{Thm}

\begin{Thm}[E-value construction from a competition procedure]
    \label{athm:evalue-construction-cp}
    The following construction is due to \citet[Theorem 1]{ren2023}. Let $\cR$ be any competition procedure that controls the FDR at a given level $\alpha_{\mathrm{cp}}\in(0,1]$, with threshold $T$ and competition statistics $(W_j,L_j)$. Define
    \[
    E_j=\frac{1}{r}\frac{|\cH|\ind\{W_j\geq T,L_j=1\}}{1+\sum_{k\in\cH}\ind\{W_k\geq T,L_k=0\}}.
    \]
    Then $E_j$ are compound e-values satisfying $\sum_{j\in\cH_0}\EE(E_j)\leq |\cH|$, and applying the eBH procedure at level $\alpha_{\mathrm{cp}}$ using this set of e-values yields the same rejection set as directly using $\cR$.
\end{Thm}

The following results are well known in the literature and are collected here for ease of reference.

\begin{proof}[Proof of Theorem~1]
By definition, $K=|\cR^{\BH}_{\alpha}(\bX)|$. Then
\begin{align*}
    \FDR&=\EE\left[\frac{\sum_{j\in\cH_0}\ind\{X_j\leq q_K^{m,\alpha}\}}{1\vee\sum_{i\in\cH}\ind\{X_i\leq q_K^{m,\alpha}\}}\right]=\sum_{j\in\cH_0}\EE\left[\frac{\ind\{X_j\leq q_K^{m,\alpha}\}}{1\vee K}\right]\\
    &=\sum_{j\in\cH_0}\EE\left[\sum_{k=1}^m\frac{\ind\{X_j\leq q_k^{m,\alpha}\}}{k}\,\ind\{K=k\}\right]\\
    &=\sum_{j\in\cH_0}\sum_{k=1}^m\frac{1}{k}\PP\!\left\{X_j\leq q_k^{m,\alpha},\,|\cR^{\BH}_{\alpha}(\bX)|=k\right\}\\
    &\leq\sum_{j\in\cH_0}\sum_{k=1}^m\frac{\alpha}{m}\PP\!\left\{|\cR^{\BH}_{\alpha}(\bX)|=k\;\big|\; X_j\leq q_k^{m,\alpha}\right\}.
\end{align*}
The inequality in the last line holds because $\bX$ are p-variables. For any $\boldsymbol{x}'\geq\boldsymbol{x}$ and $k>|\cR^{\BH}_{\alpha}(\boldsymbol{x})|$, we have $x'_{(k)}\geq x_{(k)}>q_k^{m,\alpha}$. Hence the set $\{\boldsymbol{x}:|\cR^{\BH}_{\alpha}(\boldsymbol{x})|\leq k\}$ is non-decreasing. Since $\bX$ satisfy $\PRDS_{\cH_0}$, it follows that for any $j\in\cH_0$,
\[
\PP\!\left\{|\cR^{\BH}_{\alpha}(\bX)|\leq k-1\mid X_j\leq q_{k-1}^{m,\alpha}\right\}\leq\PP\!\left\{|\cR^{\BH}_{\alpha}(\bX)|\leq k-1\mid X_j\leq q_k^{m,\alpha}\right\},
\]
and therefore
\begin{align*}
    \FDR&\leq\sum_{j\in\cH_0}\sum_{k=1}^m\frac{\alpha}{m}\PP\!\left\{|\cR^{\BH}_{\alpha}(\bX)|=k\mid X_j\leq q_k^{m,\alpha}\right\}\\
    &\leq\frac{\alpha}{m}\sum_{j\in\cH_0}\sum_{k=1}^m\Bigl(\PP\!\left\{|\cR^{\BH}_{\alpha}(\bX)|\leq k\mid X_j\leq q_k^{m,\alpha}\right\}\\
    &\qquad\qquad\qquad\qquad\quad-\PP\!\left\{|\cR^{\BH}_{\alpha}(\bX)|\leq k-1\mid X_j\leq q_{k-1}^{m,\alpha}\right\}\Bigr)\\
    &\leq\frac{\alpha}{m}\sum_{j\in\cH_0}\PP\!\left\{|\cR^{\BH}_{\alpha}(\bX)|\leq m\mid X_j\leq q_m^{m,\alpha}\right\}\leq\alpha\frac{m_0}{m}.
\end{align*}
\end{proof}

\begin{proof}[Proof of Theorem~2]
We verify the two conditions of BH-validity. First, under exchangeability, let $\bZ_{\cup j}=(X_j,\bY)$ and $\bZ_{\cup j,(\cdot)}$ be its order statistics. Then
\[
\PP\{P_j\leq t\}=\EE\bigl[\PP\{P_j\leq t\mid \bZ_{\cup j,(\cdot)},\bX_{-j}\}\bigr]=\frac{\lceil t(1+n)\rceil}{1+n}\leq t.
\]
Moreover, for any $k\in[n+1]$,
\[
\PP\!\left\{P_j=\frac{k}{1+n}\right\}=\EE\!\left[\PP\!\left\{P_j=\frac{k}{1+n}\mid\bZ_{\cup j,(\cdot)},\bX_{-j}\right\}\right]=\frac{1}{1+n}.
\]

Second, for any non-decreasing set $C\subseteq\RR^m$ and any $k\in[n+1]$,
\begin{align*}
    &\PP\{\bP_{-j}\in C\mid P_j=\tfrac{k}{1+n}\}\\
    =&\EE_{\bZ_{\cup j,(\cdot)},\bX_{-j}}\!\left[\PP\{\bP_{-j}\in C\mid P_j=\tfrac{k}{1+n},\bZ_{\cup j,(\cdot)},\bX_{-j}\}\;\Big|\;P_j=\tfrac{k}{1+n}\right],
\end{align*}
where
\begin{align*}
\PP\{\bZ_{\cup j,(\cdot)},\bX_{-j}\mid P_j=\tfrac{k}{1+n}\}
&=\frac{\PP(\bZ_{\cup j,(\cdot)},\bX_{-j})\,\PP\{P_j=\tfrac{k}{1+n}\mid\bZ_{\cup j,(\cdot)},\bX_{-j}\}}{\PP\{P_j=\tfrac{k}{1+n}\}}\\
&=\PP(\bZ_{\cup j,(\cdot)},\bX_{-j}).
\end{align*}
For any vectors $\boldsymbol{x}^j\in\RR^{m-1}$, $\boldsymbol{z}^j\in\RR^n$ and any $\boldsymbol{x}\in\RR^m$, $\boldsymbol{y}\in\RR^n$ satisfying $\boldsymbol{x}_{-j}=\boldsymbol{x}^j$, $\boldsymbol{y}_{(\cdot)}=\boldsymbol{z}^j$, define $f_j(\boldsymbol{x}^j,\boldsymbol{z}^j)=\ind\{\bP_{-j}(\boldsymbol{x},\boldsymbol{y})\in C\}$. This function is well-defined by construction of the conformal p-values. The function $f_j$ is monotone non-increasing in $\boldsymbol{z}^j$, and hence for any $w\leq k$,
\[
\PP\{\bP_{-j}\in C\mid P_j=\tfrac{w}{1+n},\bZ_{\cup j,(\cdot)},\bX_{-j}\}\leq\PP\{\bP_{-j}\in C\mid P_j=\tfrac{k}{1+n},\bZ_{\cup j,(\cdot)},\bX_{-j}\}.
\]
Taking expectations with respect to the conditional distribution derived above shows that $\PP\{\bP_{-j}\in C\mid P_j=k/(1+n)\}$ is non-decreasing in $k$. For any $w\leq k$,
\begin{align*}
    \PP\{\bP_{-j}\in C\mid P_j\leq\tfrac{k}{1+n}\}
    &=\frac{\sum_{v=1}^k\PP\{\bP_{-j}\in C\mid P_j=\frac{v}{1+n}\}\,\PP\{P_j=\frac{v}{1+n}\}}{\PP\{P_j\leq\frac{k}{1+n}\}}\\
    &=\frac{1}{k}\sum_{v=1}^k\PP\{\bP_{-j}\in C\mid P_j=\tfrac{v}{1+n}\}\\
    &\geq\frac{1}{w}\sum_{v=1}^w\PP\{\bP_{-j}\in C\mid P_j=\tfrac{v}{1+n}\}
    =\PP\{\bP_{-j}\in C\mid P_j\leq\tfrac{w}{1+n}\}.
\end{align*}
Thus the conformal p-values $\bP(\bX,\bY)$ satisfy $\PRDS_{\cH_0}$ and are BH-valid.
\end{proof}

\begin{proof}[Proof of Theorem~\ref{athm:ebh}]
By definition, $K=|\cR^{\mathrm{eBH}}_{\alpha}(\bX)|$. When $K\leq m$, the identity $f^{m,\alpha}_K=m/[(m-K+1)\alpha]$ gives $\alpha f^{m,\alpha}_K/m=1/(m-K+1)=1/(1\vee(m-K+1))$. When $K=m+1$, both $1/(1\vee(m-K+1))$ and $\alpha f^{m,\alpha}_{K\wedge m}/m$ equal $1$. Hence
\begin{align*}
    \FDR&=\EE\left[\frac{\sum_{j\in\cH_0}\ind\{X_j\geq f_K^{m,\alpha}\}}{1\vee\sum_{i\in\cH}\ind\{X_i\geq f_K^{m,\alpha}\}}\right]=\sum_{j\in\cH_0}\EE\left[\frac{\ind\{X_j\geq f_K^{m,\alpha}\}}{1\vee(m-K+1)}\right]\\
    &=\sum_{j\in\cH_0}\EE\left[\frac{\alpha f^{m,\alpha}_{K\wedge m}}{m}\,\ind\{X_j\geq f_K^{m,\alpha}\}\right]\\
    &\leq\frac{\alpha}{m}\sum_{j\in\cH_0}\EE[X_j\,\ind\{X_j\geq f_K^{m,\alpha}\}]\leq\frac{\alpha}{m}\sum_{j\in\cH_0}\EE[X_j].
\end{align*}
Since $\bX$ are compound e-values, $\FDR\leq\alpha m^{-1}\sum_{j\in\cH_0}\EE[X_j]\leq\alpha$.
\end{proof}

\section{Proofs for the Competition Procedure}
\label{asec:competition}

In this section we prove Theorem~3. We provide three distinct proofs, referred to as the martingale, combination, and flip-one proofs. The first two assume that $(\bW,\bL)$ are competition statistics, so that the label probabilities are identical fixed values. The third requires only the weaker sub-competition condition introduced in the main paper.

We begin with a common derivation. Regardless of whether $(\bW,\bL)$ are sub-competition statistics, the threshold $T$ satisfies
\begin{align*}
    \FDR&=\EE\left[\frac{\sum_{j\in\cH_0}\ind\{W_j\geq T,L_j=1\}}{1\vee\sum_{j\in\cH}\ind\{W_j\geq T,L_j=1\}}\right]\\
    &=\EE\left[\frac{r\sum_{j\in\cH}\ind\{W_j\geq T,L_j=0\}+r}{1\vee\sum_{j\in\cH}\ind\{W_j\geq T,L_j=1\}}\cdot\frac{\sum_{j\in\cH_0}\ind\{W_j\geq T,L_j=1\}}{r\sum_{j\in\cH}\ind\{W_j\geq T,L_j=0\}+r}\right]\\
    &\leq\alpha\EE\left[\frac{\sum_{j\in\cH_0}\ind\{W_j\geq T,L_j=1\}}{r\sum_{j\in\cH}\ind\{W_j\geq T,L_j=0\}+r}\right].
\end{align*}

Thus FDR control follows from the sufficient condition
\[
\EE\left[\frac{\sum_{j\in\cH_0}\ind\{W_j\geq T,L_j=1\}}{r\sum_{j\in\cH}\ind\{W_j\geq T,L_j=0\}+r}\right]\leq1,
\]
which is the focus of the three proofs below.

\subsection{Three Distinct Proofs}

\begin{proof}[Proof via martingale]

Assume that $(\bW,\bL)$ are competition statistics. Let $\bW_{(\cdot)}$ be non-increasing with $W_{(0)}=+\infty$, and recall the definitions of $V_\pm(t)$ and $\cR_\pm(t)$. Define
\[
M_t=\frac{V_+(t)}{1+V_-(t)},\qquad \widetilde{M}_k=M_{W_{(k)}}.
\]
Let $\cH_0(t)=\{j\in\cH_0:W_j\geq t\}$ and $k_t=\arg\min_j\{W_j:j\in\cH_0(t)\}$. When $V_+(t)>0$, the competition label of $k_t$ satisfies
\[
\PP\left(L_{k_t}=1\;\big|\;V_+(t),V_-(t),\cH_0(t)\right)=\frac{V_+(t)}{V_+(t)+V_-(t)}.
\]
Hence, for any $k$,
\begin{align*}
    &\EE\left[\widetilde{M}_{k}-\widetilde{M}_{k+1}\Big|\sigma\{\widetilde{M}_{k}\}\right]
    =\EE\left[\frac{V_+(W_{(k)})}{1+V_-(W_{(k)})}-\frac{V_+(W_{(k+1)})}{1+V_-(W_{(k+1)})}\Big|\sigma\{\widetilde{M}_{k}\}\right]\\
    =&\EE\Bigg[\EE\left[\frac{V_+(W_{(k)})}{1+V_-(W_{(k)})}-\frac{V_+(W_{(k+1)})}{1+V_-(W_{(k+1)})}\;\Big|\;V_+(W_{(k)}),V_-(W_{(k)}),\cH_0(W_{(k)})\right]\Bigg|\;\sigma\{\widetilde{M}_{k}\}\Bigg]\\
    \leq& 0.
\end{align*}
Thus $\{\widetilde{M}_k\}$ is a reverse supermartingale, and Doob's optional stopping theorem yields
\begin{align*}
    &\EE\left[\frac{\sum_{j\in\cH_0}\ind\{W_j\geq T,L_j=1\}}{r\sum_{j\in\cH}\ind\{W_j\geq T,L_j=0\}+r}\right]\\
    &\leq\EE\left[\frac{\sum_{j\in\cH_0}\ind\{W_j\geq 0,L_j=1\}}{r\sum_{j\in\cH_0}\ind\{W_j\geq 0,L_j=0\}+r}\right]\\
    &=\sum_{k=0}^{|\cH_0|}\frac{1}{r}\frac{k}{1+|\cH_0|-k}\binom{|\cH_0|}{k}\left(\frac{r}{1+r}\right)^{k}\left(\frac{1}{1+r}\right)^{|\cH_0|-k}<1.
\end{align*}
\end{proof}

\begin{proof}[Proof via combination]

Assume that $(\bW,\bL)$ are competition statistics. Construct a dominating variable in the sense of expectation,
\[
\delta(t)=\max_{u\geq t}\left\{V_+(u):V_-(u)=0\right\}.
\]
Clearly $\delta(0)\geq\delta(T)$. We also need the following lemma.

\begin{Lemma}
\label{slem:overvariable}
Given $V_+(T)$ and $V_-(T)$, the conditional distribution of $\delta(T)$ is
\[
\PP\left(\delta(T)=k\;\Big|\;V_+(T),V_-(T)\right)=f\left(k;V_+(T),V_-(T)\right),
\]
which does not depend on $T$, and its conditional expectation is
\[
\EE\left[\delta(T)\;\big|\;V_+(T),V_-(T)\right]=\frac{V_+(T)}{1+V_-(T)}.
\]
\end{Lemma}

Since $\delta(0)$ follows a truncated geometric distribution, we obtain
\[
\EE\left[\frac{\sum_{j\in\cH_0}\ind\{W_j\geq T,L_j=1\}}{r\sum_{j\in\cH}\ind\{W_j\geq T,L_j=0\}+r}\right]=\EE[\delta(T)]/r\leq\EE[\delta(0)]/r<1.
\]
\end{proof}

\begin{proof}[Proof via flip-one]

Assume that $(\bW,\bL)$ are competition statistics. Let
\[
T_i=\inf\left\{t\in\{W_j:L_j=1\}:\frac{r\sum_{j\neq i}\ind\{W_j\geq t,L_j=0\}+r}{1\vee\sum_{j\neq i}\ind\{W_j\geq t,L_j=0\}+\ind\{W_i\geq t\}}\leq\alpha\right\}.
\]
If $\ind\{W_i\geq T,L_i=1\}=1$, then $T=T_i$. Consequently,
\[
\EE\left[\frac{\ind\{W_i\geq T,L_i=1\}}{r\sum_{j\in\cH}\ind\{W_j\geq T,L_j=0\}+r}\right]=\EE\left[\frac{\ind\{W_i\geq T_i,L_i=1\}}{r\sum_{j\in\cH\backslash\{i\}}\ind\{W_j\geq T_i,L_j=0\}+r}\right].
\]
Since $\sigma\{T_i\}\subseteq\sigma\{\bW,\bL_{-i}\}$, we have
\begin{align*}
    &\EE\left[\frac{\ind\{W_i\geq T_i,L_i=1\}}{r\sum_{j\in\cH\backslash\{i\}}\ind\{W_j\geq T_i,L_j=0\}+r}\right]\\
    =&\EE\left[\EE\left[\frac{\ind\{W_i\geq T_i,L_i=1\}}{r\sum_{j\in\cH\backslash\{i\}}\ind\{W_j\geq T_i,L_j=0\}+r}\Big|\sigma\{\bW,\bL_{-i}\}\right]\right]\\
    =&\EE\left[\PP(L_i=1\mid\bW,\bL_{-i})\cdot\frac{\ind\{W_i\geq T_i\}}{r\sum_{j\in\cH\backslash\{i\}}\ind\{W_j\geq T_i,L_j=0\}+r}\right]\\
    =&\frac{r}{1+r}\left(\EE\left[\frac{\ind\{W_i\geq T_i,L_i=1\}}{r\sum_{j\in\cH\backslash\{i\}}\ind\{W_j\geq T_i,L_j=0\}+r}\right]\right.\\
    &\qquad\qquad\left.+\EE\left[\frac{\ind\{W_i\geq T_i,L_i=0\}}{r\sum_{j\in\cH\backslash\{i\}}\ind\{W_j\geq T_i,L_j=0\}+r}\right]\right).
\end{align*}
Hence
\[
\EE\left[\frac{\ind\{W_i\geq T,L_i=1\}}{r\sum_{j\in\cH\backslash\{i\}}\ind\{W_j\geq T,L_j=0\}+r}\right]\leq\EE\left[\frac{\ind\{W_i\geq T_i,L_i=0\}}{\sum_{j\in\cH\backslash\{i\}}\ind\{W_j\geq T_i,L_j=0\}+1}\right].
\]
A further lemma shows that $T_i$ takes at most two distinct values across different indices $i$.

\begin{Lemma}
\label{slem:biovervariable}
If $\ind\{W_i\geq T,L_i=0\}=1$, then the following hold:
\begin{enumerate}
    \item $T_i\leq T$;
    \item $\ind\{W_k\geq T_i,L_k=0\}=1\Rightarrow T_i=T_k$;
    \item $\ind\{W_k\geq T_i,L_k=0\}=0\Rightarrow T=T_k$.
\end{enumerate}
\end{Lemma}

By Lemma~\ref{slem:biovervariable},
\[
\EE\left[\frac{\sum_{j\in\cH_0}\ind\{W_j\geq T,L_j=1\}}{r\sum_{j\in\cH}\ind\{W_j\geq T,L_j=0\}+r}\right]=\sum_{i\in\cH_0}\EE\left[\frac{\ind\{W_i\geq T_i,L_i=0\}}{\sum_{j\in\cH\backslash\{i\}}\ind\{W_j\geq T_i,L_j=0\}+1}\right]\leq1.
\]
\end{proof}

\subsection{Comparison of the Three Proofs}

We have presented three proofs of FDR control for the competition procedure. They differ in the assumptions they require and in the mechanisms they exploit, and each offers distinct advantages. In this subsection we compare the three proofs from these perspectives.

The martingale proof is the classical and earliest approach. It assumes competition statistics, under which $\PP\{L_j=1\mid\bW,\bL_{-j}\}=r/(1+r)$, and constructs the process $\{\widetilde{M}_k\}$, shown to be a reverse supermartingale. The FDR bound then follows from Doob's optional stopping theorem. Its main advantage is that it identifies a sufficient condition more general than the competition-statistics framework: FDR control holds whenever the process forms a reverse supermartingale with a corresponding stopping time.

The combination proof also assumes competition statistics. Its idea is to construct an auxiliary random variable $\delta(t)$ that dominates the quantity of interest in expectation. The fraction $V_+(T)/(rR_-(T)+r)$, whose distribution is hard to obtain explicitly, is dominated by $\delta(0)/r$, a random variable with a known truncated geometric distribution. Specifically, $\delta(T)$ counts the false discoveries that occur before any negative label at threshold $t$, and its conditional expectation equals $V_+(T)/(1+V_-(T))$. Since $\delta(0)$ follows a truncated geometric distribution, whose expectation is easy to bound, this dominance replaces an intractable distributional calculation with a simple bound.

The flip-one proof is the most flexible. It works directly under sub-competition statistics, where $\PP\{L_j=1\mid\bW,\bL_{-j}\}\leq r/(1+r)$, and requires no auxiliary construction. Its strategy is a finer analysis of the threshold. After flipping the label of a single null hypothesis to $0$, the threshold $T$ changes to $T_i$. Lemma~\ref{slem:biovervariable} guarantees that $T_i$ takes at most two values and that its post-flip behaviour is fully determined. This allows the sub-competition condition to bound each null hypothesis individually, and it extends naturally to general dependence structures and asymptotic analysis.

\subsection{The Bridge between Sub-competition Statistics and Competition Statistics}
\label{sec:bridge}

First, the flip-one proof extends directly to general sub-competition statistics, without the need for any bridge construction. The martingale and combination proofs, by contrast, require the following bridge.

Suppose $(\bW,\bL)$ are sub-competition statistics. Without loss of generality, assume $W_1<\dots<W_m$. Construct a sequence of random variables $Y_1,Y_2,\dots,Y_m$ sequentially, with conditional distributions satisfying for all $j\in\cH_0$,
\[
\PP\{Y_j=1\mid \bL_{\{i:i\leq j\}},\bY_{\{i:i<j\}},\bL_{\cH_1},\bW\}=1-L_j-\frac{1}{1+r}(1-L_j)\bigl(1-p_j(\bL_{\{i:i<j\}})\bigr)^{-1}
\]
and $Y_j=0$ for all $j\in\cH_1$, where
\[
p_j(\bL_{\{i:i<j\}})=\PP\{L_j=1\mid \bL_{\{i:i<j\}},\bL_{\cH_1},\bW\}.
\]
This construction ensures that $L_j$ is independent of $\bY_{\{i:i<j\}}$ given $\bL_{\{i:i<j\}}$ for all $j\in\cH_0$, and that
\[
\PP\left\{L_j+Y_j=1\mid\bY_{\{i:i<j\}},\bL_{\{i:i<j\}},\bL_{\cH_1},\bW\right\}=\frac{r}{1+r},\quad \forall j\in\cH_0.
\]
Furthermore, for any $k>j$ with $j,k\in\cH_0$,
\begin{align*}
&\PP\left\{L_{k}+Y_{k}=1\mid\bY_{\{i:i<j\}},\bL_{\{i:i<j\}},\bL_{\cH_1},\bW\right\}\\
=&\EE\left[\PP\left\{L_{k}+Y_{k}=1\mid\bY_{\{i:i<k\}},\bL_{\{i:i<k\}},\bL_{\cH_1},\bW\right\}\mid\bY_{\{i:i<j\}},\bL_{\{i:i<j\}},\bL_{\cH_1},\bW\right]\\
=&\EE\left[\frac{r}{1+r}\mid\bY_{\{i:i<j\}},\bL_{\{i:i<j\}},\bL_{\cH_1},\bW\right]=\frac{r}{1+r}.
\end{align*}
Setting $Z_j=L_j+Y_j$, we obtain
\[
\PP\{Z_{j_1}=1\mid \bZ_{\{i:i<j_2\}},\bL_{\cH_1},\bW\}=\frac{r}{1+r},\quad\forall j_1\in\cH_0,\,j_2\leq j_1,
\]
and, for any $k>j$ with $j,k\in\cH_0$,
\begin{align*}
    &\PP\{Z_j=1\mid \bZ_{\{i:i<j\}},Z_k,\bL_{\cH_1},\bW\}\\
    =&\frac{\PP\{Z_k\mid \bZ_{\{i:i\leq j\}},\bL_{\cH_1},\bW\}\,\PP\{Z_j\mid \bZ_{\{i:i<j\}},\bL_{\cH_1},\bW\}}{\PP\{Z_k\mid \bZ_{\{i:i<j\}},\bL_{\cH_1},\bW\}}\\
    =&\PP\{Z_j\mid \bZ_{\{i:i<j\}},\bL_{\cH_1},\bW\}=\frac{r}{1+r}.
\end{align*}
Then $(\bW,\bZ)$ are competition statistics.
For $u\geq 0$, write
\[
R^{Y}_{+}(u)=\sum_{j\in\cH_0(u)}\ind\{Y_j=1\},\quad R^{Y}_{-}(u)=\sum_{j\in\cH_0(u)}\ind\{Y_j=0\}.
\]
Now consider the filtration $\{\sigma_t\}$ generated by
\[
\sigma_t=\sigma\left\{R^{Y}_{+}(u),R^{Y}_{-}(u),R_{+}(u),R_{-}(u),\bL_{\cH_1},\bW:\ u\leq t\right\}.
\]
$(\bW,\bZ)$ are competition statistics, but FDR control does not follow directly from this fact, since the rejection set of the competition procedure depends on $(\bW,\bL)$ rather than $(\bW,\bZ)$. The construction above nevertheless provides enough structure to satisfy the conditions required by the martingale proof.

\section{Proofs for the Equivalent Forms}
\label{app:ef}

\begin{proof}[Proof of the first part of Theorem~4]

If $X_j\geq T$, then by the definition of $T$,
\[
\frac{m}{1+n}\frac{1+\sum_{i=1}^n\ind\{Y_i\geq X_j\}}{1\vee\sum_{i=1}^m\ind\{X_i\geq T\}}\leq\frac{m}{1+n}\frac{1+\sum_{i=1}^n\ind\{Y_i\geq T\}}{1\vee\sum_{i=1}^m\ind\{X_i\geq T\}}\leq \alpha.
\]
Thus for any $j$ with $X_j\geq T$,
\[
P_j(\bX,\bY)=\frac{1+\sum_{i=1}^n\ind\{Y_i\geq X_j\}}{1+n}\leq\alpha\frac{\sum_{i=1}^m\ind\{X_i\geq T\}}{m},
\]
and therefore
\[
\left\{k\in[m]:X_k\geq T\right\}\subseteq\left\{k\in[m]:P_k(\bX,\bY)\leq\alpha\frac{\sum_{i=1}^m\ind\{X_i\geq T\}}{m}\right\}.
\]
Similarly, if $X_j<T$, then
\[
\frac{m}{1+n}\frac{1+\sum_{i=1}^n\ind\{Y_i\geq X_j\}}{1\vee\sum_{i=1}^m\ind\{X_i\geq T\}}\geq\frac{m}{1+n}\frac{1+\sum_{i=1}^n\ind\{Y_i\geq X_j\}}{1\vee\sum_{i=1}^m\ind\{X_i\geq X_j\}}>\alpha,
\]
so for any $j$ with $X_j<T$,
\[
P_j(\bX,\bY)=\frac{1+\sum_{i=1}^n\ind\{Y_i\geq X_j\}}{1+n}>\alpha\frac{\sum_{i=1}^m\ind\{X_i\geq X_j\}}{m},
\]
and
\[
\left\{k\in[m]:P_k(\bX,\bY)\leq\alpha\frac{\sum_{i=1}^m\ind\{X_i\geq X_j\}}{m}\right\}\subsetneq\left\{k\in[m]:X_k\geq X_j\right\}.
\]
Since the order of $P_j$ coincides with that of $X_j$,
\[
\cR^{*\cp}_{m(1+n)^{-1},\alpha}((\bX,\bY),\ind_{\cD^{\test}})=\cR^{\BH}_{\alpha}(\bP(\bX,\bY)).
\]
Moreover, since
\[
K=|\cR^{*\cp}_{m(1+n)^{-1},\alpha}((\bX,\bY),\ind_{\cD^{\test}})|=|\cR^{\BH}_{\alpha}(\bP(\bX,\bY))|,
\]
we obtain
\[
\frac{1+\sum_{i=1}^n\ind\{Y_i\geq T\}}{1+n}\leq \alpha\frac{\sum_{i=1}^m\ind\{X_i\geq T\}}{m}=\frac{\alpha K}{m}=q_K^{m,\alpha}.
\]
\end{proof}

\begin{proof}[Proof of the second part of Theorem~4]

If $R_+(0)=0$, both $\cR^{\cp}_{r,\alpha}(\bW,\bL)$ and $\cR^{*\BH}_{\alpha}(\bP(\bW,\bL);\bE(\bW,\bL))$ are empty. The two rejection sets coincide, so it suffices to consider $R_+(0)>0$. Suppose first that $P_j/E_j\leq q^{R_+(0),\alpha}_K$. Then
\[
r\frac{1+R_-(W_j)}{1\vee R_+(W_j)}=r\frac{1+R_-(W_j)}{R_+(0)}\frac{R_+(0)}{1\vee R_+(W_j)}\leq\alpha\frac{K}{R_+(0)}\frac{R_+(0)}{1\vee R_+(W_j)}=\alpha\frac{K}{1\vee R_+(W_j)}.
\]
Since the values of $E_j$ are all equal and $W_j$ and $P_j$ have the same order, the definition of $K$ gives $R_+(W_j)=|\{i:W_i\geq W_j,L_i=1\}|\leq K$ for every $j$ with $P_j/E_j\leq q^{R_+(0),\alpha}_K$. Hence there exists some $j_K$ with $R_+(W_{j_K})=K$, $P_{j_K}/E_{j_K}\leq q_K^{R_+(0),\alpha}$, and
\[
r\frac{1+R_-(W_{j_K})}{1\vee R_+(W_{j_K})}\leq \alpha\frac{K}{K}=\alpha,
\]
which implies $T\leq W_{j_K}$ and
\[
\cR^{*\BH}_{\alpha}(\bP(\bW,\bL);\bE(\bW,\bL))\subseteq\cR^{\cp}_{r,\alpha}(\bW,\bL).
\]
Conversely, for any $j$ with $P_j/E_j> q^{R_+(0),\alpha}_K$, the definition of $K$ yields $P_j/E_j>q^{R_+(0),\alpha}_{R_+(W_j)}$, hence $R_+(W_j)>K$, and
\[
r\frac{1+R_-(W_j)}{1\vee R_+(W_j)}>\alpha\frac{R_+(W_j)}{1\vee R_+(W_j)}=\alpha,
\]
so that
\[
\cR^{\cp}_{r,\alpha}(\bW,\bL)\subseteq\cR^{*\BH}_{\alpha}(\bP(\bW,\bL);\bE(\bW,\bL)).
\]
Combining the two inclusions gives
\[
\cR^{\cp}_{r,\alpha}(\bW,\bL)=\cR^{*\BH}_{\alpha}(\bP(\bW,\bL);\bE(\bW,\bL)).
\]
\end{proof}

\section{Proofs for the Generalized BH and Randomized BH Procedures}
\label{app:gbh}

\begin{proof}[Proof of Theorem~5]

By definition, $K=|\cR^{\BH}_{\alpha}(\bP')|$, where $\bP'=\bP/\bE$. Conditioning on $\bE$,
\begin{align*}
    &\EE\left[\frac{\sum_{j\in\cH_0}\ind\{P_j/E_j\leq q_K^{m,\alpha}\}}{1\vee\sum_{i\in\cH}\ind\{P_i/E_i\leq q_K^{m,\alpha}\}}\;\Big|\;\bE\right]\\
    =&\sum_{j\in\cH_0}\EE\left[\sum_{k=1}^m\frac{\ind\{P_j/E_j\leq q_k^{m,\alpha}\}}{k}\,\ind\{K=k\}\;\Big|\;\bE\right]\\
    =&\sum_{j\in\cH_0}\sum_{k=1}^m\frac{1}{k}\PP\left\{P_j/E_j\leq q_k^{m,\alpha},\,|\cR^{\BH}_{\alpha}(\bP')|=k\;\Big|\;\bE\right\}\\
    \leq&\sum_{j\in\cH_0}\sum_{k=1}^m\frac{\alpha E_j}{m}\PP\left\{|\cR^{\BH}_{\alpha}(\bP')|=k\;\Big|\;P_j\leq q_k^{m,\alpha}E_j,\,\bE\right\}.
\end{align*}
Conditioned on $\bE$, the set $\{|\cR^{\BH}_{\alpha}(\bP')|\leq k\}$ is non-decreasing in $\bP_{-j}$ for every $k$ and every $j$. Hence
\begin{align*}
    \FDR&=\EE\left[\EE\left[\FDP\mid\bE\right]\right]\\
    &\leq\EE\left[\sum_{j\in\cH_0}\sum_{k=1}^m\frac{\alpha E_j}{m}\PP\left\{|\cR^{\BH}_{\alpha}(\bP')|=k\;\Big|\;P_j\leq q_k^{m,\alpha}E_j,\,\bE\right\}\right]\\
    &\leq\frac{\alpha}{m}\sum_{j\in\cH_0}\EE[E_j].
\end{align*}
\end{proof}

\begin{proof}[Proof of Lemma~1]

By the definition of the weights, we have
\begin{align*}
    &\sum_{j\in\cH_0}\EE\left[\frac{1\vee R_+(0)}{r+rV_-(0)}\,\ind\{L_j=1\}\right]\\
    =&\sum_{j\in\cH_0}\EE\left[\frac{\sum_{i\in\cH}\ind\{L_i=1\}}{r+rV_-(0)}\,\ind\{L_j=1\}\right]\\
    =&\sum_{i\in\cH}\EE\left[\frac{\sum_{j\in\cH_0}\ind\{L_j=1\}}{r+rV_-(0)}\,\ind\{L_i=1\}\right]\\
    =&\sum_{i\in\cH}\EE\left[\frac{V_+(0)}{r+rV_-(0)}\,\ind\{L_i=1\}\right]\\
    \leq&m\,\EE\left[\frac{V_+(0)}{r+rV_-(0)}\right].
\end{align*}
The first equality holds because if $\sum_{i\in\cH}\ind\{L_i=1\}=0$, then $L_j=0$ for all $j$ and the expectation vanishes, so $1\vee R_+(0)$ may be replaced by $\sum_{i\in\cH}\ind\{L_i=1\}$. Furthermore,
\[
\EE\left[\frac{V_+(0)}{r+rV_-(0)}\right]=\sum_{k=0}^{m_0}\binom{m_0}{k}\frac{k}{r+rm_0-rk}\left(\frac{r}{1+r}\right)^k\left(\frac{1}{1+r}\right)^{m_0-k}<1,
\]
and therefore
\[
\sum_{j\in\cH_0}\EE\left[\frac{1\vee R_+(0)}{r+rR_-(0)}\,\ind\{L_j=1\}\right]\leq\sum_{j\in\cH_0}\EE\left[\frac{1\vee R_+(0)}{r+rV_-(0)}\,\ind\{L_j=1\}\right]\leq m.
\]
\end{proof}

\begin{proof}[Proof of Theorem~6]

By definition, $K=|\cR^{\BH}_{\alpha}(\bP)|$, so
\begin{align*}
    &\EE\left[\frac{\sum_{j\in\cH_0}\ind\{P_j\leq q_K^{m,\alpha}\}}{1\vee\sum_{i\in\cH}\ind\{P_i\leq q_K^{m,\alpha}\}}\right]=\sum_{j\in\cH_0}\EE\left[\sum_{k=1}^m\frac{\ind\{P_j\leq q_k^{m,\alpha}\}}{k}\,\ind\{K=k\}\right]\\
    =&\sum_{j\in\cH_0}\sum_{k=1}^m\frac{1}{k}\PP\left\{P_j\leq q_k^{m,\alpha},\,|\cR^{\BH}_{\alpha}(\bP)|=k\right\}\\
    \leq&\sum_{j\in\cH_0}\sum_{k=1}^m\frac{1}{k}\PP\left\{|\cR^{\BH}_{\alpha}(\bP)|=k,\,Q_j\leq q_k^{m,\alpha}\right\}.
\end{align*}
The third line holds because $\PP(Q_j\leq P_j)=1$. Unlike in the standard BH proof, the condition $Q_j\leq q_k^{m,\alpha}$ does not immediately yield a non-decreasing set, so we refine the event first:
\begin{align*}
    &\left\{P_j\leq q_k^{m,\alpha},\,|\cR^{\BH}_{\alpha}(\bP)|=k\right\}\\
    =&\left\{P^j_{(k-1)}\leq q^{m,\alpha}_k,\,P^j_{(k)}>q^{m,\alpha}_{k+1},\,\dots,\,P^j_{(m-1)}>q^{m,\alpha}_m,\,P_j\leq q^{m,\alpha}_m\right\},
\end{align*}
where $P^j_{(i)}$ denotes the $i$-th order statistic of $\bP$ with $P_j$ excluded. The set $\{\boldsymbol{p}\in\RR_*^{m-1}:p_{(k)}>q^{m,\alpha}_{k+1},\dots,p_{(m-1)}>q^{m,\alpha}_m\}$ is non-decreasing. Hence
\begin{align*}
    &\FDR=\sum_{j\in\cH_0}\sum_{k=1}^m\frac{1}{k}\PP\left\{|\cR^{\BH}_{\alpha}(\bP)|=k,\,P_j\leq q_k^{m,\alpha}\right\}\\
    \leq&\sum_{j\in\cH_0}\sum_{k=1}^m\frac{1}{k}\PP\left\{P^j_{(k-1)}\leq q^{m,\alpha}_k,\dots,P^j_{(m-1)}>q^{m,\alpha}_m,\,Q_j\leq q^{m,\alpha}_k\right\}\\
    \leq&\sum_{j\in\cH_0}\sum_{k=1}^m\frac{\alpha}{m}\PP\left\{P^j_{(k-1)}\leq q^{m,\alpha}_k,\dots,P^j_{(m-1)}>q^{m,\alpha}_m\;\Big|\;Q_j\leq q^{m,\alpha}_k\right\}\\
    =&\sum_{j\in\cH_0}\sum_{k=1}^m\frac{\alpha}{m}\left(\PP\left\{P^j_{(k)}>q^{m,\alpha}_{k+1},\dots,P^j_{(m-1)}>q^{m,\alpha}_m\;\Big|\;Q_j\leq q^{m,\alpha}_k\right\}\right.\\
    &\quad\left.-\PP\left\{P^j_{(k-1)}>q^{m,\alpha}_k,P^j_{(k)}>q^{m,\alpha}_{k+1},\dots,P^j_{(m-1)}>q^{m,\alpha}_m\;\Big|\;Q_j\leq q^{m,\alpha}_k\right\}\right)\\
    \leq&\sum_{j\in\cH_0}\frac{\alpha}{m}\leq\alpha\frac{m_0}{m}.
\end{align*}
\end{proof}

\begin{proof}[Proof of Theorem~7]

By definition, $K^\eta=|\cR^{\BH}_{\alpha}(\bP_\eta)|$. Conditioning on $\eta$,
\begin{align*}
    &\EE\left[\frac{\sum_{j\in\cH_0\cap\eta}\ind\{P_j\leq q_{K^\eta}^{|\eta|,\alpha}\}}{1\vee\sum_{i\in\eta}\ind\{P_i\leq q_{K^\eta}^{|\eta|,\alpha}\}}\;\Big|\;\eta\right]\\
    =&\sum_{j\in\cH_0\cap\eta}\EE\left[\sum_{k=1}^{|\eta|}\frac{\ind\{P_j\leq q_k^{|\eta|,\alpha}\}}{k}\,\ind\{K^\eta=k\}\;\Big|\;\eta\right]\\
    =&\sum_{j\in\cH_0\cap\eta}\sum_{k=1}^{|\eta|}\frac{1}{k}\PP\left\{P_j\leq q_k^{|\eta|,\alpha},\,|\cR^{\BH}_{\alpha}(\bP_\eta)|=k\;\Big|\;\eta\right\}\\
    \leq&\sum_{j\in\cH_0\cap\eta}\sum_{k=1}^{|\eta|}\frac{\alpha}{|\eta|}\PP\left\{|\cR^{\BH}_{\alpha}(\bP_\eta)|=k\;\Big|\;P_j\leq q_k^{|\eta|,\alpha},\,\eta\right\}.
\end{align*}
Conditioned on $\eta$, the set $\{|\cR^{\BH}_{\alpha}(\bP_\eta)|\leq k\}$ is non-decreasing in $\bP_{\eta,-j}$ for every $k$ and every $j\in\cH_0\cap\eta$. Hence
\begin{align*}
    \FDR&=\EE\left[\EE\left[\frac{|\cH_0\cap\eta\cap\cR^{\BH}_{\alpha}(\bP_\eta)|}{1\vee|\eta\cap\cR^{\BH}_{\alpha}(\bP_\eta)|}\;\Big|\;\eta\right]\right]\\
    &\leq\EE\left[\sum_{j\in\cH_0\cap\eta}\sum_{k=1}^{|\eta|}\frac{\alpha}{|\eta|}\PP\left\{|\cR^{\BH}_{\alpha}(\bP_\eta)|=k\;\Big|\;P_j\leq q_k^{|\eta|,\alpha},\,\eta\right\}\right]\\
    &\leq\alpha\EE\left[\frac{|\cH_0\cap\eta|}{1\vee|\eta|}\right]\leq\alpha.
\end{align*}
\end{proof}

\begin{proof}[Proof of Theorem~8]

By definition, $K=|\cR^{\BH}_{\alpha}(\bP)|$, so
\begin{align*}
    \FDR&=\EE\left[\frac{\sum_{j\in\cH_0}\ind\{P_j\leq q^{m,\alpha}_K\}}{1\vee\sum_{i\in\cH}\ind\{P_i\leq q^{m,\alpha}_K\}}\right]=\sum_{j\in\cH_0}\EE\left[\frac{\ind\{P_j\leq q^{m,\alpha}_K\}}{1\vee\sum_{i\in\cH}\ind\{P_i\leq q^{m,\alpha}_K\}}\right]\\
    &=\sum_{j\in\cH_0}\EE\left[\sum_{k=1}^m\frac{\ind\{P_j\leq q^{m,\alpha}_k\}}{k}\,\ind\{K=k\}\right]\\
    &=\sum_{j\in\cH_0}\sum_{k=1}^m\frac{1}{k}\PP\{P_j\leq q^{m,\alpha}_k,\,|\cR^{\BH}_{\alpha}(\bP)|=k\}.
\end{align*}
We now reformulate this expression in terms of the randomly selected null index $J$. Since $\PP(J=j)=1/|\cH_0|$ for $j\in\cH_0$,
\begin{align*}
    &\PP\left\{P_J\leq q^{m,\alpha}_m,\,|\cR^{\BH}_{\alpha}(\bP)|=k\right\}\\
    =&\sum_{j\in\cH_0}\PP\left\{P_J\leq q^{m,\alpha}_m,\,|\cR^{\BH}_{\alpha}(\bP)|=k,\,J=j\right\}\\
    =&\sum_{j\in\cH_0}\PP\left\{P_j\leq q^{m,\alpha}_m,\,|\cR^{\BH}_{\alpha}(\bP)|=k\;\big|\;J=j\right\}\PP(J=j)\\
    =&\frac{1}{|\cH_0|}\sum_{j\in\cH_0}\PP\left\{P_j\leq q^{m,\alpha}_m,\,|\cR^{\BH}_{\alpha}(\bP)|=k\right\},
\end{align*}
so that, under randomized super-uniformity,
\begin{align*}
    \FDR&=\sum_{j\in\cH_0}\sum_{k=1}^m\frac{1}{k}\PP\{P_j\leq q^{m,\alpha}_k,\,|\cR^{\BH}_{\alpha}(\bP)|=k\}\\
    &=\sum_{k=1}^m\frac{|\cH_0|}{k}\PP\{P_J\leq q^{m,\alpha}_k,\,|\cR^{\BH}_{\alpha}(\bP)|=k\}\\
    &\leq\sum_{k=1}^m\frac{\alpha m_0}{m}\PP\{|\cR^{\BH}_{\alpha}(\bP)|=k\;\big|\;P_J\leq q^{m,\alpha}_k\}.
\end{align*}
Because of the randomness of $J$, the randomized PRDS condition must be applied to the order statistics $\bP^J_{(\cdot)}$, so the set $\{|\cR^{\BH}_{\alpha}(\bP)|\leq k\}$ is not directly non-decreasing and needs to be refined. Define
\[
\cD_k=\left\{\boldsymbol{p}\in\RR_*^{m-1}:p_{(k-1)}\leq q^{m,\alpha}_k,\,p_{(k)}>q^{m,\alpha}_{k+1},\dots,p_{(m-1)}>q^{m,\alpha}_m\right\},
\]
so that $\{\bP^J_{(\cdot)}\in\cD_k\}=\{|\cR^{\BH}_{\alpha}(\bP)|=k\}$ given $P_J\leq q_k^{m,\alpha}$, and define
\[
\cC_k=\left\{\boldsymbol{p}\in\RR_*^{m-1}:p_{(k)}>q^{m,\alpha}_{k+1},\dots,p_{(m-1)}>q^{m,\alpha}_m\right\},
\]
so that $\cD_k=\cC_k\setminus\cC_{k-1}$, with $\cC_0=\emptyset$ and $\cC_m=\RR_*^{m-1}$. Each $\cC_k$ is non-decreasing. Hence, under randomized $\PRDS_{\cH_0}$,
\begin{align*}
    \FDR&\leq\sum_{k=1}^m\frac{\alpha m_0}{m}\PP\{|\cR^{\BH}_{\alpha}(\bP)|=k\;\big|\;P_J\leq q^{m,\alpha}_k\}\\
    &=\sum_{k=1}^m\frac{\alpha m_0}{m}\PP\{\bP^J_{(\cdot)}\in\cD_k\;\big|\;P_J\leq q^{m,\alpha}_k\}\\
    &=\sum_{k=1}^m\frac{\alpha m_0}{m}\left(\PP\{\bP^J_{(\cdot)}\in\cC_k\;\big|\;P_J\leq q^{m,\alpha}_k\}-\PP\{\bP^J_{(\cdot)}\in\cC_{k-1}\;\big|\;P_J\leq q^{m,\alpha}_k\}\right)\\
    &\leq\sum_{k=1}^m\frac{\alpha m_0}{m}\left(\PP\{\bP^J_{(\cdot)}\in\cC_k\;\big|\;P_J\leq q^{m,\alpha}_k\}-\PP\{\bP^J_{(\cdot)}\in\cC_{k-1}\;\big|\;P_J\leq q^{m,\alpha}_{k-1}\}\right)\\
    &\leq\alpha\frac{m_0}{m}.
\end{align*}
\end{proof}

\begin{proof}[Proof of Lemma~2]
By conditioning on 
\[
(Z_{p+1},\dots,Z_{p+q})=(X_{p+1},\dots,X_{p+q}),
\]
it suffices to treat the case $q=0$, for which $p=m$. Let $\pi$ be any permutation of $[m]$, and fix a measurable set $A\subseteq\RR_*^m$. Since $\mathbf{J}$ is a uniformly random permutation of $[m]$ independent of $\bX$, the vector $(J_{\pi(1)},\dots,J_{\pi(m)})$ has the same distribution as $(J_1,\dots,J_m)$. Hence
\begin{align*}
\PP\{\bZ\in A\}=&\PP\{(X_{J_1},\dots,X_{J_m})\in A\}=\PP\{(X_{J_{\pi(1)}},\dots,X_{J_{\pi(m)}})\in A\}\\
=&\PP\{(Z_{\pi(1)},\dots,Z_{\pi(m)})\in A\},
\end{align*}
where the last equality holds because $Z_i=X_{J_i}$. As $\pi$ is arbitrary, the distribution of $\bZ$ is invariant under coordinate permutations, so $\bZ$ is exchangeable.
\end{proof}

\begin{proof}[Proof of Lemma~3]
Without loss of generality, we may assume that $\bX$ is degenerate, i.e., $X_1<\dots<X_m$ are fixed strictly ordered constants, so that the claim reduces to showing that the induced map $\mathbf{S}$ is a uniformly random permutation of $[m]$. For $\tau=0,\dots,m$, let $\boldsymbol{l}=(l_1,\dots,l_m)$ be any permutation of $[m]$, and condition on $|\cA|=\tau$. The event $\{S_1=l_1,\dots,S_\tau=l_\tau\}$ occurs precisely when $\cA=\{l_1,\dots,l_\tau\}$ and $\mathbf{J}$ is the ordering $(l_1,\dots,l_\tau)$ of $\cA$. By the definition of a uniformly random split,
\[
\PP\{S_1=l_1,\dots,S_\tau=l_\tau\mid|\cA|=\tau\}
=\binom{m}{\tau}^{-1}\frac{1}{\tau!}
=\frac{(m-\tau)!}{m!}.
\]
Conditional on this event, $(S_{\tau+1},\dots,S_m)$ is a uniformly random ordering of $\cB=[m]\setminus\{l_1,\dots,l_\tau\}$, whence
\[
\PP\{S_{\tau+1}=l_{\tau+1},\dots,S_m=l_m\mid S_1=l_1,\dots,S_\tau=l_\tau\}
=\frac{1}{(m-\tau)!}.
\]
Multiplying the two displays,
\[
\PP\{\mathbf{S}=\boldsymbol{l}\mid|\cA|=\tau\}
=\frac{(m-\tau)!}{m!}\cdot\frac{1}{(m-\tau)!}
=\frac{1}{m!}.
\]
As $\boldsymbol{l}$ and $\tau$ are arbitrary,
\[
\PP(\mathbf{S}=\boldsymbol{l})=\EE[\PP\{\mathbf{S}=\boldsymbol{l}\mid|\cA|\}]=\frac{1}{m!},
\]
it follows that $\mathbf{S}$ is a uniformly random permutation of $[m]$.
\end{proof}

\begin{proof}[Proof of Corollary~1]
If $\cH_0^{\test}=\emptyset$, then all null hypotheses lie in the calibration set and have p-values $+\infty$, so that the BH procedure rejects none of them and the FDR is zero, making the bound trivial. We henceforth assume $\cH_0^{\test}\neq\emptyset$, which in particular gives $\cH^{\test}\neq\emptyset$ and ensures that the ratio $|\cH_0^{\test}|/|\cH^{\test}|$ below is well defined. Condition on $|\cH_0^{\test}|=m_{0t}$, so $|\cH_0^{\rcal}|=m_0-m_{0t}$. Let $J_1,\dots,J_{m_{0t}}$ be a uniformly random permutation of $\cH_0^{\test}$ and $J_{m_{0t}+1},\dots,J_{m_0}$ a uniformly random permutation of $\cH_0^{\rcal}$, independent of $\bX$. Set $Z_i=X_{J_i}$ for $i=1,\dots,m_0$ and $Z_i=X_i$ for $i=m_0+1,\dots,m$.

Define conformal p-values $\bU\in\RR_*^{m-m_0+m_{0t}}$ with calibration set $\{m_{0t}+1,\dots,m_0\}$ and test set $\{1,\dots,m_{0t}\}\cup\{m_0+1,\dots,m\}$. For the null test coordinates,
\[
U_j=\frac{1+\sum_{i=m_{0t}+1}^{m_0}\ind\{Z_i\geq Z_j\}}{1+m_0-m_{0t}},\quad j=1,\dots,m_{0t},
\]
and for the non-null test coordinates,
\[
U_j=\frac{1+\sum_{i=m_{0t}+1}^{m_0}\ind\{Z_i\geq Z_{j+m_0-m_{0t}}\}}{1+m_0-m_{0t}},\quad j=m_{0t}+1,\dots,m_1+m_{0t}.
\]
By Lemma~3, conditionally on $|\cH_0^{\test}|=m_{0t}$, the sequence $(J_1,\dots,J_{m_0})$ is a uniformly random permutation of $\cH_0$. By Lemma~2, $(Z_1,\dots,Z_{m_0})$ are exchangeable conditional on $(Z_{m_0+1},\dots,Z_m)$. Hence $U_1,\dots,U_{m_{0t}}$ are super-uniform and $\bU$ satisfies the PRDS condition with respect to the null indices $\{1,\dots,m_{0t}\}$. Moreover, since $J_1$ is uniformly distributed over $\cH_0^{\test}$ and independent of $\bX$, and $P_{J_1}=U_1$ with $\bP^{J_1}_{(\cdot)}=\bU^{1}_{(\cdot)}$, we have
\[
\PP\{P_{J_1}\leq t\mid |\cH_0^{\test}|\}=\PP\{U_1\leq t\mid |\cH_0^{\test}|\}\leq t,\quad t\in[0,1],
\]
and for any non-decreasing set $C\subseteq\RR_*^{m-m_0+m_{0t}-1}$,
\[
\PP\{\bP^{J_1}_{(\cdot)}\in C\mid P_{J_1}\leq t,\ |\cH_0^{\test}|\}=\PP\{\bU^{1}_{(\cdot)}\in C\mid U_1\leq t,\ |\cH_0^{\test}|\}
\]
is non-decreasing in $t$. Then let $J$ be $J_1$ with probability $|\cH_0^{\test}|/m_0$ and $J_{m_{0t}+1}$ with probability $|\cH_0^{\rcal}|/m_0$, so that $J$ is uniformly distributed over $\cH_0$ and independent of $\bX$. Define $\bQ$ by $Q_j=P_j$ for $j\in\cH^{\test}$ and $Q_j=+\infty$ for $j\in\cH_0^{\rcal}$. Since $Q_{J_{m_{0t}+1}}=+\infty$, the event $\{Q_J\leq t\}$ forces $J=J_1$ almost surely, so that
\[
\PP\{Q_J\leq t\mid |\cH_0^{\test}|\}=\frac{m_{0t}}{m_0}\,\PP\{P_{J_1}\leq t\mid |\cH_0^{\test}|\}\leq \frac{m_{0t}}{m_0}t,\quad t\in[0,1].
\]
For any non-decreasing set $C\subseteq\RR_*^{m_1+m_{0t}-1}$, define
\[
D=\{\boldsymbol{x}\in\RR_*^{m-1}:(x_1,\dots,x_{m_1+m_{0t}-1})\in C\}.
\]
Then $D$ is non-decreasing. On the event $\{Q_J\leq t\}$ we have $J=J_1$, so that
\[
\bQ^{J}_{(\cdot)}=(\bP^{J_1}_{(\cdot)},+\infty,\dots,+\infty),
\]
and hence
\[
\PP\{\bQ^{J}_{(\cdot)}\in D\mid Q_{J}\leq t,\ |\cH_0^{\test}|\}=\PP\{\bP^{J_1}_{(\cdot)}\in C\mid P_{J_1}\leq t,\ |\cH_0^{\test}|\}
\]
is non-decreasing in $t$. Therefore $\bQ$ is randomized BH-valid with respect to $\cH_0$ and weight $m_{0t}/m_0$, i.e., with the super-uniformity bound $\PP\{Q_J\leq t\}\leq m_{0t}t/m_0$ in place of $t$.
\begin{Lemma}
\label{slem:level_dim}
Let $\bP\in\RR_*^m$ and $S\subseteq[m]$. If $P_j=+\infty$ for all $j\notin S$, then for any $\alpha\in(0,1)$,
\[
\cR_{m\alpha/|S|}(\bP)=\cR_{\alpha}(\bP_S),
\]
where $\cR_{\alpha}(\bP_S)$ denotes the BH procedure applied to the subvector $\bP_S$.
\end{Lemma}
By Lemma~\ref{slem:level_dim} applied to $\bQ$ with $S=\cH^{\test}$, we have $\cR_\alpha(\bP)=\cR_\beta(\bQ)$, where $\beta=m\alpha/(m_1+m_{0t})$. Conditionally on $|\cH_0^{\test}|$, the argument of Theorem~8, with $\PP\{Q_J\leq t\}\leq m_{0t}t/m_0$ in place of $\PP\{P_J\leq t\}\leq t$, bounds the conditional FDR of $\bQ$ at level $\beta$ by
\[
\frac{\beta}{m}\sum_{j\in\cH_0}\frac{m_{0t}}{m_0}.
\]
Hence
\begin{align*}
    \FDR&=\EE\left[\frac{|\cH_0\cap\cR_\alpha(\bP)|}{1\vee|\cR_\alpha(\bP)|}\right]
    =\EE\left[\frac{|\cH_0\cap\cR_\beta(\bQ)|}{1\vee|\cR_\beta(\bQ)|}\right]\\
    &\leq\EE\left[\frac{\sum_{j\in\cH_0}\frac{m_{0t}}{m_0}}{m}\,\beta\right]
    =\alpha\,\EE\left[\frac{m_{0t}}{m_{0t}+m_1}\right]=\alpha\,\EE\left[\frac{|\cH^{\test}_0|}{|\cH^{\test}|}\right].
\end{align*}
\end{proof}

\section{Embedding Theorems}
\label{asec:embedding}

\subsection{Proofs of the Three Embedding Theorems}

\begin{proof}[Proof of Theorem~9]
(i)~For any $J$ uniformly distributed over $\cH_0$ and independent of $\bP$, and any non-decreasing set $C\subseteq\RR_*^{m-1}$,
\begin{align*}
    \PP\{\bP^J_{(\cdot)}\in C\mid P_J\leq t\}
    &=\sum_{j\in\cH_0}\PP\{\bP^j_{(\cdot)}\in C\mid P_J\leq t,J=j\}\PP\{J=j\mid P_J\leq t\}\\
    &=\sum_{j\in\cH_0}\PP\{\bP^j_{(\cdot)}\in C\mid P_j\leq t\}\frac{\PP\{J=j\}\PP\{P_j\leq t\}}{\sum_{k\in\cH_0}\PP\{J=k\}\PP\{P_k\leq t\}}.
\end{align*}
Since $J$ is uniform over $\cH_0$ and the $P_j$ are identically distributed for $j\in\cH_0$, the weight reduces to $1/|\cH_0|$, so that
\[
\PP\{\bP^J_{(\cdot)}\in C\mid P_J\leq t\}=\frac{1}{|\cH_0|}\sum_{j\in\cH_0}\PP\{\bP^j_{(\cdot)}\in C\mid P_j\leq t\}.
\]
Each term on the right is non-decreasing in $t$ by the PRDS condition, hence so is their average, and $\bP$ is randomized PRDS with respect to $\cH_0$. Moreover, by super-uniformity,
\[
\PP\{P_J\leq t\}=\frac{1}{|\cH_0|}\sum_{j\in\cH_0}\PP\{P_j\leq t\}\leq t,
\]
so that $\bP$ is randomized BH-valid with respect to $\cH_0$.
(ii)~Theorem~2 shows that conformal p-values are BH-valid with respect to $\cH_0$. Moreover, by the definition of conformal p-values and the exchangeability assumption, for any $j\in\cH_0$,
\[
\PP\{P_j=\frac{1+k}{1+n}\}=\frac{1}{1+n},\quad k=0,\dots,n,
\]
so that the $P_j$ are identically distributed for $j\in\cH_0$. Part~(i) now yields that conformal p-values are randomized BH-valid with respect to $\cH_0$.
\end{proof}

\begin{proof}[Proof of Theorem~10]
Fix $j\in\cH_0$, and let $F_j$ be the cumulative distribution function of $P_j$. Since $P_j$ is super-uniform, $F_j(t)\leq t$ for all $t\in[0,1]$. Let $\mathcal{I}_j\subseteq[0,1]$ be the set of discontinuity points of $F_j$, and define
\[
Q_j=
\begin{cases}
F_j(P_j), & P_j\in[0,1]\setminus\mathcal{I}_j,\\[4pt]
\mathrm{Unif}\bigl[\lim_{t\to P_j-}F_j(t),\;\lim_{t\to P_j+}F_j(t)\bigr], & P_j\in\mathcal{I}_j,\\[4pt]
\mathrm{Unif}\bigl[F_j(1),\;1\bigr], & P_j\in(1,\infty).
\end{cases}
\]
By construction $Q_j\leq P_j$ almost surely, so $\PP(\bP\geq\bQ)=1$, and $Q_j$ is uniformly distributed on $[0,1]$. Moreover, $Q_j$ depends on $(\bP,\bQ)$ only through $P_j$ and the independent randomization used in the second and third cases, so $Q_j\perp\!\!\!\perp(\bP_{-j},\bQ_{-j})\mid P_j$. The quantile function of $F_j$ recovers $P_j\wedge1$ from $Q_j$, so that $\sigma(\bP\wedge\mathbf{1})\subseteq\sigma(\bQ)$. In what follows we restrict $t$ to $[0,1]$, as the BH critical values lie in this interval.

Now let $J$ be uniformly distributed over $\cH_0$ and independent of $(\bP,\bQ)$, and fix a non-decreasing set $C\subseteq[0,1]^{m-1}$. Then
\begin{align*}
    \PP\{\bP^J_{(\cdot)}\in C\mid Q_J\leq t\}
    &=\sum_{j\in\cH_0}\PP\{\bP^j_{(\cdot)}\in C\mid Q_J\leq t,J=j\}\PP\{J=j\mid Q_J\leq t\}\\
    &=\sum_{j\in\cH_0}\PP\{\bP^j_{(\cdot)}\in C\mid Q_j\leq t\}\frac{\PP\{J=j\}\PP\{Q_j\leq t\}}{\sum_{k\in\cH_0}\PP\{J=k\}\PP\{Q_k\leq t\}}.
\end{align*}
Since the $Q_j$ are identically distributed for $j\in\cH_0$, the weight reduces to $1/|\cH_0|$, so that
\[
\PP\{\bP^J_{(\cdot)}\in C\mid Q_J\leq t\}=\frac{1}{|\cH_0|}\sum_{j\in\cH_0}\PP\{\bP^j_{(\cdot)}\in C\mid Q_j\leq t\}.
\]

It remains to show that $H_j(t)=\PP\{\bP^j_{(\cdot)}\in C\mid Q_j\leq t\}$ is non-decreasing in $t$. Since $Q_j$ is uniform on $[0,1]$, $H_j$ is continuous on $(0,1]$; indeed, for $0<t<t'$,
\[
|H_j(t')-H_j(t)|\leq\frac{t'-t}{t'}\to 0\quad\text{as }t'\downarrow t.
\]
It therefore suffices to examine the behaviour of $H_j$ on the jump intervals of $F_j$, which correspond to the (countably many) atoms of $P_j$. Fix a jump point $t_0$, write $a=F_j(t_0-)<F_j(t_0)=b$, and $P_j=t_0$ carries mass $b-a$. Since $Q_j\perp\!\!\!\perp\bP^j\mid P_j$, for $u\in(a,b)$ the event $\{Q_j\leq u\}$ splits as $\{P_j<t_0\}\cup\{P_j=t_0,\;Q_j\leq u\}$, whence
\begin{align*}
    H_j(u)&=\frac{a\,\PP\{\bP^j_{(\cdot)}\in C\mid P_j<t_0\}+(u-a)\,\PP\{\bP^j_{(\cdot)}\in C\mid P_j=t_0\}}{u},\\
    H_j'(u)&=\frac{a}{u^2}\Bigl(\PP\{\bP^j_{(\cdot)}\in C\mid P_j=t_0\}-\PP\{\bP^j_{(\cdot)}\in C\mid P_j<t_0\}\Bigr).
\end{align*}
Since the two conditional probabilities on the right-hand side are constant in $u\in(a,b)$, $H_j$ is differentiable on $(a,b)$ and its derivative is as displayed.
The difference in parentheses is nonnegative: by the $\PRDS_{\cH_0}$ condition, $t\mapsto\PP\{\bP^j_{(\cdot)}\in C\mid P_j\leq t\}$ is non-decreasing, then for $\varepsilon>0$ small enough that $t_0-\varepsilon\in(a,b)$, we have
\begin{align*}
    &\PP\{\bP^j_{(\cdot)}\in C\mid t_0-\varepsilon<P_j\leq t_0\}-\PP\{\bP^j_{(\cdot)}\in C\mid P_j\leq t_0-\varepsilon\}\\
    =&\frac{\PP\{\bP^j_{(\cdot)}\in C,\,t_0-\varepsilon<P_j\leq t_0\}}{\PP\{t_0-\varepsilon<P_j\leq t_0\}}-\frac{\PP\{\bP^j_{(\cdot)}\in C,\,P_j\leq t_0-\varepsilon\}}{\PP\{P_j\leq t_0-\varepsilon\}}\\
    =&\frac{\PP(\bP^j_{(\cdot)}\in C,\,P_j\leq t_0)-\frac{\PP(P_j\leq t_0)}{\PP(P_j\leq t_0-\varepsilon)}\PP(\bP^j_{(\cdot)}\in C,\,P_j\leq t_0-\varepsilon)}{\PP(t_0-\varepsilon<P_j\leq t_0)}\geq 0.
\end{align*}
Hence $H_j'\geq 0$ on each jump interval, and $H_j$ is non-decreasing on $[0,1]$. By the averaging identity above, the same holds for $t\mapsto\PP\{\bP^J_{(\cdot)}\in C\mid Q_J\leq t\}$, and together with the uniformity of $Q_j$ and the domination $\PP(\bP\geq\bQ)=1$, this shows that $\bP$ is dominated randomized BH-valid with respect to $\cH_0$ and $\bQ$.
\end{proof}

\begin{proof}[Proof of Lemma~4]
Since $W_j\in\RR_*=[0,\infty]$, we have $W_j\geq 0$, so that $\cR_+(0)=\{j\in\cH:L_j=1\}$ and hence $\cV_+(0)=\{j\in\cH_0:L_j=1\}$. Conditional on $\sigma_{\cH_1}$, the set $\cV_+(0)$ is uniformly distributed over all subsets of $\cH_0$ of a given size, in the sense that
\[
\PP\{\cV_+(0)=A\mid\sigma_{\cH_1}\}=\binom{m_0}{|A|}^{-1},\quad\forall A\subseteq\cH_0.
\]
By the competition condition, $\PP\{L_j=1\mid\bW,\bL_{-j}\}=r/(1+r)$ for every $j\in\cH_0$. Iterating over $j\in\cH_0$, conditional on $\bW$ and $\bL_{\cH_1}$, the labels $\{L_j\}_{j\in\cH_0}$ are independent Bernoulli$(r/(1+r))$ variables, so that for any $\boldsymbol{l}\in\{0,1\}^{m_0}$,
\[
\PP\{\bL_{\cH_0}=\boldsymbol{l}\mid\bW,\bL_{\cH_1}\}
=\left(\frac{r}{1+r}\right)^{\sum_{i=1}^{m_0}l_i}\left(\frac{1}{1+r}\right)^{m_0-\sum_{i=1}^{m_0}l_i}.
\]
Let $\mathcal{L}_{k}=\{\boldsymbol{l}\in\{0,1\}^{m_0}:\sum_{i=1}^{m_0}l_i=k\}$, so that $|\mathcal{L}_k|=m_0!/(k!(m_0-k)!)$. Since $V_+(0)=\sum_{i=1}^{m_0}l_i$ on the event $\{\bL_{\cH_0}=\boldsymbol{l}\}$, we have $\{V_+(0)=k,V_-(0)=m_0-k\}=\{\bL_{\cH_0}\in\mathcal{L}_k\}$. Summing the display over $\mathcal{L}_k$ gives
\[
\PP\{V_+(0)=k,V_-(0)=m_0-k\mid\bW,\bL_{\cH_1}\}
=\binom{m_0}{k}\left(\frac{r}{1+r}\right)^{k}\left(\frac{1}{1+r}\right)^{m_0-k}.
\]
Dividing the product form by this binomial mass shows that, conditional on $V_+(0)=k$, the labels are uniform over $\mathcal{L}_k$. Hence, conditional on $\sigma_{\cH_1}$,
\[
\PP\{\bL_{\cH_0}=\boldsymbol{l}\mid\sigma_{\cH_1}\}=
\begin{cases}
\dfrac{1}{|\mathcal{L}_{V_+(0)}|}, & \boldsymbol{l}\in\mathcal{L}_{V_+(0)},\\[6pt]
0, & \boldsymbol{l}\in\{0,1\}^{m_0}\setminus\mathcal{L}_{V_+(0)}.
\end{cases}
\]
In other words, the set $\cV_+(0)=\{i\in\cH_0:L_i=1\}$ is uniformly distributed over all subsets of $\cH_0$ of size $V_+(0)$, and the labels $\bL$ induce a uniformly random split of $\cH_0$ into $\cH_0=\cV_+(0)+\cV_-(0)$.
\end{proof}

\begin{proof}[Proof of Theorem~11]
We consider the $\sigma$-algebra $\sigma_{w,\cH_1}=\sigma(\bW,\bL_{\cH_1})$ generated by all scores and the non-null labels as the given condition. By the definition of competition statistics, conditional on $\sigma_{w,\cH_1}$, the null labels $\{L_j\}_{j\in\cH_0}$ are independent Bernoulli$(r/(1+r))$ variables. Lemma~4 then implies that $\bL_{\cH_0}$ induces a uniformly random split $\cH_0=\cV_+(0)+\cV_-(0)$ conditional on $\sigma_{w,\cH_1}$. By Lemmas~2 and~3 and conditional on $\sigma_v=\sigma(V_+(0))$, the random split generates a random permutation of $\cH_0$, under which $\bQ$ is randomized BH-valid with respect to $\cH_0$ conditional on $\sigma_{w,\cH_1}$.

It remains to lift the randomized BH-validity from $\bQ$ to $\bP$. From the definitions,
\begin{align*}
    P_j&=\frac{1+R_-(W_j)}{1+V_-(0)} \\
       &=\frac{1+V_-(W_j)+\sum_{i\in\cH_1}\ind\{W_i\geq W_j,L_i=0\}}{1+V_-(0)} \\
       &=Q_j+\frac{\sum_{i\in\cH_1}\ind\{W_i\geq W_j,L_i=0\}}{1+V_-(0)}, && L_j=1,\\[4pt]
    P_j&=Q_j=+\infty, && L_j=0.
\end{align*}
The correction term is nonnegative. Moreover,
\[
\sigma\Bigl(\sum_{i\in\cH_1}\ind\{W_i\geq W_j,L_i=0\}\Bigr)\subseteq\sigma_{w,\cH_1},
\quad
\sigma(V_-(0))\subseteq\sigma_v,
\]
so the correction term is $\sigma_v\vee\sigma_{w,\cH_1}$-measurable. Hence $\bP\geq\bQ$ almost surely and the map $\phi\colon\bQ\mapsto\bP$ is a non-decreasing bijection given $\sigma_{w,\cH_1}$.

\begin{Lemma}
\label{slem:rank_indep}
Let $X_1,\dots,X_m$ be continuous random variables, and let $M\in\{1,\dots,m-1\}$ be independent of $(X_1,\dots,X_m)$. Choose a random subset $I\subseteq\{1,\dots,m\}$ of size $M$ uniformly over all such subsets, independently of the $X_i$'s given $M$, and an index $j$ uniformly from $[m]\setminus I$, independently of all previous choices. Set $S=\{X_i:i\in I\}$, $Y=X_j$, and
\[
R = 1 + \sum_{i \in I} \ind_{\{X_i < Y\}}\;\in\;\{1,\dots,M+1\}.
\]
Then $R$ is independent of $(X_1,\dots,X_m)$.
\end{Lemma}
By this lemma, the random variable $Q_J$ is independent of $(\bW,\bL_{\cH_1})$. Consequently
\[
\PP(\bW,\bL_{\cH_1}\mid Q_J\leq t,\sigma_v)=\PP(\bW,\bL_{\cH_1}\mid\sigma_v),\quad t\in[0,1],
\]
and for any non-decreasing set $C\subseteq\RR_*^{m-1}$, the preimage $\phi^{-1}(C)$ is also non-decreasing. Thus
\begin{align*}
&\PP(\bP^J_{(\cdot)}\in C\mid Q_J\leq t,\sigma_v)\\
=&\int_{\bW,\bL_{\cH_1}}\PP(\bP^J_{(\cdot)}\in C\mid Q_J\leq t,\sigma_{w,\cH_1})\rd \PP(\bW,\bL_{\cH_1}\mid Q_J\leq t,\sigma_v)\\
=&\int_{\bW,\bL_{\cH_1}}\PP(\bP^J_{(\cdot)}\in C\mid Q_J\leq t,\sigma_{w,\cH_1})\rd \PP(\bW,\bL_{\cH_1}\mid\sigma_v)\\
\end{align*}
which is non-decreasing in $t$ because the inner conditional probability is non-decreasing by the conditional randomized PRDS of $\bQ$ and the outer expectation preserves monotonicity. Together with $\PP(Q_J\leq t)=\EE[\PP(Q_J\leq t\mid \sigma_{w,\cH_1})]\leq t$, this shows that $\bP$ is dominated randomized BH-valid with respect to $\cH_0$ and $\bQ$ when conditioned on $\sigma_v$.

To establish the FDR bound, we consider the conditional expectation first. Conditional on $\sigma_v$, the random index $J$ is uniform over $\cH_0$, and $Q_j<+\infty$ only for $j\in\cV_+(0)$. Hence
\[
\PP(Q_J\leq t\mid\sigma_v)\leq t\,\frac{V_+(0)}{m_0}\leq t,\qquad t\in[0,1].
\]
Let $K=|\cR^{\BH}(\bP')|$. Define the threshold vector
\[
\boldsymbol{q}^{E,\alpha}
=\Bigl(q_1^{R_+(0),\alpha}\frac{1\vee R_+(0)}{r+rV_-(0)},\dots,
q_m^{R_+(0),\alpha}\frac{1\vee R_+(0)}{r+rV_-(0)}\Bigr),
\]
which is $\sigma_v\vee\sigma_{w,\cH_1}$-measurable and hence deterministic under our conditioning. Introduce the level sets
\[
\cC_k(\boldsymbol{q})=\{\boldsymbol{p}\in\RR_*^{m-1}:p_{(k-1)}\leq q_k,\;p_{(k)}>q_{k+1},\dots,p_{(m-1)}>q_m\},
\]
\[
\cD_k(\boldsymbol{q})=\{\boldsymbol{p}\in\RR_*^{m-1}:p_{(k)}>q_{k+1},\dots,p_{(m-1)}>q_m\},
\]
so that $\cC_k(\boldsymbol{q})=\cD_k(\boldsymbol{q})\setminus\cD_{k-1}(\boldsymbol{q})$ with $\cD_0=\emptyset,\cD_m = \RR_*^{m-1}$; each $\cD_k$ is non-decreasing. By definition of the BH procedure, $K=k$ if and only if $\bP_{(\cdot)}\in\cC_k(\boldsymbol{q}^{E,\alpha})$ conditional on $P_J\leq q^{R_+(0),\alpha}_kE_J$.

Conditioning on $\sigma_v\vee\sigma_{w,\cH_1}$ and applying the tower property,
\begin{align*}
    \FDR
    &=\EE\left[\frac{\sum_{j\in\cH_0}\ind\{P_j\leq q^{R_+(0),\alpha}_KE_j\}}{1\vee K}\right]\\
    &=\EE\left[\EE\left[\frac{\sum_{j\in\cH_0}\ind\{P_j\leq q^{R_+(0),\alpha}_KE_j\}}{1\vee K}\;\Big|\;\sigma_v\vee\sigma_{w,\cH_1}\right]\right]\\
    &\leq\EE\left[\sum_{k=1}^m\frac{m_0}{k}\,
    \PP\bigl(Q_J\leq q^{R_+(0),\alpha}_kE_J,\;
    \bP^J_{(\cdot)}\in\cC_k(\boldsymbol{q}^{E,\alpha})
    \mid\sigma_v\vee\sigma_{w,\cH_1}\bigr)\right].
\end{align*}
Factor the joint probability as
\begin{align*}
    &\PP\bigl(Q_J\leq q^{R_+(0),\alpha}_kE_J,\;
    \bP^J_{(\cdot)}\in\cC_k(\boldsymbol{q}^{E,\alpha})
    \mid\sigma_v\vee\sigma_{w,\cH_1}\bigr)\\
    =&\PP\bigl(\bP^J_{(\cdot)}\in\cC_k(\boldsymbol{q}^{E,\alpha})
    \mid Q_J\leq q^{R_+(0),\alpha}_kE_J,\;\sigma_v\vee\sigma_{w,\cH_1}\bigr)\,
    \PP\bigl(Q_J\leq q^{R_+(0),\alpha}_kE_J\mid\sigma_v\vee\sigma_{w,\cH_1}\bigr).
\end{align*}
From the super-uniformity bound at the beginning of this section,
\[
\PP\bigl(Q_J\leq q^{R_+(0),\alpha}_kE_J\mid\sigma_v\vee\sigma_{w,\cH_1}\bigr)
\leq\frac{V_+(0)}{m_0}\,q^{R_+(0),\alpha}_kE_J
=\frac{k\alpha}{m_0}\,\frac{V_+(0)}{r+rV_-(0)}.
\]
Substituting and cancelling the factor $m_0/k$ yields
\[
\FDR\leq\EE\left[\sum_{k=1}^m\alpha\frac{V_+(0)}{r+rV_-(0)}\,
\PP\bigl(\bP^J_{(\cdot)}\in\cC_k(\boldsymbol{q}^{E,\alpha})
\mid Q_J\leq q^{R_+(0),\alpha}_kE_J,\;\sigma_v\vee\sigma_{w,\cH_1}\bigr)\right].
\]
Now $\cC_k=\cD_k\setminus\cD_{k-1}$ with $\cD_k$ non-decreasing. By the dominated randomized $\PRDS_{\cH_0}$ of $\bP$ with respect to $\bQ$, for $t_{k-1}\leq t_k$,
\[
\PP(\bP^J_{(\cdot)}\in\cD_{k-1}\mid Q_J\leq t_k,\sigma_v\vee\sigma_{w,\cH_1})
\geq\PP(\bP^J_{(\cdot)}\in\cD_{k-1}\mid Q_J\leq t_{k-1},\sigma_v\vee\sigma_{w,\cH_1}).
\]
Hence we have
\begin{align*}
    \FDR
    &\leq\EE\left[\alpha\frac{V_+(0)}{r+rV_-(0)}
    \sum_{k=1}^m\PP(\bP^J_{(\cdot)}\in\cC_k\mid Q_J\leq q^{R_+(0),\alpha}_kE,\;\sigma_v\vee\sigma_{w,\cH_1})\right]\\
    &\leq\EE\left[\alpha\frac{V_+(0)}{r+rV_-(0)}\right]\leq\alpha.
\end{align*}
This completes the FDR bound.
\end{proof}

\subsection{Embedding of Sub-competition Statistics}

Sub-competition statistics lack the label-independence of competition statistics: a null label may obey the sub-competition bound $\PP\{L_j=1\mid\bW,\bL_{-j}\}\leq r/(1+r)$ while still depending on the other labels, for instance as
\[
\PP\{L_j=1\mid \bW,\bL_{-j}\}
= \frac{1}{2}\frac{r}{1+r}+f_j(\bW,\bL_{-j})\frac{1}{2}\frac{1}{1+r},
\qquad j\in\cH_0,
\]
with $f_j(\bW,\bL_{-j})\in[0,1]$. Subsection~\ref{sec:bridge} constructs $\bZ=\bL+\bY$ whose null labels satisfy the competition condition $\PP\{Z_j=1\mid\bW,\bZ_{-j}\}=r/(1+r)$. Since the rejection set depends on $(\bW,\bL)$ rather than $(\bW,\bZ)$, a second step is required to couple $\bZ$ back to $\bL$.

However, the label-dependence may cause issues in the analysis. We suppose in addition that for each $j\in\cH_0$ there exists $r_j\leq r$ such that
\[
\PP\{L_j=1\mid \bW,\bL_{-j}\}=r_j\,\PP\{L_j=0\mid \bW,\bL_{-j}\},\qquad j\in\cH_0,
\]
where $r_j$ is a constant but not a function of $\bW$ and $\bL_{-j}$. Construct two independent families of labels. First, let $\bZ$ satisfy
\[
\PP\{Z_j=1\mid \bW,\bZ_{-j}\}=\frac{r}{1+r},\qquad j\in\cH_0,
\]
with $\bZ_{\cH_1}=\bL_{\cH_1}$. Second, let $\bZ'$ be independent of $\bZ$ and satisfy
\[
\PP\{Z'_j=1\mid \bZ'_{-j}\}=1-\frac{r_j(1+r)}{r(1+r_j)},\qquad
\PP\{Z'_j=0\mid \bZ'_{-j}\}=\frac{r_j(1+r)}{r(1+r_j)}
\]
for $j\in\cH_0$, with $\bZ'_{\cH_1}=\mathbf{0}$.
Define $\bZ^*=\bZ-\bZ\circ\bZ'$, where $\boldsymbol{a}\circ\boldsymbol{b}=(a_1b_1,\dots,a_mb_m)$ denotes the Hadamard product, so that $Z^*_j=1$ if and only if $Z_j=1$ and $Z'_j=0$. For any $j\in\cH_0$, using the independence of $\bZ$ and $\bZ'$,
\begin{align*}
\PP\{Z^*_j=1\mid \bW,\bZ^*_{-j}\}
&=\PP\{Z_j=1,Z'_j=0\mid \bW,\bZ^*_{-j}\}\\
&=\frac{\sum_{\bZ_{-j},\bZ'_{-j}}\PP\{Z_j=1,\bZ_{-j}\mid \bW\}\,\PP\{Z'_j=0,\bZ'_{-j}\}}
{\sum_{\bZ_{-j},\bZ'_{-j}}\PP\{\bZ_{-j}\mid \bW\}\,\PP\{\bZ'_{-j}\}}\\[2pt]
&=\frac{r}{1+r}\cdot\frac{r_j(1+r)}{r(1+r_j)}
=\frac{r_j}{1+r_j}.
\end{align*}
Hence $(\bW,\bZ^*)\overset{d}{=}(\bW,\bL)$. By the following elementary lemma, the FDR of any procedure depending only on $(\bW,\bL)$ equals that of the same procedure applied to $(\bW,\bZ^*)$.

\begin{Lemma}
\label{slem:fdr_eq}
If $(\bW,\bL)\overset{d}{=}(\bW',\bL')$, then for the competition procedure $\cR^{\cp}_{r,\alpha}$,
\[
\EE\left[\frac{|\cH_0\cap\cR^{\cp}_{r,\alpha}(\bW,\bL)|}
{1\vee|\cR^{\cp}_{r,\alpha}(\bW,\bL)|}\right]
=\EE\left[\frac{|\cH_0\cap\cR^{\cp}_{r,\alpha}(\bW',\bL')|}
{1\vee|\cR^{\cp}_{r,\alpha}(\bW',\bL')|}\right].
\]
\end{Lemma}

With this bridge in place, the randomized BH embedding proceeds in parallel with the proof of Theorem~11, replacing $(\bW,\bL)$ by $(\bW,\bZ^*)$ and using $(\bW,\bZ)$ as the dominating competition structure. Define the relevant counts
\[
R^{Z}_{+}(t)=\sum_{j}\ind\{W_j\geq t,\,Z_j=1\},\quad
V^{Z}_{+}(t)=\sum_{j\in\cH_0}\ind\{W_j\geq t,\,Z_j=1\},
\]
and analogously $R^{Z}_{-},V^{Z}_{-}$ using $Z_j=0$, and $R^{*}_{+},V^{*}_{+},R^{*}_{-},V^{*}_{-}$ using $Z^*_j$. Construct, for $j\in\cR^{*}_{+}(0)$,
\[
P_j=\frac{1+R^{*}_{-}(W_j)}{1+V^{Z}_{-}(0)},\qquad
Q_j=\frac{1+V^{Z}_{-}(W_j)}{1+V^{Z}_{-}(0)},\qquad
E_j=\frac{1\vee R^{Z}_{+}(0)}{r+rV^{Z}_{-}(0)},
\]
and for $j\in\cR^{*}_{-}(0)$, set $P_j=Q_j=+\infty$. Since $Z^*_j=1$ implies $Z_j=1$ (because $Z^*_j=Z_j(1-Z'_j)$), we have $\cR^{*}_{+}(t)\subseteq\cR^{Z}_{+}(t)$ for all $t$, and consequently $R^{*}_{-}(W_j)\geq R^{Z}_{-}(W_j)$. For any $j\in\cR^{*}_{+}(0)$,
\[
P_j=\frac{1+R^{*}_{-}(W_j)}{1+V^{Z}_{-}(0)}
\geq\frac{1+R^{Z}_{-}(W_j)}{1+V^{Z}_{-}(0)}
=Q_j+\frac{\sum_{i\in\cH_1}\ind\{W_i\geq W_j,\,Z_i=0\}}{1+V^{Z}_{-}(0)}
\geq Q_j,
\]
and $\bQ$ is randomized BH-valid. Notice that
\[
r\frac{1+R^*_-(W_j)}{1\vee R^*_+(W_j)}\leq \alpha \quad\text{and}\quad Z^*_j=1\iff\frac{P_j}{E_j}\leq q^{R^Z_+(0),\alpha}_{R^*_+(W_j)}=q^{R^*_+(0),\alpha}_{R^*_+(W_j)}\frac{R^*_+(0)}{R^Z_+(0)},
\]
which implies that, with the analogous reformulation theorem, the rejection set $\cR^{\cp}_{r,\alpha}(\bW,\bZ^*)$ is equal to the rejection set of a BH procedure
\[
\cR=\left\{j\in[m]:P_j\leq q^{R^Z_+(0),\alpha}_{K}\frac{1\vee R^*_+(0)}{r+rV^Z_-(0)}\right\},
\]
where
\[
K=\sup\left\{k:P_{(k)}\leq q^{R^Z_+(0),\alpha}_k\frac{1\vee R^*_+(0)}{r+rV^Z_-(0)}\right\}.
\]
Then we construct the $\sigma$-algebra
\[
\sigma_{v,z',\cH_1}=\sigma(V^Z_+(0),V^Z_-(0),\bZ',\bW,\bZ_{\cH_1}),
\]
and we have
\begin{align*}
    \FDR=&\EE\left[\sum_{j\in\cH_0}\frac{\ind\{P_j\leq q^{R^Z_+(0),\alpha}_KE_j\}}{1\vee K}\right]\\
    =&\EE\left[\EE\left[\sum_{k=1}^m\frac{m_0}{k}\PP\left(K=k,P_J\leq q^{R^Z_+(0),\alpha}_K\frac{1\vee R^*_+(0)}{r+rV^Z_-(0)}\mid\sigma_{v,z',\cH_1}\right)\right]\right].
\end{align*}
Let $\cC_k(\boldsymbol{q})$ and $\cD_k(\boldsymbol{q}^{E,\alpha})$ be the same level sets as in the proof of Theorem~11, with
\[
\boldsymbol{q}^{E,\alpha}=(q^{R^Z_+(0),\alpha}_1\frac{1\vee R^*_+(0)}{r+rV^Z_-(0)},\dots,q^{R^Z_+(0),\alpha}_m\frac{1\vee R^*_+(0)}{r+rV^Z_-(0)}).
\]
Then
\begin{align*}
    \FDR\leq\EE\left[\EE\left[\sum_{k=1}^m\frac{m_0}{k}\PP\left(Q_J\leq q^{R^Z_+(0),\alpha}_k\frac{1\vee R^*_+(0)}{r+rV^Z_-(0)},\bP^J_{(\cdot)}\in\cC_k(\boldsymbol{q}^{E,\alpha})\mid\sigma_{v,z',\cH_1}\right)\right]\right].
\end{align*}
Because $\bZ'$ is independent of $(\bW,\bZ)$, the conditional probability factors as
\begin{align*}
    &\PP\left(Q_J\leq q^{R^Z_+(0),\alpha}_k\frac{1\vee R^*_+(0)}{r+rV^Z_-(0)},\bP^J_{(\cdot)}\in\cC_k(\boldsymbol{q}^{E,\alpha})\mid\sigma_{v,z',\cH_1}\right)\\
    =&\PP\left(\bP^J_{(\cdot)}\in\cC_k(\boldsymbol{q}^{E,\alpha})\mid Q_J\leq q^{R^Z_+(0),\alpha}_k\frac{1\vee R^*_+(0)}{r+rV^Z_-(0)},\sigma_{v,z',\cH_1}\right)\\
    &\times \PP\left(Q_J\leq q^{R^Z_+(0),\alpha}_k\frac{1\vee R^*_+(0)}{r+rV^Z_-(0)} \mid \sigma_{v,z',\cH_1}\right).
\end{align*}
The first factor is non-decreasing in $k$ by the dominated randomized $\PRDS_{\cH_0}$ of $\bP$ with respect to $\bQ$, and the second factor is bounded by
\[
\PP\left(Q_J\leq q^{R^Z_+(0),\alpha}_k\frac{1\vee R^*_+(0)}{r+rV^Z_-(0)} \mid \sigma_{v,z',\cH_1}\right)
\leq \frac{k\alpha}{m_0}\frac{V^Z_+(0)}{r+rV^Z_-(0)}\frac{R^*_+(0)}{1\vee R^Z_+(0)},
\]
and $R^*_+(0)\leq R^Z_+(0)$. Thus we have the FDR control.

We emphasize that the label-independence condition, while not required by the definition of sub-competition statistics, is essential for the above analysis. Without label-independence we cannot construct auxiliary labels $\bZ'$ that are independent of $\bZ$ and serve to couple $\bZ$ back to $\bL$. (The bridge variables $\bY$ of Section~\ref{sec:bridge} are inherently dependent on both $\bZ$ and $\bL$.) In the dependent case, only one of the two factors in the final decomposition can be controlled: conditioning on $\sigma_{v,z',\cH_1}$ preserves the monotonicity of the first (PRDS) factor, but the second factor
\[
\PP\Bigl(Q_J\leq q^{R^Z_+(0),\alpha}_k\frac{1\vee R^*_+(0)}{r+rV^Z_-(0)} \Bigm| \sigma_{v,z',\cH_1}\Bigr)
\]
now depends on the auxiliary variables through the conditioning algebra, and the super-uniformity bound used above may no longer follow. Hence this embedding argument does not establish FDR control for general sub-competition statistics.

This construction differs structurally from that in the proof of Theorem~10. For classical p-values, the PRDS property is intrinsic to the statistics themselves, and one constructs a deterministic correction: in the PRDS conditional probability, the argument before the conditioning bar is held fixed while the threshold after it is adjusted. For sub-competition statistics, by contrast, all favourable properties originate from $\bQ$, whose PRDS follows from the exchangeability of the random split; here one proceeds in the opposite direction, adjusting the argument before the bar via a random non-decreasing map $\bP=\phi(\bQ)$ while the conditioning event $Q_J\leq t$ remains. $\bQ$ and $\bP$ are determined by the exchangeability structure and the rejection region respectively, and their gap $\bP-\bQ$ is inherently random, which must be accounted for in the condition.

\section{Analysis of Compound p-Values under the Randomized BH Framework}
\label{app:compound}

\subsection{Violation of PRDS by Compound p-Values}

We consider the example of independent compound p-values from~\citet{barber2025falsediscoveryratecontrol}. The original construction includes a $\pi_0$ correction, which we retain as is. Let $k = 1/(2\alpha)$, assumed much smaller than $m$, and define the null set $\cH_0 = \{1, \dots, m - 3k + 2\} \subseteq [m] = \cH$. Let $\bP = (P_1, \dots, P_m)$ have the following distribution:

    The null $P_i$'s have the distribution
    \[
    P_1 \sim \operatorname{Unif}\left\{\frac{1.5}{m}, 1\right\}, \quad 
    P_2 = \frac{1}{m}, \quad 
    P_3 = \dots = P_{m-3k+2} = 1.
    \]
    whereas the non-null $P_i$'s are deterministic, with values
    \[
    P_{m-3k+3} = \dots = P_{m-k+1} = \frac{1}{m}, \quad 
    P_{m-k+2} = \dots = P_m = \frac{1.5}{m}.
    \]

The $P_i$'s are mutually independent. This construction satisfies the definition of compound p-values, since
\[
\sum_{i \in \cH_0} \PP\left\{P_i \leq \frac{1}{m}\right\} = \PP\left\{P_2 \leq \frac{1}{m}\right\} = 1,
\]
and
\[
\sum_{i \in \cH_0} \PP\left\{P_i \leq \frac{1.5}{m}\right\} = \PP\left\{P_1 \leq \frac{1.5}{m}\right\} + \PP\left\{P_2 \leq \frac{1.5}{m}\right\} = \frac{1}{2} + 1 = 1.5.
\]
In this example, only $P_1$ is random. When $P_1 = 1.5/m$, we have $\cR_{\alpha}(\bP) = \{1, 2\} \cup \cH_1$; when $P_1 = 1$, we have $\cR_{\alpha}(\bP) = \{2\} \cup \{m-3k+3, \dots, m-k+1\}$. Consequently, the BH procedure fails to control the FDR at level $\alpha$; in fact $\FDR \geq 7\alpha/6$. The failure occurs because the randomized $\PRDS_{\cH_0}$ condition is violated. To verify this, consider the non-decreasing set $C=\{\boldsymbol{p}\in\RR_*^{m-1}:p_{(3k-1)}\geq 2/m\}$ and two values of $t$:
\begin{align*}
    &\PP\left\{\bP^J_{(\cdot)}\in C\mid P_J\leq\frac{1}{m}\right\}\\
    =&\sum_{j\in\cH_0}\PP\left\{\bP^j_{(\cdot)}\in C,J=j\mid P_J\leq\frac{1}{m}\right\}\\
    =&\sum_{j\in\cH_0}\PP\left\{\bP^j_{(\cdot)}\in C\mid J=j\right\}\PP\left\{J=j\mid P_J\leq\frac{1}{m}\right\}\\
    =&\PP\left\{P_1\geq \frac{2}{m}\mid J=2\right\}=\frac{1}{2}
\end{align*}
while
\begin{align*}
    &\PP\left\{\bP^J_{(\cdot)}\in C\mid P_J\leq\frac{1.5}{m}\right\}\\
    =&\sum_{j\in\cH_0}\PP\left\{\bP^j_{(\cdot)}\in C,J=j\mid P_J\leq\frac{1.5}{m}\right\}\\
    =&\sum_{j\in\cH_0}\PP\left\{\bP^j_{(\cdot)}\in C\mid J=j\right\}\PP\left\{J=j\mid P_J\leq\frac{1.5}{m}\right\}\\
    =&\frac{1}{3}\PP\left\{\bP^1_{(\cdot)}\in C\mid J=1\right\}+\frac{2}{3}\PP\left\{P_1\geq\frac{2}{m}\mid J=2\right\}\\
    =&\frac{1}{3}\cdot0+\frac{2}{3}\cdot\frac{1}{2}=\frac{1}{3}<\frac{1}{2}.
\end{align*}
This example shows that independence alone does not imply the randomized $\PRDS_{\cH_0}$ condition when the null p-values lack identical marginal distributions. The key distinction between compound p-values (or randomized p-variables) and ordinary p-values lies precisely in this balance of marginals.

More generally, the FDR inflation caused by the failure of PRDS can be bounded as follows.

Suppose the coordinates of $\bP\in\RR_*^m$ are independent randomized super-uniform with respect to $\cH_0$. Define
\[
\cC_k=\left\{\boldsymbol{p}\in\RR_*^{m-1}:p_{(k)}>q^{m,\alpha}_{k+1},\dots,p_{(m-1)}>q^{m,\alpha}_{m}\right\},
\]
with $p_{(0)}=0$ and $\cC_0=\emptyset$, and let $\cD_k=\cC_k\setminus\cC_{k-1}$. Then the BH procedure $\cR^{\BH}_{\alpha}(\bP)$ satisfies
\begin{align*}
    &\FDR\leq\alpha\frac{m_0}{m}\sum_{k=1}^m\PP\left\{|\cR|=k\mid P_J\leq q^{m,\alpha}_k\right\}\\
    =&\alpha\frac{m_0}{m}\sum_{k=1}^m\PP\left\{\bP^J_{(\cdot)}\in\cD_k\mid P_J\leq q^{m,\alpha}_k\right\}\\
    =&\alpha\frac{m_0}{m}\sum_{k=1}^m\Bigl(\PP\left\{\bP^J_{(\cdot)}\in\cC_k\mid P_J\leq q^{m,\alpha}_k\right\}
      -\PP\left\{\bP^J_{(\cdot)}\in\cC_{k-1}\mid P_J\leq q^{m,\alpha}_k\right\}\Bigr)\\
    =&\alpha\frac{m_0}{m}+\alpha\frac{m_0}{m}\sum_{k=1}^m\Bigl(\PP\left\{\bP^J_{(\cdot)}\in\cC_{k-1}\mid P_J\leq q^{m,\alpha}_{k-1}\right\}
      -\PP\left\{\bP^J_{(\cdot)}\in\cC_{k-1}\mid P_J\leq q^{m,\alpha}_k\right\}\Bigr).
\end{align*}
Under independence of the coordinates, $P_j$ is independent of $\bP^j$ for each $j$, and $J$ is independent of $\bP$. Hence for any threshold $t$,
\[
\PP\left\{\bP^J_{(\cdot)}\in\cC_{k-1},J=j\mid P_J\leq t\right\}
=\PP\left\{\bP^j_{(\cdot)}\in\cC_{k-1}\right\}\,
\PP\left\{J=j\mid P_J\leq t\right\}.
\]
Applying this with $t=q^{m,\alpha}_{k-1}$ and $t=q^{m,\alpha}_k$, we obtain
\begin{align*}
    \PP\left\{\bP^J_{(\cdot)}\in\cC_{k-1}\mid P_J\leq q^{m,\alpha}_{k-1}\right\}
    &=\sum_{j\in\cH_0}\PP\left\{\bP^j_{(\cdot)}\in\cC_{k-1}\right\}\,
    \PP\left\{J=j\mid P_J\leq q^{m,\alpha}_{k-1}\right\},\\
    \PP\left\{\bP^J_{(\cdot)}\in\cC_{k-1}\mid P_J\leq q^{m,\alpha}_k\right\}
    &=\sum_{j\in\cH_0}\PP\left\{\bP^j_{(\cdot)}\in\cC_{k-1}\right\}\,
    \PP\left\{J=j\mid P_J\leq q^{m,\alpha}_k\right\}.
\end{align*}
Define
\[
A_{j,k}=\PP\left\{\bP^j_{(\cdot)}\in\cC_{k}\right\},\qquad
B_{j,k}=\PP\left\{J=j\mid P_J\leq q^{m,\alpha}_{k}\right\}.
\]
Substituting and applying summation by parts,
\begin{align*}
    &\sum_{k=1}^m\Bigl(\PP\left\{\bP^J_{(\cdot)}\in\cC_{k-1}\mid P_J\leq q^{m,\alpha}_{k-1}\right\}-\PP\left\{\bP^J_{(\cdot)}\in\cC_{k-1}\mid P_J\leq q^{m,\alpha}_k\right\}\Bigr)\\
    &\quad=\sum_{k=1}^m\sum_{j\in\cH_0}\bigl(A_{j,k-1}B_{j,k-1}-A_{j,k-1}B_{j,k}\bigr)\\
    &\quad=\sum_{j\in\cH_0}\Bigl(-A_{j,m-1}B_{j,m}+\sum_{k=1}^m(A_{j,k}-A_{j,k-1})B_{j,k}\Bigr).
\end{align*}
The terms $A_{j,k}$ and $B_{j,k}$ are determined by the marginal distributions of $\bP$. The resulting expression bounds the FDR in terms of these marginals. Evaluating the bound numerically is feasible in principle, though the underlying optimization over the joint distribution may be computationally demanding.

\subsection{Negative Dependence and Compound p-Values}

\begin{proof}[Proof of Lemma~5]
For $k\in[m]$, let $F_k(x)=\PP(X_k\le x)$.

Fix $A\subseteq[m]$ with $|A|=r$ and $x\in\RR$. Conditioning on $\bJ=\pi$ for a fixed permutation $\pi$ of $[m]$,
\[
\begin{aligned}
\PP\Bigl(\bigcap_{i\in A}\{Z_i\le x\}\Bigm| \bJ=\pi\Bigr)
&=\PP\Bigl(\bigcap_{i\in A}\{X_{\pi(i)}\le x\}\Bigr)\\
&\le \prod_{i\in A}\PP(X_{\pi(i)}\le x)
=\prod_{i\in A}F_{\pi(i)}(x),
\end{aligned}
\]
where the inequality follows from the WND property of $\bX$ applied to $\pi(A)$.

Averaging over $\bJ$ and using its independence from $\bX$,
\[
\PP\Bigl(\bigcap_{i\in A}\{Z_i\le x\}\Bigr)
\le \EE\Bigl[\prod_{i\in A}F_{\bJ_i}(x)\Bigr]
= \EE\Bigl[\prod_{k\in \bJ(A)}F_k(x)\Bigr].
\]
Since $\bJ(A)$ is a uniformly random $r$-subset of $[m]$,
\[
\EE\Bigl[\prod_{k\in \bJ(A)}F_k(x)\Bigr]
=
\frac{1}{\binom{m}{r}}
\sum_{\substack{B\subseteq[m]\\ |B|=r}}
\prod_{k\in B}F_k(x).
\]
By Maclaurin's inequality for the elementary symmetric functions of the nonnegative numbers $F_1(x),\dots,F_m(x)$,
\[
\frac{1}{\binom{m}{r}}
\sum_{\substack{B\subseteq[m]\\ |B|=r}}
\prod_{k\in B}F_k(x)
\le
\left(\frac{1}{m}\sum_{k=1}^m F_k(x)\right)^r.
\]
On the other hand, for each $i\in[m]$,
\[
\PP(Z_i\le x)
=
\sum_{k=1}^m \PP(\bJ_i=k)\PP(X_k\le x)
=
\frac{1}{m}\sum_{k=1}^m F_k(x).
\]
Thus
\[
\prod_{i\in A}\PP(Z_i\le x)
=
\left(\frac{1}{m}\sum_{k=1}^m F_k(x)\right)^r.
\]
Combining these bounds yields
\[
\PP\Bigl(\bigcap_{i\in A}\{Z_i\le x\}\Bigr)
\le
\prod_{i\in A}\PP(Z_i\le x).
\]
Since $A\subseteq[m]$ and $x\in\RR$ were arbitrary, $\bZ$ is WND.
\end{proof}

\begin{proof}[Proof of Theorem~13]
Construct $\bQ\in\RR_*^m$ with $Q_i=P_i$ for $i\in\cH_1$ and $Q_i=P_{J_i}$ for $i\in\cH_0$, where $\cH_1=\{m_0+1,\dots,m\}$, $\cH_0=[m_0]$, and $\bJ$ is an independent uniformly random permutation of $[m_0]$.

Since $\bQ$ differs from $\bP$ only by a permutation of the null coordinates, the two vectors share the same multiset of values, and the BH procedure depends on its input only through this multiset. Hence $|\cR^{\BH}_{\alpha}(\bQ)|=|\cR^{\BH}_{\alpha}(\bP)|$ and $|\cR^{\BH}_{\alpha}(\bQ)\cap\cH_0|=|\cR^{\BH}_{\alpha}(\bP)\cap\cH_0|$, so testing with $\bQ$ yields the same FDR as testing with $\bP$, and it suffices to control the FDR of $\bQ$.

The FDR-Linking theorem \citep{su2018} gives
\[
\FDR\leq\alpha+\alpha\int_{\alpha}^1\frac{\FDR_0(x)}{x^2}\rd x,
\]
where $\FDR_0(x)$ is the FDR of the BH procedure when all hypotheses are null and the target level is $x$. It remains to bound $\FDR_0(x)$ in the FDR-Linking formula, which we do via \citet{chi2024}. By Lemma~5, $\bQ_{\cH_0}$ is WND: indeed, $\bP_{\cH_0}$ is WND by the assumption of the theorem, and $\bQ_{\cH_0}$ is obtained from $\bP_{\cH_0}$ by an independent uniformly random permutation. Under WND, \citet{chi2024} yields
\[
\FDR_0(x)=\PP\{S_m(\bQ)\leq x\}\leq x+\sum_{k=2}^m\binom{m}{k}\Bigl(\frac{xk}{m}\Bigr)^k\leq x+2x^2+\frac{9}{2}x^3+\sum_{k=4}^{\infty}\frac{(\re x)^k}{\sqrt{2\pi k}},
\]
with $S_m(\boldsymbol{q})=\wedge_{k=1}^m mq_{(k)}/k$. Substituting this bound into the FDR-Linking formula and using $\FDR_0(x)\leq 1$ completes the proof.
\end{proof}

\subsection{Effect of Marginal Heterogeneity on the FDR Bound}

The analysis of the preceding subsection shows that under independence the $\FDR$ can be expressed as
\[
\FDR=\pi_0\alpha+\pi_0\alpha\sum_{k=1}^{m}\Delta_k,
\qquad
\Delta_k=\sum_{j\in\cH_0}A_{j,k-1}\bigl(B_{j,k-1}-B_{j,k}\bigr),
\]
where $A_{j,k}=\PP\{\bP^j_{(\cdot)}\in\cC_k\}$ and $B_{j,k}=\PP\{J=j\mid P_J\leq q^{m,\alpha}_k\}$ are as defined above. Since $\pi_0\alpha$ is the nominal level, the deviation from exact FDR control is governed by the $\Delta_k$ terms. To quantify how heterogeneity in the marginal distributions of the $\bP$ affects these terms, we decompose $A_{j,k}$ and $B_{j,k}$ around their averages over $\cH_0$. Define
\begin{gather*}
\bar A_k=\frac{1}{m_0}\sum_{j\in\cH_0}A_{j,k},
\qquad
d_{j,k}=A_{j,k}-\bar A_k,\\
\bar F_k=\bar F(q^{m,\alpha}_k)=\frac{1}{m_0}\sum_{j\in\cH_0}\PP\{P_j\leq q^{m,\alpha}_k\},\\
h_{j,k}=\PP\{P_j\leq q^{m,\alpha}_k\}-\bar F_k,
\qquad
V_k=\frac{1}{m_0}\sum_{j\in\cH_0}(h_{j,k})^2.
\end{gather*}
Using
\[
\PP\{J=j\mid P_J\leq q^{m,\alpha}_k\}=\frac{\PP\{P_j\leq q^{m,\alpha}_k\}}{m_0\,\PP\{P_J\leq q^{m,\alpha}_k\}},
\]
we have
\[
B_{j,k}=\frac{1}{m_0}\Bigl(1+\frac{h_{j,k}}{\bar F_k}\Bigr),
\qquad
B_{j,k-1}-B_{j,k}=\frac{1}{m_0}\Bigl(\frac{h_{j,k-1}}{\bar F_{k-1}}-\frac{h_{j,k}}{\bar F_k}\Bigr).
\]
Since $\sum_{j\in\cH_0}(B_{j,k-1}-B_{j,k})=0$, the mean $\bar A_{k-1}$ cancels in $\Delta_k$, giving
\begin{align*}
\Delta_k
&=\sum_{j\in\cH_0}(A_{j,k-1}-\bar A_{k-1})\bigl(B_{j,k-1}-B_{j,k}\bigr)
   =\sum_{j\in\cH_0}\frac{d_{j,k-1}}{m_0}\,
     \Bigl(\frac{h_{j,k-1}}{\bar F_{k-1}}-\frac{h_{j,k}}{\bar F_k}\Bigr).
\end{align*}
Applying the Cauchy--Schwarz inequality,
\begin{align*}
|\Delta_k|
&\leq \frac{1}{m_0}\,
   \sqrt{\sum_{j\in\cH_0}d_{j,k-1}^{\,2}}\,
   \sqrt{\sum_{j\in\cH_0}\Bigl(\frac{h_{j,k-1}}{\bar F_{k-1}}-\frac{h_{j,k}}{\bar F_k}\Bigr)^{2}}\\
&\leq \frac{1}{m_0}\,
   \sqrt{\frac{m_0}{4}}\,
   \sqrt{\sum_{j\in\cH_0}\Bigl(\frac{h_{j,k-1}}{\bar F_{k-1}}-\frac{h_{j,k}}{\bar F_k}\Bigr)^{2}}
   \qquad\bigl(\text{since }A_{j,k-1}\in[0,1]\Rightarrow d_{j,k-1}^{\,2}\leq\tfrac14\bigr)\\
&=\frac{1}{2}\,
   \sqrt{\frac{1}{m_0}\sum_{j\in\cH_0}
   \Bigl(\frac{h_{j,k-1}}{\bar F_{k-1}}-\frac{h_{j,k}}{\bar F_k}\Bigr)^{2}}.
\end{align*}
By the triangle inequality for the $\ell^2$-norm,
\[
\sqrt{\frac{1}{m_0}\sum_{j\in\cH_0}\Bigl(\frac{h_{j,k-1}}{\bar F_{k-1}}-\frac{h_{j,k}}{\bar F_k}\Bigr)^{2}}
\leq \frac{\sqrt{V_{k-1}}}{\bar F_{k-1}}+\frac{\sqrt{V_k}}{\bar F_k}.
\]
Hence
\[
|\Delta_k|\leq \frac{1}{2}\Bigl(\frac{\sqrt{V_{k-1}}}{\bar F_{k-1}}+\frac{\sqrt{V_k}}{\bar F_k}\Bigr).
\]
This bound shows that the deviation $|\Delta_k|$ is controlled by the coefficient of variation $\sqrt{V_k}/\bar F_k$ of the marginal exceedance probabilities $\PP\{P_j\leq q^{m,\alpha}_k\}$ across $j\in\cH_0$. When the null p-values have identical marginal distributions, $V_k=0$ for all $k$ and the bound collapses to $\Delta_k=0$, recovering exact $\FDR$ control at $\pi_0\alpha$. Heterogeneity in the marginals inflates $V_k$, and the resulting $\FDR$ deviation is proportional to the dispersion of these marginal probabilities.

If $\Delta_k$ were always nonnegative, then marginal heterogeneity would necessarily inflate the $\FDR$ above the nominal level $\pi_0\alpha$. This is not the case, however: the following example shows that $\Delta_k$ can be strictly negative, so that heterogeneity may instead render the procedure more conservative. With $m=m_0=2$ and
\[
P_1=\begin{cases}
1, & \text{with probability } 1-\alpha/3,\\[2pt]
\alpha/3, & \text{with probability } \alpha/3,
\end{cases}
\qquad
P_2=\begin{cases}
1, & \text{with probability } 1-2\alpha/3,\\[2pt]
2\alpha/3, & \text{with probability } 2\alpha/3,
\end{cases}
\]
a direct calculation using $A_{j,1}=\PP\{P_{3-j}\geq \alpha\}$ gives
\begin{align*}
\Delta_2
&=A_{1,1}B_{1,1}+A_{2,1}B_{2,1}-A_{1,1}B_{1,2}-A_{2,1}B_{2,2}\\[2pt]
&=\Bigl(1-\frac{2\alpha}{3}\Bigr)\frac{\alpha/3}{\alpha/3+0}
  +\Bigl(1-\frac{\alpha}{3}\Bigr)\frac{0}{\alpha/3+0}\\
&\qquad-\Bigl(1-\frac{2\alpha}{3}\Bigr)\frac{\alpha/3}{\alpha/3+2\alpha/3}
  -\Bigl(1-\frac{\alpha}{3}\Bigr)\frac{2\alpha/3}{\alpha/3+2\alpha/3}\\[2pt]
&=1-\frac{2\alpha}{3}-\frac{1}{3}+\frac{2\alpha}{9}-\frac{2}{3}+\frac{2}{9}\alpha
 =-\frac{2}{9}\,\alpha<0.
\end{align*}
Thus marginal heterogeneity does not systematically inflate the $\FDR$; depending on the configuration of the marginal distributions, the deviation may be positive or negative, and its magnitude in either direction is controlled by the coefficient of variation bound derived above.

Furthermore, the result shows that the deviation of the $\FDR$ from its nominal level can be bounded in terms of the coefficient of variation of the marginal exceedance probabilities, which is a measure of the heterogeneity of the null distributions. A sharper bound is obtained before applying the triangle inequality: the deviation is governed by
\[
\sqrt{\frac{1}{m_0}\sum_{j\in\cH_0}\Bigl(\frac{h_{j,k-1}}{\bar F_{k-1}}-\frac{h_{j,k}}{\bar F_k}\Bigr)^{2}},
\]
the $\ell^2$-norm of the change in the relative deviation $h_{j,k}/\bar F_k$ between adjacent thresholds. In particular, if $\PP\{P_j\leq t\}=c_j f(t)$ for positive constants $c_j$ and a common function $f$, then
\[
\frac{h_{j,k}}{\bar F_k}=\frac{(c_j-\bar c)\,f(q^{m,\alpha}_k)}{\bar c\,f(q^{m,\alpha}_k)}=\frac{c_j-\bar c}{\bar c},
\qquad \bar c=\frac{1}{m_0}\sum_{j\in\cH_0}c_j,
\]
which is independent of $k$. Hence $h_{j,k-1}/\bar F_{k-1}=h_{j,k}/\bar F_k$, so $\Delta_k=0$ and the $\FDR$ is controlled at exactly $\pi_0\alpha$ with no inflation. This is consistent with the analysis in the main text, where $\PP\{P_j\leq t\}=c_j f(t)$ is identified as the precise condition under which the selection weights $\PP\{J=j\mid P_J\leq t\}$ do not vary with $t$, so that FDR is preserved.

\section{Analysis of Integrated Tests under the Randomized BH Framework}
\label{asec:integration}

\subsection{A Counterexample to FDR Control under Integration}

We construct a counterexample to demonstrate that FDR control can be violated under integration. Consider two groups of competition statistics with parameters $r_1=r_2=1$. Each group contains $11$ hypotheses, of which the first $9$ are alternatives and the last $2$ are null. Under the alternative the label is $1$ and the score exceeds every null score. At target level $\alpha=0.1$, each group admits three null configurations that determine the competition outcome: (a) both null hypotheses carry label $1$; (b) only the higher-scoring null carries label $1$; and (c) the higher-scoring null carries label $0$, in which case no hypothesis is rejected. Since $r=1$, these configurations occur with probabilities $0.25$, $0.25$, and $0.5$, and yield $(11,2)$, $(10,1)$, and $(0,0)$ rejections, each pair being (total rejections, false rejections). Hence
\begin{align*}
    \FDR&=\EE\left[\frac{|(\cH^1_0\cap\cR^{\cp}_{1,0.1}(\bW^1,\bL^1))\cup(\cH^2_0\cap\cR^{\cp}_{1,0.1}(\bW^2,\bL^2))|}{1\vee|\cR^{\cp}_{1,0.1}(\bW^1,\bL^1)\cup\cR^{\cp}_{1,0.1}(\bW^2,\bL^2)|}\right]\\
    &=\frac{1}{4}\times\frac{1}{4}\times\frac{4}{22}
     +\frac{1}{4}\times\frac{1}{4}\times\frac{2}{20}
     +2\times\frac{1}{4}\times\frac{1}{4}\times\frac{3}{21}\\
    &\qquad+2\times\frac{1}{4}\times\frac{1}{2}\times\frac{2}{11}
     +2\times\frac{1}{4}\times\frac{1}{2}\times\frac{1}{10}
    =\frac{261}{2464}>0.1=\alpha.
\end{align*}


\subsection{The Proofs of Integration Tests}

We first record an equivalent characterization of the BH rejection count, which underlies the FDR decomposition below.
\begin{Lemma}
\label{slem:rejection}
Let $\cR^{\BH}_{\alpha}(\bP)$ denote the rejection set of the BH procedure at level $\alpha$ applied to $\bP\in\RR_*^{m}$, and let $K=|\cR^{\BH}_{\alpha}(\bP)|$ be the number of rejections. With
\[
K'=\sup\left\{k\in[m]:\sum_{j\in\cH}\ind\{P_j\leq q^{m,\alpha}_k\}\geq k\right\},
\]
we have $K'=K$.
\end{Lemma}

\begin{proof}[Proof of Theorem~14]
Let $m=m^1+m^2$ and $q_k=q^{m,\alpha}_k=\alpha k/m$. By Lemma~\ref{slem:rejection} applied to the concatenated vector $(\bP,\bQ)$,
\[
K=\sup\left\{k:\sum_{j_1\in\cH^1}\ind\{P_{j_1}\leq q_k\}+\sum_{j_2\in\cH^2}\ind\{Q_{j_2}\leq q_k\}\geq k\right\},
\]
and $K=|\cR^{\BH}_{\alpha}((\bP,\bQ))|$. Let $V_1=\sum_{j_1\in\cH^1_0}\ind\{P_{j_1}\leq q_K\}$ and $V_2=\sum_{j_2\in\cH^2_0}\ind\{Q_{j_2}\leq q_K\}$ denote the numbers of false discoveries in the two groups. Then the FDR can be expressed as
\begin{align*}
\FDR
  &=\EE\Biggl[\frac{\sum_{j_1\in\cH^1_0}\ind\{P_{j_1}\leq q^{m^1+m^2,\alpha}_K\}+\sum_{j_2\in\cH^2_0}\ind\{Q_{j_2}\leq q^{m^1+m^2,\alpha}_K\}}{1\vee K}\Biggr]\\
  &=\EE\Biggl[\frac{\sum_{j_1\in\cH^1_0}\ind\{P_{j_1}\leq q^{m^1+m^2,\alpha}_K\}}{1\vee K}
       +\frac{\sum_{j_2\in\cH^2_0}\ind\{Q_{j_2}\leq q^{m^1+m^2,\alpha}_K\}}{1\vee K}\Biggr]\\
  &=\sum_{k=1}^{m^1+m^2}\frac{1}{k}\sum_{j_1\in\cH^1_0}\PP\bigl\{P_{j_1}\leq q^{m^1+m^2,\alpha}_k,\,K=k\bigr\}\\
  &\quad+\sum_{k=1}^{m^1+m^2}\frac{1}{k}\sum_{j_2\in\cH^2_0}\PP\bigl\{Q_{j_2}\leq q^{m^1+m^2,\alpha}_k,\,K=k\bigr\}.
\end{align*}
Let $J_1$ and $J_2$ be independent of $(\bP,\bQ)$ and uniformly distributed over $\cH^1_0$ and $\cH^2_0$, respectively. For any non-decreasing set $D\in\RR_*^{m-1}$ and any $t$,
\begin{align*}
&\PP\bigl\{(\bP^{J_1},\bQ)_{(\cdot)}\in D,\ P_{J_1}\leq t \mid \bQ\bigr\}\\
&\qquad=\frac{1}{1\vee|\cH^1_0|}\sum_{j_1\in\cH^1_0}\PP\bigl\{(\bP^{j_1},\bQ)_{(\cdot)}\in D,\ P_{j_1}\leq t \mid \bQ\bigr\},\\
&\PP\bigl\{(\bP,\bQ^{J_2})_{(\cdot)}\in D,\ Q_{J_2}\leq t \mid \bP\bigr\}\\
&\qquad=\frac{1}{1\vee|\cH^2_0|}\sum_{j_2\in\cH^2_0}\PP\bigl\{(\bP,\bQ^{j_2})_{(\cdot)}\in D,\ Q_{j_2}\leq t \mid \bP\bigr\}.
\end{align*}
where $(\bP^{j_1},\bQ)_{(\cdot)}$ and $(\bP,\bQ^{j_2})_{(\cdot)}$ denote the order statistics of the combined vectors obtained by removing $P_{j_1}$ from $\bP$ or $Q_{j_2}$ from $\bQ$ respectively. Define
\[
\cC_k=\bigl\{\boldsymbol{u}\in\RR_*^{m-1}:u_{(k)}>q_{k+1},\dots,u_{(m-1)}>q_m\bigr\}
\]
with $u_{(0)}=0$ and $\cC_0=\emptyset$, $\cC_m=\RR_*^{m-1}$, and let $\cD_k=\cC_k\setminus\cC_{k-1}$. Then by the independence of $\bP$ and $\bQ$,
\[
\PP\{P_{j_1}\leq q_k,\ K=k\}=\EE_{\bQ}\bigl[\PP\bigl\{P_{j_1}\leq q_k,\ (\bP^{j_1},\bQ)_{(\cdot)}\in\cD_k\mid\bQ\bigr\}\bigr],
\]
\[
\PP\{Q_{j_2}\leq q_k,\ K=k\}=\EE_{\bP}\bigl[\PP\bigl\{Q_{j_2}\leq q_k,\ (\bP,\bQ^{j_2})_{(\cdot)}\in\cD_k\mid\bP\bigr\}\bigr].
\]
Now set
\[
\cC^1_k(\boldsymbol{q})=\bigl\{\boldsymbol{p}\in\RR_*^{m^1-1}:(\boldsymbol{p},\boldsymbol{q})\in\cC_k\bigr\},\qquad
\cC^2_k(\boldsymbol{p})=\bigl\{\boldsymbol{q}\in\RR_*^{m^2-1}:(\boldsymbol{p},\boldsymbol{q})\in\cC_k\bigr\},
\]
\[
\cD^1_k(\boldsymbol{q})=\bigl\{\boldsymbol{p}\in\RR_*^{m^1-1}:(\boldsymbol{p},\boldsymbol{q})\in\cD_k\bigr\},\qquad
\cD^2_k(\boldsymbol{p})=\bigl\{\boldsymbol{q}\in\RR_*^{m^2-1}:(\boldsymbol{p},\boldsymbol{q})\in\cD_k\bigr\}.
\]
For any given $\boldsymbol{p},\boldsymbol{q}$ and any $k\in[m]$, we have $\cD^1_k(\boldsymbol{q})=\cC^1_k(\boldsymbol{q})\setminus\cC^1_{k-1}(\boldsymbol{q})$ and $\cD^2_k(\boldsymbol{p})=\cC^2_k(\boldsymbol{p})\setminus\cC^2_{k-1}(\boldsymbol{p})$. Moreover $\cC^1_k(\boldsymbol{q})$ and $\cC^2_k(\boldsymbol{p})$ are non-decreasing sets. For brevity, set
\[
a_{n_1,n_2}(\bQ)=\PP\bigl\{P_{J_1}\leq q_{n_1},
    \bP^{J_1}_{(\cdot)}\in\cC^1_{n_2}(\bQ)\mid\bQ\bigr\},
\]
\[
b_{n_1,n_2}(\bP)=\PP\bigl\{Q_{J_2}\leq q_{n_1},
    \bQ^{J_2}_{(\cdot)}\in\cC^2_{n_2}(\bP)\mid\bP\bigr\},
\]
with $a_{n,0}(\bQ)=b_{n,0}(\bP)=0$. Then
\begin{align*}
\FDR
  &=\sum_{k=1}^{m}\frac{1}{k}
    \Bigl(\sum_{j_1\in\cH^1_0}\PP\bigl\{P_{j_1}\leq q_k,\;K=k\bigr\}
       +\sum_{j_2\in\cH^2_0}\PP\bigl\{Q_{j_2}\leq q_k,\;K=k\bigr\}\Bigr)\\
  &=\sum_{k=1}^{m}\frac{1}{k}
    \Bigl(|\cH^1_0|\,\EE_{\bQ}\bigl[\PP\bigl\{P_{J_1}\leq q_k,\ \bP^{J_1}_{(\cdot)}\in\cD^1_k(\bQ)\mid\bQ\bigr\}\bigr]\\
  &\qquad+|\cH^2_0|\,\EE_{\bP}\bigl[\PP\bigl\{Q_{J_2}\leq q_k,\ \bQ^{J_2}_{(\cdot)}\in\cD^2_k(\bP)\mid\bP\bigr\}\bigr]\Bigr)\\
  &=\EE_{\bQ}\Bigg[\sum_{k=1}^{m}\frac{|\cH^1_0|}{k}\bigl(a_{k,k}(\bQ)-a_{k,k-1}(\bQ)\bigr)\Bigg]
   +\EE_{\bP}\Bigg[\sum_{k=1}^{m}\frac{|\cH^2_0|}{k}\bigl(b_{k,k}(\bP)-b_{k,k-1}(\bP)\bigr)\Bigg]\\
  &\leq\alpha\frac{|\cH^1_0|+|\cH^2_0|}{m}.
\end{align*}
\end{proof}

\begin{proof}[Proof of Theorem~15]
The theorem has two parts, one for each algorithm. We prove only the part for Algorithm~2. The proof of the first part is similar and is omitted.
Let $m=m^1+m^2$ and $q_k=q^{m,\alpha}_k=\alpha k/m$. Define
\[
Y_j=\frac{1+\sum_{i=1}^{m^2}\ind\{W_i\geq W_j,\ L_i=0\}}{1+\sum_{i\in\cH^2_0}\ind\{L_i=0\}},\qquad
Z_j=\frac{1+\sum_{i\in\cH^2_0}\ind\{W_i\geq W_j,\ L_i=0\}}{1+\sum_{i\in\cH^2_0}\ind\{L_i=0\}},
\]
\[
E_j=\frac{m^2}{r+r\sum_{i\in\cH^2_0}\ind\{L_i=0\}},\qquad \text{for all } j \text{ with } L_j=1,
\]
and set $Y_j=Z_j=+\infty$ and $E_j=0$ for all $j$ with $L_j=0$. For any $k$, the equivalence
\[
X_j\leq q_k,\ L_j=1\;\Longleftrightarrow\;Y_j\leq q_kE_j
\]
holds. Then
\begin{align*}
\FDR
  &=\EE\left[\frac{|\cR_1\cap\cH^1_0|+|\cR_2\cap\cH^2_0|}{1\vee(|\cR_1|+|\cR_2|)}\right]\\
  &=\EE\Biggl[\EE\left[\frac{\sum_{j\in\cH^1_0}\ind\{P_j\leq q_K\}}{1\vee K}\mid\bW,\bL\right]\\
  &\qquad+\EE\left[\frac{\sum_{j\in\cH^2_0}\ind\{X_j\leq q_K,\ L_j=1\}}{1\vee K}\mid\bP\right]\Biggr].
\end{align*}

We bound the two conditional expectations separately. Since $\sigma\{\bX,\bL\}\supseteq\sigma\{\bW,\bL\}$, we may take $X_j=+\infty$ whenever $L_j=0$. From the proof of Theorem~14,
\[
\EE\left[\frac{\sum_{j\in\cH^1_0}\ind\{P_j\leq q_K\}}{1\vee K}\mid\bW,\bL\right]\leq\alpha\frac{|\cH^1_0|}{m},
\]
and the bound does not depend on $\bW$ or $\bL$.

For the second conditional expectation, we similarly take $Y_j=Z_j=+\infty$ whenever $L_j=0$ and consider the $\sigma$-algebra $\sigma_{\cH^2_1}=\sigma\{V_+(0),V_-(0),\bW,\bL_{\cH^2_1}\}$. Since $\bP$ is independent of $(\bW,\bL)$, the conditional controlled randomized $\PRDS_{\cH^2_0}$ condition holds with respect to $\bZ$, that is, for any non-decreasing set $C\subseteq\RR^{m^2-1}$ the conditional probability $\PP\{\bY^J_{(\cdot)}\in C\mid Z_J\leq t,\sigma_{\cH^2_1},\bP\}$ is non-decreasing in $t$. By Theorem~11, it remains to verify that, given $Z_J\leq q_k$, $\sigma_{\cH^2_1}$, and $\bP$, the event $\{K=k\}$ is a non-decreasing set with respect to $\bY^J_{(\cdot)}$. The proof of Theorem~14 shows that this condition holds as well. Hence
\[
\EE\left[\frac{\sum_{j\in\cH^2_0}\ind\{Y_j\leq q_KE_j\}}{1\vee K}\mid\bP,\sigma_{\cH^2_1}\right]\leq\alpha\frac{m^2}{m}\frac{V_+(0)}{r+rV_-(0)}.
\]
Combining the two bounds gives
\[
\FDR\leq\alpha\,\EE\left[\frac{|\cH^1_0|}{m}+\frac{m^2}{m}\frac{V_+(0)}{r+rV_-(0)}\right]\leq\alpha.
\]
\end{proof}

\begin{proof}[Proof of Theorem~16]
Let $m=m^1+m^2$ and $q_k=q^{m,\alpha}_k=\alpha k/m$. Rewrite the FDR as
\[
\FDR=\EE\Biggl[\sum_{g=1}^2
\frac{1+\sum_{j\in\cR^g_-(0)}\ind\bigl\{\frac{1+R^g_-(W_j)}{m^g c(r_g)}\leq q_K\bigr\}}{K}
\cdot
\frac{\sum_{j\in\cH^g_0}\ind\{X^g_j\leq q_K\}}
{1+\sum_{j\in\cR^g_-(0)}\ind\bigl\{\frac{1+R^g_-(W_j)}{m^g c(r_g)}\leq q_K\bigr\}}\Biggr].
\]

Define the minimal negative statistic
\[
W^g_{\mathrm{negmin}} = \min\Bigl\{ W^g_j : j\in\cR^g_-(0),\ \frac{1+R^g_-(W^g_j)}{m^g c(r_g)} \leq q_K \Bigr\}.
\]
For any $W^g_j < W^g_{\mathrm{negmin}}$ with $j\in\cR^g_-(0)$,
\[
\frac{1+R^g_-(W^g_j)}{m^g c(r_g)} \ge \frac{2+R^g_-(W^g_{\mathrm{negmin}})}{m^g c(r_g)} > q_K.
\]
Set
\[
W^g_{\mathrm{negmin}+1} = \max\{ W^g_j : j\in\cR^g_-(0),\ W^g_j < W^g_{\mathrm{negmin}} \}.
\]
Then
\[
X^g_j \leq q_K \quad\Longleftrightarrow\quad W^g_j > W^g_{\mathrm{negmin}+1}\ \text{and}\ L^g_j=1,
\]
and consequently
\begin{align*}
\frac{1+\sum_{j\in\cR^g_-(0)}\ind\bigl\{\frac{1+R^g_-(W_j)}{m^g c(r_g)}\leq q_K\bigr\}}{K}
&=\frac{1+R^g_-(W^g_{\mathrm{negmin}})}{K}\\
&=\frac{m^g c(r_g)}{K}\sup_{j\in\cH^g_0}\{X^g_j:X^g_j\leq q_K\}.
\end{align*}
Therefore,
\begin{align*}
\FDR
&\leq\alpha\sum_{g=1}^2\frac{m^g c(r_g)}{m}
  \,\EE\Biggl[\frac{\sum_{j\in\cH^g_0}\ind\{X^g_j\leq q_K\}}
  {1+\sum_{j\in\cR^g_-(0)}\ind\bigl\{\frac{1+R^g_-(W_j)}{m^g c(r_g)}\leq q_K\bigr\}}\Biggr]\\
&=\alpha\sum_{g=1}^2\frac{m^g c(r_g)}{m}
  \,\EE\Biggl[\frac{\sum_{j\in\cH^g_0}\ind\{W^g_j<W^g_{\mathrm{negmin}+1},\ L^g_j=1\}}
  {1+\sum_{j\in\cH^g}\ind\{W^g_j<W^g_{\mathrm{negmin}+1},\ L^g_j=0\}}\Biggr],
\end{align*}
where $W^g_{\mathrm{negmin}+1}$ depends on $K$, and $K$ depends on both groups, so direct computation is infeasible. By the properties of competition statistics, we may replace $W^g_{\mathrm{negmin}+1}$ with an arbitrary threshold $t$ and take the supremum over $t$ to obtain a bound independent of the inter-group dependence structure. This motivates the definition
\[
c(r_g)\,\EE\Biggl[\sup_{t\in\RR}\frac{\sum_{j\in\cH^g_0}\ind\{W^g_j<t,\ L^g_j=1\}}{1+\sum_{j\in\cH^g}\ind\{W^g_j<t,\ L^g_j=0\}}\Biggr]\leq1.
\]
Finally,
\begin{align*}
&\EE\Biggl[\sup_{t\in\RR}\frac{\sum_{j\in\cH^g_0}\ind\{W^g_j<t,\ L^g_j=1\}}
{1+\sum_{j\in\cH^g}\ind\{W^g_j<t,\ L^g_j=0\}}\Biggr]\\
&\qquad\leq\EE\Biggl[\sup_{t\in\RR}\frac{\sum_{j\in\cH^g_0}\ind\{W^g_j<t,\ L^g_j=1\}}
{1+\sum_{j\in\cH^g_0}\ind\{W^g_j<t,\ L^g_j=0\}}\Biggr]\\
&=\EE\Bigl[\sup_{k\in\{0\}\cup[|\cH^g_0|]}\frac{S_k}{1+k-S_k}\Bigr]
\leq\EE\Bigl[\sup_{k\in\mathbb{N}}\frac{S_k}{1+k-S_k}\Bigr],
\end{align*}
where $S_k$ follows a binomial distribution with parameters $k$ and $r_g/(1+r_g)$. The bound is independent of the statistics. To maximize the number of rejections, we should choose $c(r_g)$ as large as possible, so we take
\[
c(r)=\Bigl(\EE\Bigl[\sup_{k\in\mathbb{N}}\frac{S_k}{1+k-S_k}\Bigr]\Bigr)^{-1},
\]
where $S_k\sim\mathrm{Binomial}(k, r/(1+r))$. When $r=1$,
\[
\EE\Bigl[\sup_{k\in\mathbb{N}}\frac{S_k}{1+k-S_k}\Bigr]\leq 1.9227.
\]
This numerical bound is evaluated in \citet{10.1111/rssb.12430} and \citet{eugene2019}.
\end{proof}

\subsection{More Analysis of the Inflation of FDR under Integration}
First, without loss of generality, consider competition statistics with parameter $r=1$, and fix $\bW$ and $\bL_{\cH_1}$, thereby determining the ranking of the scores. By interleaving null hypotheses with alternative hypotheses carrying labels $0$ and $1$, we can achieve
\[
\frac{1+R_-(T)}{1\vee R_+(T)} = \alpha.
\]
Consider the following example, ordered by decreasing score. It begins with $10$ alternative hypotheses with label $1$, followed by a loop. In each iteration of the loop, we place one alternative hypothesis with label $0$, then one null hypothesis, and finally nine alternative hypotheses with label $1$. Setting $\alpha=0.1$ makes the equality above hold. The resulting sequence of labels is shown below, where $x$ denotes a null hypothesis whose label is left undetermined:
\[
1,1,1,1,1,1,1,1,1,1,\underbrace{0,x,1,1,1,1,1,1,1,1,1}_{\text{loop}},\dots.
\]
If $x=1$, the procedure proceeds to the next iteration. If $x=0$, the procedure rejects the most recent alternative hypothesis with label $1$ that precedes this null hypothesis, thereby attaining the equality. By the proof of FDR control for the competition procedure, it suffices to consider
\[
\EE\left[\frac{V^1_+(T^1)+V^2_+(T^2)}{1+V^1_-(T^1)+V^2_-(T^2)}\right],
\]
where $T^1$ and $T^2$ are the thresholds from the competition procedure. Applying the flip-one method,
\begin{align*}
    &\EE\left[\frac{V^1_+(T^1)+V^2_+(T^2)}{1+V^1_-(T^1)+V^2_-(T^2)}\right]
    =\EE\left[\frac{\sum_{g=1}^2\sum_{i\in\cH^g_0}\ind\{W^g_i\geq T^g,L^g_i=1\}}{1+\sum_{h=1}^2\sum_{j\in\cH^h_0}\ind\{W^h_j\geq T^h,L^h_j=0\}}\right]\\
    &=\EE\left[\sum_{g=1}^2\sum_{i\in\cH^g_0}\EE\left[\frac{\ind\{W^g_i\geq T^g,L^g_i=1\}}{1+V^1_-(T^1)+V^2_-(T^2)}\mid \bW^{3-g},\bL^{3-g},\bW^g,\bL^g_{-i}\right]\right]\\
    &=\sum_{g=1}^2\sum_{i\in\cH^g_0}\EE\Bigg[\EE\Bigg[\frac{\ind\{W^g_i\geq T^g_i,L^g_i=1\}}{1+V^g_-(T^g_i)+V^{3-g}_-(T^{3-g})}\mid \sigma_i^g\Bigg]\Bigg]
\end{align*}
where $\sigma^g_i=\sigma\{\bW^{3-g},\bL^{3-g},\bW^g,\bL^g_{-i}\}$, $V^{g,i}_-(t)=\sum_{j\in\cH^g_0\setminus\{i\}}\ind\{W^g_j\geq t,L^g_j=0\}$, and
\[
T^g_i=\inf\left\{t\in\{W^g_j:L_j=1\}:\frac{V^{g,i}_-(t)+1}{1\vee\left(V^{g,i}_-(t)+\ind\{W^g_i\geq t\}\right)}\leq\alpha\right\}.
\]
By the definition of competition statistics,
\begin{align*}
    &\EE\left[\frac{\ind\{W^g_i\geq T^g,L^g_i=1\}}{1+V^g_-(T^g)+V^{3-g}_-(T^{3-g})}\right]\\
    &=\EE\left[\frac{\ind\{W^g_i\geq T^g_i,L^g_i=0\}}{1+V^{g,i}_-(T^g_i)+V^{3-g}_-(T^{3-g})}\right]
\end{align*}
and by Lemma~\ref{slem:biovervariable},
\begin{align*}
    &\EE\left[\frac{V^1_+(T^1)+V^2_+(T^2)}{1+V^1_-(T^1)+V^2_-(T^2)}\right]\\
    &=\sum_{g=1}^2\EE\left[\frac{\sum_{j\in\cH^g_0}\ind\{W^g_j\geq T^g_j,L^g_j=0\}}{1\vee\left(\sum_{j\in\cH^g_0}\ind\{W^g_j\geq T^g_j,L^g_j=0\}+V^{3-g}_-(T^{3-g})\right)}\right].
\end{align*}
In this example, every flip decreases the threshold and increases the number of rejections, except when it involves the last null hypothesis. Consequently,
\[
\sum_{g=1}^2\frac{\sum_{j\in\cH^g_0}\ind\{W^g_j\geq T^g_j,L^g_j=0\}}{1\vee\left(\sum_{j\in\cH^g_0}\ind\{W^g_j\geq T^g_j,L^g_j=0\}+V^{3-g}_-(T^{3-g})\right)}>1.
\]
Thus, the constructed example can lead to a loss of FDR control.

Second, we use the proof of Theorem~15 to analyse why the counterexample occurs, that is, why direct integration of multiple testing results fails. Stated without proof, if in Algorithm~2 the statistic $X_j$ and the BH threshold are replaced as follows,
\begin{align*}
    &X_j\leftarrow\frac{r+r\sum_{i=1}^{m^2}\ind\{W_i\geq W_j,\;L_i=0\}}{m^2}\\&\quad\Longrightarrow\quad
X_j\leftarrow\frac{r+r\sum_{i=1}^{m^2}\ind\{W_i\geq W_j,\;L_i=0\}}{1\vee\sum_{i=1}^{m^2}\ind\{L_i=1\}},
\end{align*}

and
\[
q^{m^1+m^2,\alpha}_k \;\Longrightarrow\; q^{m^1+R_+(0),\alpha}_k,
\]
where $R_+(0)=\sum_{i=1}^{m^2}\ind\{L_i=1\}$, then
\[
\FDR\leq\alpha\,\EE\left[\frac{|\cH^1_0|}{m^1+R_+(0)}+\frac{R_+(0)}{m^1+R_+(0)}\frac{V_+(0)}{r+rV_-(0)}\right].
\]
For a fixed multiple testing problem, $R_+(0)$ and $V_+(0)$ are strongly positively correlated. Both factors are increasing functions of their arguments: $t\mapsto t/(m^1+t)$ is increasing in $t$, and since $V_+(0)+V_-(0)=|\cH^2_0|$, the map $V_+(0)\mapsto V_+(0)/(r+rV_-(0))$ is increasing as well. Hence, by the rearrangement inequality,
\[
\EE\left[\frac{R_+(0)}{m^1+R_+(0)}\frac{V_+(0)}{r+rV_-(0)}\right]>\EE\left[\frac{R_+(0)}{m^1+R_+(0)}\right]\EE\left[\frac{V_+(0)}{r+rV_-(0)}\right],
\]
and, when $m^2$ is sufficiently large,
\[
\EE\left[\frac{V_+(0)}{r+rV_-(0)}\right]\to1.
\]
This analysis suggests that the loss of global FDR control when integrating multiple testing results may be understood through the rearrangement inequality: results with higher FDP tend to have more rejections, occupy a larger share of the global FDR, and thereby push the overall FDR above the target level.

In brief, consistent with both the constructed example and the rearrangement inequality, if each group of tests tends to reject more hypotheses when its local FDP is large, then the global FDR may exceed the target level.

\subsection{More Analysis for Theorem~16}

First, consider why a correction factor $c(r)$ is needed for inter-group dependence. For competition statistics $(\bW^1,\bL^1)$ and $(\bW^2,\bL^2)$, Algorithm~1 yields $\bX^1$, $\bX^2$, and the threshold $K$. The FDR is
\begin{align*}
    \FDR = \EE\left[ \frac{1}{1\vee K} \sum_{g=1}^2 \sum_{j\in\cH^g_0} \ind\bigl\{ X^g_j \leq q^{m^1+m^2,\alpha}_K \bigr\} \right].
\end{align*}

Under inter-group independence, the integration test remains valid because the conditional probability
\[
\PP\bigl\{K\leq k \mid X^g_j\leq t,\; \bX^{3-g},\; \bL^1,\; \bL^2\bigr\}
\]
is non-decreasing in $t$ for any $j\in\cH^g_0$; as the event $\{K\leq k\}$ is non-decreasing in the p-value ordering, this conditional probability coincides with the marginal probability
\[
\PP\bigl\{K\leq k \mid X^g_j\leq t,\; \bL^g\bigr\}.
\]
Inter-group dependence, however, breaks this monotonicity and invalidates the integration test. To overcome this, we must apply the correction factor.

Second, does the correction factor $c(r)$ restore the conditional PRDS property? It does not. The correction factor $c(r)$ is designed to adapt to the special dependence structure of competition statistics, which is not PRDS. If this correction factor is applied to the pBH part of the integrated test, the FDR remains uncontrolled. Hence $c(r)$ addresses the dependence structure specific to competition statistics without restoring conditional PRDS.

Third, what property of the randomized p-variable generated by competition statistics allows FDR control under arbitrary dependence with the correction factor $c(r)$? From the proof of Theorem~16, if $\bP$ are continuous p-values (generated p-values, such as conformal p-values, can be made continuous by adding independent uniform noise), we construct a series of random variables $U_k\in[0,1]$ satisfying
\[
\sum_{j\in\cH_0}\ind\{P_j\leq U_k\} + \lfloor mU_k \rfloor = k.
\]
In other words, $U_k$ solves the equation
\[
\sum_{j\in\cH_0}\ind\{P_j\leq t\} + \lfloor mt \rfloor = k.
\]
The solution is unique in terms of $\lfloor mt \rfloor$ and always exists when $k\leq m$, because the left-hand side is a step function that jumps by one as $t$ increases from $0$ to $1$. A sufficient condition for FDR control with the correction factor $c(r)$ is that $k-\lfloor mU_k \rfloor$ is stochastically dominated by a binomial distribution with parameters $k$ and $r/(1+r)$.

\subsection{Integrated Test with Different Weights}

We present the following algorithms. They are almost identical to the algorithms given in the main text, except that the value $m$ for the integrated competition testing is replaced by $R_-(0)$. The details are as follows.

\begin{algorithm}[htbp]
\caption{CP-CP integration with $R_-(0)$}
\LinesNotNumbered
\label{salgo:int_cp_cp}
\emph{Input:} $(\bW^1,\bL^1)\in\RR_*^{m^1}\times\{0,1\}^{m^1}$, $(\bW^2,\bL^2)\in\RR_*^{m^2}\times\{0,1\}^{m^2}$, $r_1,r_2>0$, $\alpha\in(0,1)$. \emph{Output:} Rejection sets $\cR^1,\cR^2$.
\begin{enumerate}
  \item For $g=1,2$, set $R^g_-(0)=\sum_{i=1}^{m^g}\ind\{L^g_i=0\}$, and compute for each $j=1,\dots,m^g$:
        \[
          X^g_j=\frac{r_g+r_g\sum_{i=1}^{m^g}\ind\{W^g_i\geq W^g_j,\,L^g_i=0\}}{1\vee R^g_-(0)}.
        \]
  \item Find the largest $K\in\{1,\dots,R^1_-(0)+R^2_-(0)\}$ such that
    \[
    \sum_{g=1}^2\sum_{j=1}^{m^g}\ind\{X^g_j\leq q^{R^1_-(0)+R^2_-(0),\alpha}_K,\,L^g_j=1\}\geq K.
    \]
  \item Return $\cR^g=\{j: X^g_j\leq q^{R^1_-(0)+R^2_-(0),\alpha}_K,\,L^g_j=1\}$, $g=1,2$.
\end{enumerate}
\end{algorithm}

\begin{algorithm}[htbp]
\caption{CP-P integration with $R_-(0)$}
\LinesNotNumbered
\label{salgo:int_p_cp}
\emph{Input:} $\bP\in\RR_*^{m^1}$, $(\bW,\bL)\in\RR_*^{m^2}\times\{0,1\}^{m^2}$, $r>0$, $\alpha\in(0,1)$. \emph{Output:} Rejection sets $\cR^1,\cR^2$.
\begin{enumerate}
  \item Set $R_-(0)=\sum_{i=1}^{m^2}\ind\{L_i=0\}$, and compute for each $j=1,\dots,m^2$:
  \[
  X_j=\frac{r+r\sum_{i=1}^{m^2}\ind\{W_i\geq W_j,\,L_i=0\}}{1\vee R_-(0)}.
  \]      
  \item Find the largest $K\in\{1,\dots,m^1+m^2\}$ such that
    \[
    \sum_{j=1}^{m^1}\ind\{P_j\leq q^{m^1+R_-(0),\alpha}_K\}+\sum_{j=1}^{m^2}\ind\{X_j\leq q^{m^1+R_-(0),\alpha}_K,\,L_j=1\}\geq K.
    \]
  \item Return $\cR^1=\{j: P_j\leq q^{m^1+R_-(0),\alpha}_K\}$, $\cR^2=\{j: X_j\leq q^{m^1+R_-(0),\alpha}_K,\,L_j=1\}$.
\end{enumerate}
\end{algorithm}
Similarly, we have the following theorem to show they can control the FDR at the target level.
\begin{Thm}[Integration controls global FDR]
\label{sthm:fdr_int}
The integrated procedures described in
Algorithms~\ref{salgo:int_cp_cp}--\ref{salgo:int_p_cp}
control the global FDR under the following conditions:
\begin{enumerate}
\renewcommand{\labelenumi}{(\roman{enumi})}
  \item For Algorithm~\ref{salgo:int_cp_cp}:
        $(\bW^1,\bL^1)$ and $(\bW^2,\bL^2)$ are independent and
        are competition statistics with respect to $r_1,\cH^1_0\subseteq\cH^1$
        and $r_2,\cH^2_0\subseteq\cH^2$, respectively;
  \item For Algorithm~\ref{salgo:int_p_cp}:
        $\bP$ and $(\bW,\bL)$ are independent, $\bP$ are randomized BH-valid
        with respect to $\cH^1_0\subseteq\cH^1$, and $(\bW,\bL)$ are
        competition statistics with respect to $r,\cH^2_0\subseteq\cH^2$.
\end{enumerate}
Under either condition, for any $\alpha\in(0,1)$,
\[
\FDR = \EE\left[
\frac{|(\cR^1\cap\cH^1_0)\cup(\cR^2\cap\cH^2_0)|}
{1\vee|\cR^1\cup\cR^2|}
\right] \leq \alpha.
\]
\end{Thm}
We present the proof directly.
\begin{proof}[Proof of Theorem~\ref{sthm:fdr_int}]
As before, we focus only on Algorithm~\ref{salgo:int_p_cp}. From the proof of Theorem~15, it follows directly that
\[
\EE\left[\frac{\sum_{j\in\cH^1_0}\ind\{P_j\leq q_K^{m^1+R_-(0),\alpha}\}}{1\vee K}\;\Bigg|\;\bW,\bL\right]\leq\alpha\frac{|\cH^1_0|}{m^1+R_-(0)}
\]
and, with $\sigma_{\cH^2_1}=\sigma\{V_+(0),V_-(0),\bW,\bL_{\cH^2_1}\}$,
\[
\EE\left[\frac{\sum_{j\in\cH^2_0}\ind\{Y_j\leq q_K^{m^1+R_-(0),\alpha}E_j,\,L_j=1\}}{1\vee K}\;\Bigg|\;\bP,\sigma_{\cH^2_1}\right]\leq\alpha\frac{R_-(0)}{m^1+R_-(0)}\frac{V_+(0)}{r+rV_-(0)}.
\]
Given $\bL_{\cH^2_1}$, we have $R_-(0)=C-V_+(0)$ for a constant $C$, so that $R_-(0)/(m^1+R_-(0))$ is non-increasing in $V_+(0)$ while $V_+(0)/(r+rV_-(0))$ is non-decreasing in $V_+(0)$. By the rearrangement inequality,
\[
\EE\left[\frac{R_-(0)}{m^1+R_-(0)}\frac{V_+(0)}{r+rV_-(0)}\;\Bigg|\;\bL_{\cH^2_1}\right]\leq\EE\left[\frac{R_-(0)}{m^1+R_-(0)}\;\Bigg|\;\bL_{\cH^2_1}\right]\EE\left[\frac{V_+(0)}{r+rV_-(0)}\right],
\]
where the last factor is unconditional, $V_+(0)$ and $V_-(0)$ being independent of $\bL_{\cH^2_1}$. Since $\EE[V_+(0)/(r+rV_-(0))]\leq1$, it follows that
\begin{align*}
    \FDR &\leq \EE\left[\alpha\frac{|\cH^1_0|}{m^1+R_-(0)}+\alpha\frac{R_-(0)}{m^1+R_-(0)}\frac{V_+(0)}{r+rV_-(0)}\right] \\
    &= \EE\left[\EE\left[\alpha\frac{|\cH^1_0|}{m^1+R_-(0)}+\alpha\frac{R_-(0)}{m^1+R_-(0)}\frac{V_+(0)}{r+rV_-(0)}\;\Bigg|\;\bL_{\cH^2_1}\right]\right] \\
    &\leq \EE\left[\EE\left[\alpha\frac{|\cH^1_0|}{m^1+R_-(0)}+\alpha\frac{R_-(0)}{m^1+R_-(0)}\;\Bigg|\;\bL_{\cH^2_1}\right]\right]\leq\alpha,
\end{align*}
where the last inequality uses $|\cH^1_0|+R_-(0)\leq m^1+R_-(0)$.
\end{proof}

\section{Proofs for Lemmas in the Supplementary Material}
\label{asec:lemmas}

\begin{proof}[Proof of Lemma~\ref{slem:overvariable}]
First, for a fixed threshold $t$,
\begin{align*}
    &\PP\left(\delta(t)=k_3,V_+(t)=k_1,V_-(t)=k_2\right)\\
    &= \sum_{l\in \Omega(k_1,k_2,k_3)}\sum_{j_1,\dots,j_{k_1+k_2}\in\cH_0}\\
    &\qquad \PP\left(W_{j_1}>\dots>W_{j_{k_1+k_2}}>t>W_{\max}(j_1,\dots,j_{k_1+k_2})\right)\\
    &\qquad \times \PP\left(L_{j_s}=l_s;s=1,\dots,k_1+k_2\right)
\end{align*}
by the independence of $\bW$ and $\bL_{\cH_0}$, where
\[
W_{\max}(j_1,\dots,j_{k_1+k_2})=\max\left\{W_{j}:j\in\cH_0\setminus\{j_s:s=1,\dots,k_1+k_2\}\right\}
\]
and
\[
\Omega(k_1,k_2,k_3)=\left\{l\in\{0,1\}^{k_1+k_2}:\sum_{s=1}^{k_1+k_2}l_s=k_1,\sum_{s=1}^{k_3}l_s=k_3,l_{k_3+1}=0\right\}.
\]
Similarly,
\begin{align*}
    &\PP\left(V_+(t)=k_1,V_-(t)=k_2\right)\\
    &= \sum_{l\in \Omega(k_1,k_2)}\sum_{j_1,\dots,j_{k_1+k_2}\in\cH_0}\\
    &\qquad \PP\left(W_{j_1}>\dots>W_{j_{k_1+k_2}}>t>W_{\max}(j_1,\dots,j_{k_1+k_2})\right)\\
    &\qquad \times \PP\left(L_{j_s}=l_s;s=1,\dots,k_1+k_2\right)
\end{align*}
where
\[
\Omega(k_1,k_2)=\left\{l\in\{0,1\}^{k_1+k_2}:\sum_{s=1}^{k_1+k_2}l_s=k_1\right\}=\bigcup_{k_3=0}^{k_1}\Omega(k_1,k_2,k_3).
\]

By the definition of competition statistics, for any $l\in\Omega(k_1,k_2,k_3)$ and $j_1,\dots,j_{k_1+k_2}\in\cH_0$,
\[
\PP\left(L_{j_s}=l_s;s=1,\dots,k_1+k_2\right)=\left(\frac{r}{r+1}\right)^{k_1}\left(\frac{1}{r+1}\right)^{k_2}.
\]
Therefore,
\begin{align*}
    &\PP\left(\delta(t)=k_3\,\middle|\,V_+(t)=k_1,V_-(t)=k_2\right)
      = \frac{\PP\left(\delta(t)=k_3,V_+(t)=k_1,V_-(t)=k_2\right)}
             {\PP\left(V_+(t)=k_1,V_-(t)=k_2\right)}\\[4pt]
    &= \frac{\displaystyle\sum_{l\in \Omega(k_1,k_2,k_3)}
                     \sum_{j_1,\dots,j_{k_1+k_2}\in\cH_0}
                     \PP\left(W_{j_1}>\dots>W_{j_{k_1+k_2}}>t>W_{\max}(j_1,\dots,j_{k_1+k_2})\right)}
            {\displaystyle\sum_{l\in \Omega(k_1,k_2)}
                     \sum_{j_1,\dots,j_{k_1+k_2}\in\cH_0}
                     \PP\left(W_{j_1}>\dots>W_{j_{k_1+k_2}}>t>W_{\max}(j_1,\dots,j_{k_1+k_2})\right)}\\[4pt]
    &= \frac{|\Omega(k_1,k_2,k_3)|}{|\Omega(k_1,k_2)|}
       = \frac{\binom{k_1}{k_3}\cdot k_3!\cdot k_2\cdot (k_1-k_3+k_2-1)!}{(k_1+k_2)!}\\
    &= \binom{k_1-k_3+k_2-1}{k_2-1}\Bigg/\binom{k_1+k_2}{k_2}\\[4pt]
    &= f(k_3;k_1,k_2),
\end{align*}
which is independent of $t$.

Next we show that the same result holds for the data-driven threshold $T$.
\begin{align*}
    &\PP\left(\delta(T)=k_3,V_+(T)=k_1,V_-(T)=k_2\right)\\
    &= \sum_{t}\PP\left(\delta(t)=k_3,V_+(t)=k_1,V_-(t)=k_2,T=t\right)\\
    &= \sum_{t}\PP\left(T=t\;\middle|\;\delta(t)=k_3,V_+(t)=k_1,V_-(t)=k_2\right)
       \PP\left(\delta(t)=k_3,V_+(t)=k_1,V_-(t)=k_2\right)\\
    &= \sum_{t}\PP\left(T=t\;\middle|\;V_+(t)=k_1,V_-(t)=k_2\right)
       \PP\left(\delta(t)=k_3,V_+(t)=k_1,V_-(t)=k_2\right)\\
    &= f(k_3;k_1,k_2)\,\PP\left(V_+(T)=k_1,V_-(T)=k_2\right).
\end{align*}
Hence,
\[
\PP\left(\delta(T)=k\mid V_+(T),V_-(T)\right)=f\left(k;V_+(T),V_-(T)\right).
\]
With the exact conditional distribution known, we can compute the conditional expectation by a combinatorial argument. Treat $\{j\in\cH_0:W_j\geq T,\,L_j=0\}$ as $V_-(T)$ separation points that partition $\{j\in\cH_0:W_j\geq T,\,L_j=1\}$ into $V_-(T)+1$ bins. The variable $\delta(T)$ counts the number of labelled null variables in the first bin. Because the permutation of null variables is uniformly random, the expected count in each bin is the same. Since the $V_-(T)+1$ bins partition $V_+(T)$ total labelled null variables,
\[
\EE\left[\delta(T)\mid V_+(T),V_-(T)\right]=\frac{V_+(T)}{1+V_-(T)}.
\]
\end{proof}

\begin{proof}[Proof of Lemma~\ref{slem:biovervariable}]
The first fact is obvious.

If $\ind\{W_i\geq T,L_i=0\}\ind\{W_k\geq T_i,L_k=0\}=1$, for any $t\leq T_i$,
\begin{align*}
    \sum_{j\in\cH\setminus\{i\}}\ind\{W_j\geq t,L_j=1\}&=\sum_{j\in\cH}\ind\{W_j\geq t,L_j=1\}=\sum_{j\in\cH\setminus\{k\}}\ind\{W_j\geq t,L_j=1\},\\
    \sum_{j\in\cH\setminus\{i\}}\ind\{W_j\geq t,L_j=0\}&=\sum_{j\in\cH}\ind\{W_j\geq t,L_j=0\}-1=\sum_{j\in\cH\setminus\{k\}}\ind\{W_j\geq t,L_j=0\},
\end{align*}
which implies $T_k\leq T_i$ and $T_k\geq T_i$, so $T_k=T_i$.

If $\ind\{W_k\geq T_i,L_k=0\}=0$, for any $t$,
\begin{align*}
    \sum_{j\in\cH\setminus\{k\}}\ind\{W_j\geq t,L_j=1\}+1&\geq\sum_{j\in\cH}\ind\{W_j\geq t,L_j=1\}\geq\sum_{j\in\cH\setminus\{k\}}\ind\{W_j\geq t,L_j=1\},\\
    \sum_{j\in\cH}\ind\{W_j\geq t,L_j=0\}-1&\leq\sum_{j\in\cH\setminus\{k\}}\ind\{W_j\geq t,L_j=0\}\leq\sum_{j\in\cH}\ind\{W_j\geq t,L_j=0\},
\end{align*}
and the right part of the second inequality is equal when $t\geq T_i$, so $T_k\geq T$ and $T_k\leq T$ naturally, hence $T=T_k$.
\end{proof}

\begin{proof}[Proof of Lemma~\ref{slem:level_dim}]
The BH threshold for $\bP$ at level $m\alpha/|S|$ is $q_k^{m,m\alpha/|S|}=k(m\alpha/|S|)/m=k\alpha/|S|$. Since $P_j=+\infty$ for all $j\notin S$, the order statistics satisfy $P_{(k)}=P_{S,(k)}$ for $k\leq|S|$ and $P_{(k)}=+\infty$ for $k>|S|$. Hence the BH rejection count for $\bP$ is
\[
K=\max\left\{k: P_{(k)}\leq\frac{k\alpha}{|S|}\right\}
  =\max\left\{k\leq|S|: P_{S,(k)}\leq\frac{k\alpha}{|S|}\right\},
\]
which coincides with the rejection count for $\bP_S$ at level $\alpha$. Therefore
\[
\cR_{m\alpha/|S|}(\bP)=\left\{j:P_j\leq\frac{K\alpha}{|S|}\right\}
=\left\{j\in S:P_j\leq\frac{K\alpha}{|S|}\right\}
=\cR_{\alpha}(\bP_S).
\]
\end{proof}

\begin{proof}[Proof of Lemma~\ref{slem:rank_indep}]
Condition on $\bX=(X_1,\dots,X_m),M=m_0\in\{1,\dots,m-1\}$. Given $(\bX,M=m_0)$, the set $I$ is uniform over the $m_0$-subsets of $[m]$, and $j$ is uniform over $[m]\setminus I$. For $r=0,\dots,m_0$, let $\rho_j$ denote the rank of $X_j$ among $(X_1,\dots,X_m)$;
\[
\PP(R=1+r\mid\bX,M=m_0)
=\frac{1}{m}\sum_{j=1}^{m}\frac{\binom{\rho_j-1}{r}\binom{m-\rho_j}{m_0-r}}{\binom{m-1}{m_0}}.
\]
Since the $X_j$'s are continuous, they are distinct almost surely, so $\{\rho_j:j\in[m]\}=\{1,\dots,m\}$. By the hockey-stick identity
\[
\sum_{\rho=1}^{m}\binom{\rho-1}{r}\binom{m-\rho}{m_0-r}=\binom{m}{m_0+1},
\]
we obtain
\[
\PP(R=1+r\mid\bX,M=m_0)
=\frac{1}{m}\frac{\binom{m}{m_0+1}}{\binom{m-1}{m_0}}=\frac{1}{m_0+1},
\]
independently of $r$ and of $\bX$. Hence $R$ is uniform on $\{1,\dots,M+1\}$ conditionally on $M$, and in particular $R$ is independent of $\bX$.
\end{proof}

\begin{proof}[Proof of Lemma~\ref{slem:fdr_eq}]
The map $(\bW,\bL)\mapsto|\cH_0\cap\cR^{\cp}_{r,\alpha}(\bW,\bL)|/(1\vee|\cR^{\cp}_{r,\alpha}(\bW,\bL)|)$ is a deterministic functional of $(\bW,\bL)$, so its expectation depends on the pair only through its distribution. Since $(\bW,\bL)\overset{d}{=}(\bW',\bL')$, the two expectations coincide.
\end{proof}

\begin{proof}[Proof of Lemma~\ref{slem:rejection}]
We first show that $K\leq K'$. It is obvious that $K$ equals the number of rejections:
\[
K=\sum_{j\in\cH}\ind\{P_j\leq q^{m,\alpha}_K\},
\]
so $K\leq K'$. We next show that $K'\leq K$. Because 
\[
K'\leq\sum_{j\in\cH}\ind\{P_j\leq q^{m,\alpha}_{K'}\},
\]
we have
\[
P_{(K')}\leq P_{(\sum_{j\in\cH}\ind\{P_j\leq q^{m,\alpha}_{K'}\})}\leq q^{m,\alpha}_{K'},
\]
so $K'\leq K$. Therefore, $K=K'$.
\end{proof}

\section{More Numerical Results}
\label{app:morefig}

\subsection{More Analysis of the Prescreen Test}

\begin{figure}[htbp]
\centering
\includegraphics[width=369.11pt]{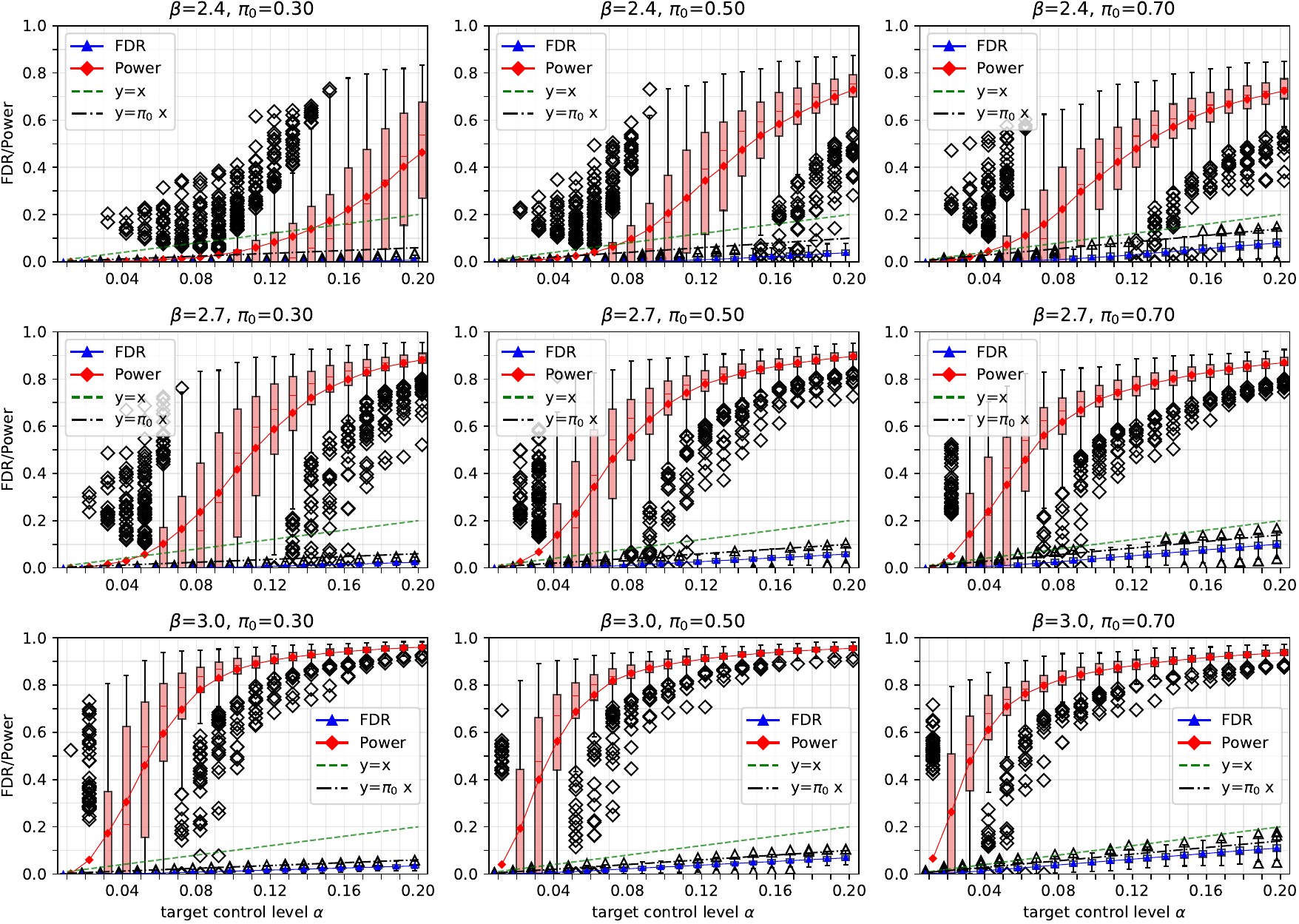}
\caption{The box plots of FDR under different signal strength $\beta$, null proportion $\pi_0$, and target level setting $\alpha$.}
\label{fig:prescreen_cal_box}
\end{figure}

We report additional simulation results in Figure~\ref{fig:prescreen_cal_box}, complementing Figure 1 of the main text.

\begin{figure}[htbp]
\centering
\includegraphics[width=271.34pt]{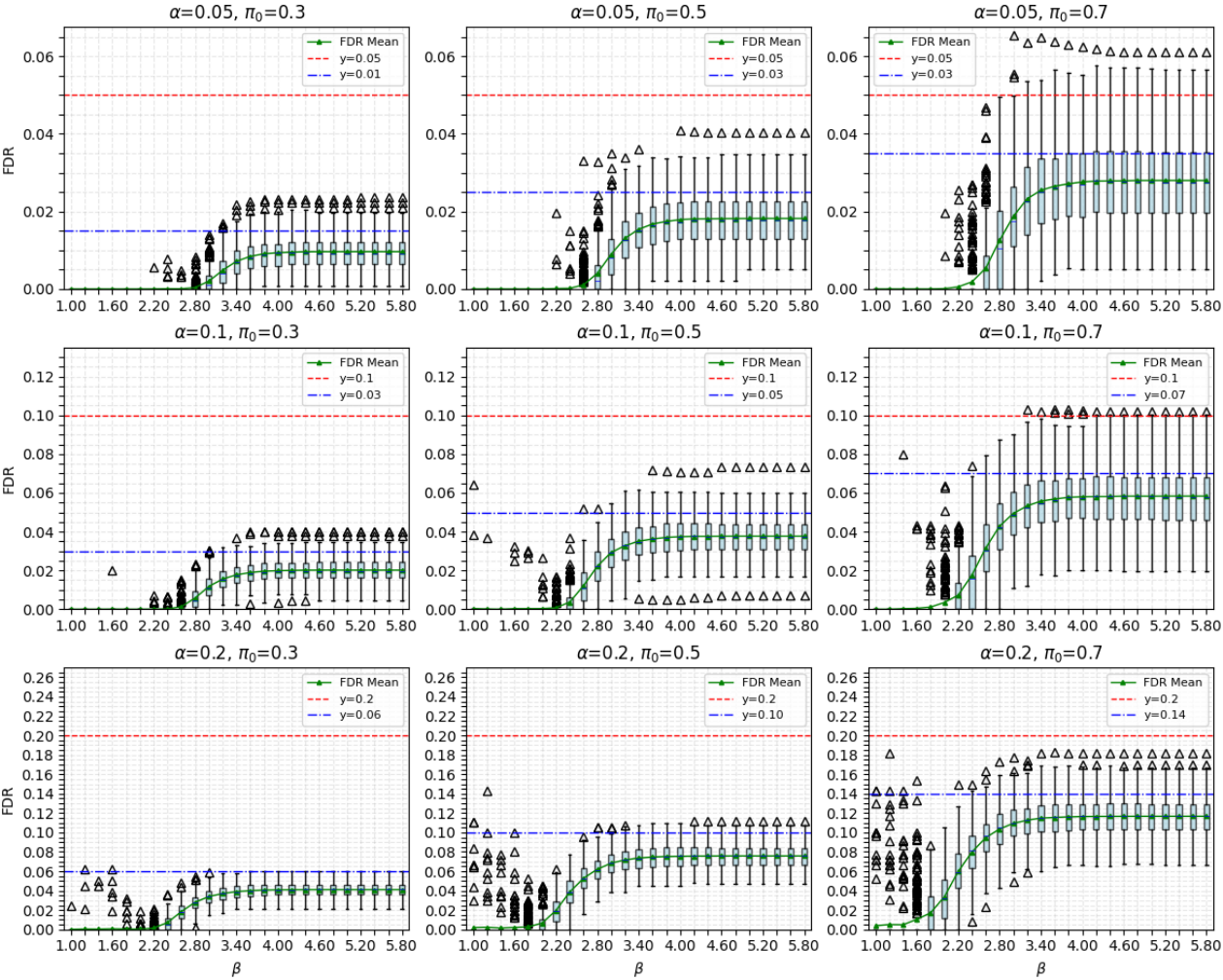}
\caption{The box plots of FDR under different signal strength $\beta$, null proportion $\pi_0$, and target level setting $\alpha$.}
\label{fig:prescreen_cal_focus_fdr}
\end{figure}

An important caveat of the conformal BH procedure with a calibration set constructed by prescreening is that the prescreening step need not be perfectly accurate: the resulting calibration set $\cD^{\rcal}=\cD^{\rcal}_0+\cD^{\rcal}_1$ is only partially exchangeable with the null test statistics. Exchangeability holds between the null hypotheses and the samples in $\cD^{\rcal}_0$, but not for those in $\cD^{\rcal}_1$. By the dominated BH formulation in the main text, the contaminated samples enter the procedure only through the weight $(1+|\cD^{\rcal}|)/(1+|\cD^{\rcal}_0|)$, so that FDR control is preserved whenever the exchangeable subset is sufficiently large. To illustrate this behaviour empirically, we examine the FDR under the same experimental setting and report the results in Figure~\ref{fig:prescreen_cal_focus_fdr}.

\begin{figure}[htbp]
\centering
\includegraphics[width=334.48pt]{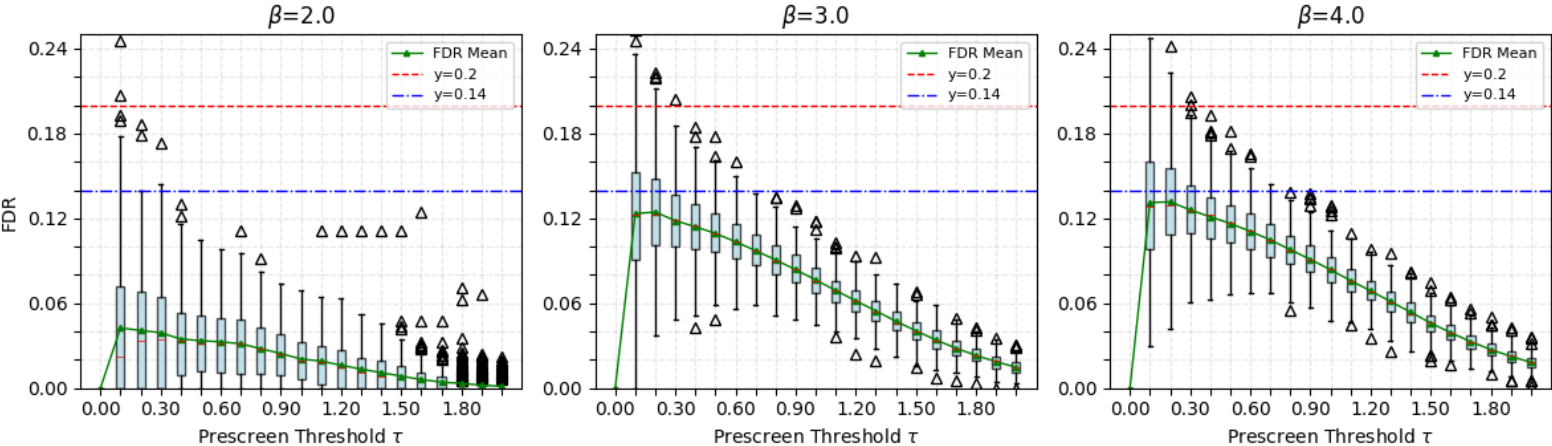}
\caption{The box plots of FDR under different prescreening thresholds and signal strength setting $\beta$.}
\label{fig:prescreen_cal_focus_threshold}
\end{figure}

We provide a theoretical analysis of the FDR control of the prescreen test.
\[
\begin{aligned}
\FDR&=\EE\left[\EE\left[\frac{|\cR_{\alpha}(\bP(\bX_{\cD^{\test}},\bX_{\cD^{\rcal}}))\cap\cD^{\test}\cap\cH_0|}{1\vee|\cR_{\alpha}(\bP(\bX_{\cD^{\test}},\bX_{\cD^{\rcal}}))|}\mid\bX^{\mathrm{pre}}\right]\right]\\
&\leq\alpha\,\EE_{\bX^{\mathrm{pre}}}\left[\frac{|\cD^{\test}\cap\cH_0|}{1\vee|\cD^{\test}|}\cdot\frac{1+|\cD^{\rcal}|}{1+|\cD^{\rcal}_0|}\right].
\end{aligned}
\]
This bound differs slightly from the corresponding bound for the dominated BH procedure in Section~3.2 of the main text, because it concerns a random sub-problem determined by the prescreening process. Notice that
\[
\frac{|\cD^{\test}\cap\cH_0|}{1\vee|\cD^{\test}|}\cdot\frac{1+|\cD^{\rcal}|}{1+|\cD^{\rcal}_0|}
=\frac{|\cH_0|-|\cD^{\rcal}_0|}{1\vee(|\cH|-|\cD^{\rcal}|)}\cdot\frac{1+|\cD^{\rcal}|}{1+|\cD^{\rcal}_0|}.
\]
Without loss of generality, we assume that the prescreening process ensures that $\cH\setminus\cD^{\rcal}$ is always nonempty (for example, by guaranteeing that the most significant hypothesis is never screened out). Writing $|\cH|=h$, $|\cD^{\rcal}|=c$, and $|\cD^{\rcal}_0|=c_0$, we have for any $\varepsilon\geq 0$ with $c_0\geq(1+\varepsilon)\pi_0c$ that
\[
\frac{|\cH_0|-|\cD^{\rcal}_0|}{1\vee(|\cH|-|\cD^{\rcal}|)}\cdot\frac{1+|\cD^{\rcal}|}{1+|\cD^{\rcal}_0|}
=\frac{\pi_0h-c_0}{h-c}\cdot\frac{1+c}{1+c_0}
\leq\pi_0\cdot\frac{1}{(1+\varepsilon)\pi_0}=\frac{1}{1+\varepsilon},
\]
and consequently
\[
\begin{aligned}
\FDR&\leq\inf_{\varepsilon\geq0}\Biggl\{\frac{\alpha}{1+\varepsilon}\,\PP_{\bX^{\mathrm{pre}}}\left\{|\cD^{\rcal}_0|\geq(1+\varepsilon)\,\pi_0\,|\cD^{\rcal}|\right\}\\
&\qquad+\PP_{\bX^{\mathrm{pre}}}\left\{|\cD^{\rcal}_0|<(1+\varepsilon)\,\pi_0\,|\cD^{\rcal}|\right\}\Biggr\}.
\end{aligned}
\]
This inequality follows from the data splitting method, which ensures that the samples used for testing and those used for prescreening are independent. The second term can be made sufficiently small by choosing an appropriate screening procedure. To illustrate the effect of the prescreening threshold on FDR, we fix the signal strength to $\beta\in\{2,3,4\}$, set the target significance level to $\alpha=0.2$ and the null proportion to $\pi_0=0.7$, and define the prescreened calibration set as $\cD^{\rcal}=\{i:|X_i^{\mathrm{pre}}|<\tau\}$. The results are reported in Figure~\ref{fig:prescreen_cal_focus_threshold}.
In practice, whether prescreening compromises FDR control should be examined case by case. However, the results above suggest that in most routine cases this concern can be disregarded. A reasonable prescreening procedure tends to first screen out null hypotheses and then non-null ones, so that the ratio of the two remains below the null proportion $\pi_0$. When this ratio approaches $\pi_0$, most of the nulls have been screened out, which makes the first fractional term in the bound extremely small and the resulting test more conservative.

\subsection{Effect of Null p-Value Heterogeneity on Power}

For each Monte Carlo repetition we draw a Gaussian vector $X \sim N_m(0, I_m)$
of dimension $m = m_1 + m_{0,\mathrm{test}} + m_{0,\mathrm{cal}}$, where the
first $m_0 = m_{0,\mathrm{test}} + m_{0,\mathrm{cal}}$ coordinates are null
($X_i \sim N(0,1)$) and the remaining $m_1$ coordinates are signals
($X_i \sim N(\beta,1)$). A periodic weight vector
$\boldsymbol{w} = (w_1,\dots,w_L)$ is applied to the null coordinates only, so
that the null coordinates inherit a heterogeneous scale while the signal
coordinates keep unit weight. Exactly $m_{0,\mathrm{cal}}$ of the $m_0$ null
coordinates are drawn uniformly at random into the calibration set $\cD^{\rcal}$;
the remaining $m_{0,\mathrm{test}}$ null coordinates together with the $m_1$
signal coordinates form the test set $\cD^{\test}$.

Reported FDR and power are the averages over $n_{\mathrm{rep}} = 300$ repetitions, at each of 20 levels $\alpha \in \{0.01,\dots,0.20\}$.

Two complementary schemes generate the weight vector, both with $L \in \{1,3,5,7\}$; see Table~\ref{tab:weight-schemes}.

\begin{table}[htbp]
\centering
\caption{Weight-vector schemes. EXTEND widens the weight range at constant step $0.2$; INSERT keeps the range $[0.4,1.6]$ fixed and adds interior points. Both reduce to the unweighted case at $L=1$ and coincide at $L=7$.}
\label{tab:weight-schemes}
\begin{tabular}{@{}lll@{}}
\toprule
$L$ & EXTEND (step $0.2$, range widens) & INSERT (range $[0.4,1.6]$ fixed) \\
\midrule
1 & $[1.0]$                                        & $[1.0]$ \\
3 & $[0.8, 1.0, 1.2]$                              & $[0.4, 1.0, 1.6]$ \\
5 & $[0.6, 0.8, 1.0, 1.2, 1.4]$                    & $[0.4, 0.7, 1.0, 1.3, 1.6]$ \\
7 & $[0.4, 0.6, 0.8, 1.0, 1.2, 1.4, 1.6]$          & $[0.4, 0.6, 0.8, 1.0, 1.2, 1.4, 1.6]$ \\
\bottomrule
\end{tabular}
\end{table}

The remaining parameters are fixed at $m_1 = 500$, $m_{0,\mathrm{test}} = 1500$, $m_{0,\mathrm{cal}} = 1000$, and $\beta = 2.4$. EXTEND thus widens the weight \emph{range} as $L$ grows, whereas INSERT adds interior weights within the fixed range $[0.4,1.6]$, isolating these two mechanisms of null p-value heterogeneity.

\begin{figure}[htbp]
\centering
\includegraphics[width=378.01pt]{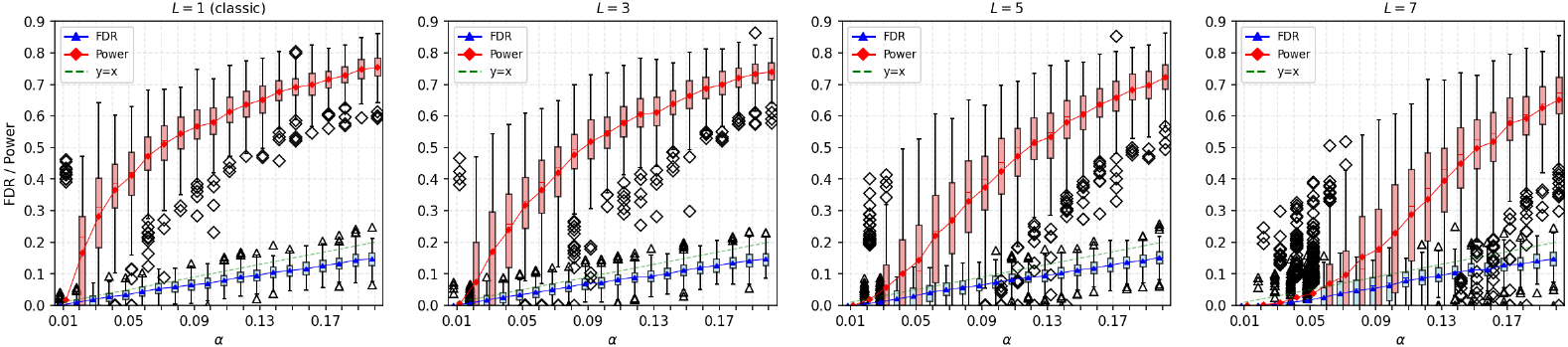}
\caption{The box plots of FDR and power versus $\alpha$ for the EXTEND scheme under different weight-vector lengths $L$.}
\label{fig:hetero-extend}
\end{figure}

\begin{figure}[htbp]
\centering
\includegraphics[width=378.01pt]{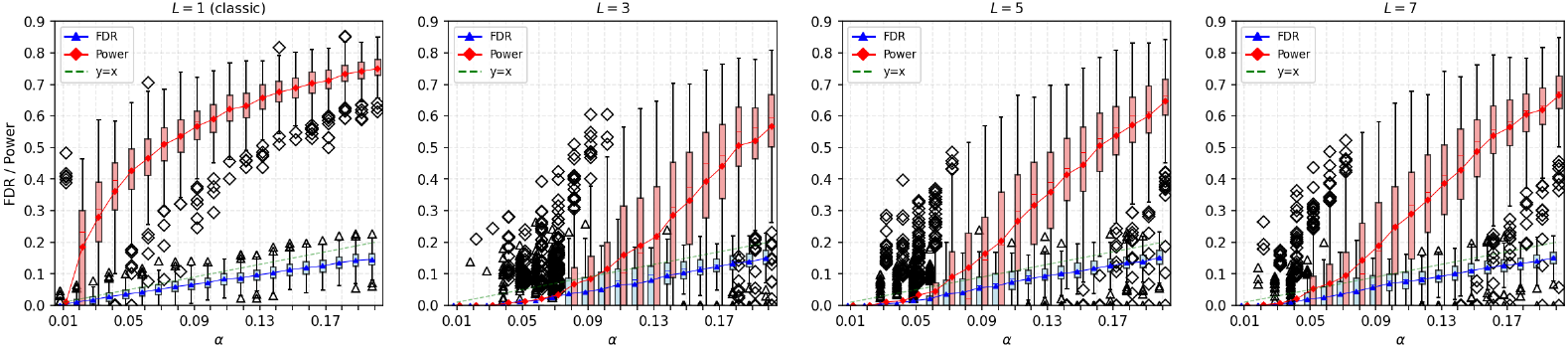}
\caption{The box plots of FDR and power versus $\alpha$ for the INSERT scheme with fixed weight range $[0.4,1.6]$ under different lengths $L$.}
\label{fig:hetero-insert}
\end{figure}

\begin{table}[htbp]
\centering
\caption{Averaged FDR / power at selected levels $\alpha$
($n_{\mathrm{rep}} = 300$). The $L=1$ rows coincide across schemes, and the two schemes share the same $L=7$ weight vector.}
\label{tab:hetero-summary}
\begin{tabular}{@{}llcccc@{}}
\toprule
Scheme & $L$ & $0.05$ & $0.10$ & $0.15$ & $0.20$ \\
\midrule
EXTEND & 1 & $0.034 / 0.417$ & $0.072 / 0.588$ & $0.111 / 0.683$ & $0.148 / 0.753$ \\
EXTEND & 3 & $0.032 / 0.306$ & $0.074 / 0.544$ & $0.110 / 0.663$ & $0.146 / 0.739$ \\
EXTEND & 5 & $0.030 / 0.165$ & $0.069 / 0.410$ & $0.111 / 0.598$ & $0.148 / 0.711$ \\
EXTEND & 7 & $0.017 / 0.041$ & $0.063 / 0.241$ & $0.108 / 0.486$ & $0.151 / 0.654$ \\
\midrule
INSERT & 1 & $0.035 / 0.417$ & $0.074 / 0.587$ & $0.111 / 0.689$ & $0.149 / 0.755$ \\
INSERT & 3 & $0.008 / 0.009$ & $0.055 / 0.121$ & $0.102 / 0.335$ & $0.147 / 0.572$ \\
INSERT & 5 & $0.015 / 0.034$ & $0.069 / 0.208$ & $0.110 / 0.466$ & $0.149 / 0.643$ \\
INSERT & 7 & $0.021 / 0.048$ & $0.067 / 0.243$ & $0.111 / 0.490$ & $0.147 / 0.651$ \\
\bottomrule
\end{tabular}
\end{table}

Figures~\ref{fig:hetero-extend} and~\ref{fig:hetero-insert}, together with Table~\ref{tab:hetero-summary}, show that FDR is controlled at the nominal level for every $L$ in both schemes, consistent with the validity of the randomized conformal p-values. The two schemes differ markedly in power.

\textbf{EXTEND: monotone power loss.}
Power decreases monotonically with $L$: at $\alpha = 0.10$ it falls
$0.588 \to 0.544 \to 0.410 \to 0.241$ as $L$ runs $1 \to 3 \to 5 \to 7$, and the
same ordering holds at every level $\alpha$ in the table.

\textbf{INSERT: a non-monotone effect.}
With the range fixed at $[0.4,1.6]$, power behaves non-monotonically: it
collapses from $L=1$ to $L=3$ ($0.587 \to 0.121$ at $\alpha = 0.10$), the step
in which the extreme weights $0.4$ and $1.6$ are introduced, and then recovers
as $L$ grows further ($0.208$ at $L=5$, $0.243$ at $L=7$). The INSERT panel is
thus U-shaped in $L$, with the minimum at $L=3$.

\textbf{Key takeaway.}
The two schemes coincide at $L=7$, and their $L=7$ numbers agree within random error (power $0.243$ vs.\ $0.241$, FDR $0.067$ vs.\ $0.063$ at $\alpha=0.10$), a consistency check of the simulation. We conjecture that the contrast between the two panels is driven by how the periodic weights distort the null p-values: weights below $1$ shrink the corresponding null scores and push their conformal p-values towards small values, so these ``deflated'' null p-values crowd the left tail of the BH ordering and raise the threshold $p_{(K)}$, leaving fewer signals to survive. EXTEND widens the range and thereby exposes a growing fraction of the null coordinates to such extreme weights, aggravating the crowding; INSERT keeps the range fixed and only dilutes the extreme weights---each acts on $1/L$ of the null coordinates---so the crowding weakens as $L$ grows and power recovers. Power is therefore governed by the \emph{span} of the periodic weights and the \emph{fraction of null coordinates} covered by the extreme weights, rather than by the vector length $L$ per se.

\subsection{FDR and Power Comparison between Conformal and Competition Methods}

According to the analysis in this work, the conformal test and the competition test reside within a unified framework: the conformal test can be regarded as a special case of the competition test. This naturally raises the question of how the two tests differ when the requirements of both are satisfied.

\begin{figure}[htbp]
\centering
\includegraphics[width=336.1pt]{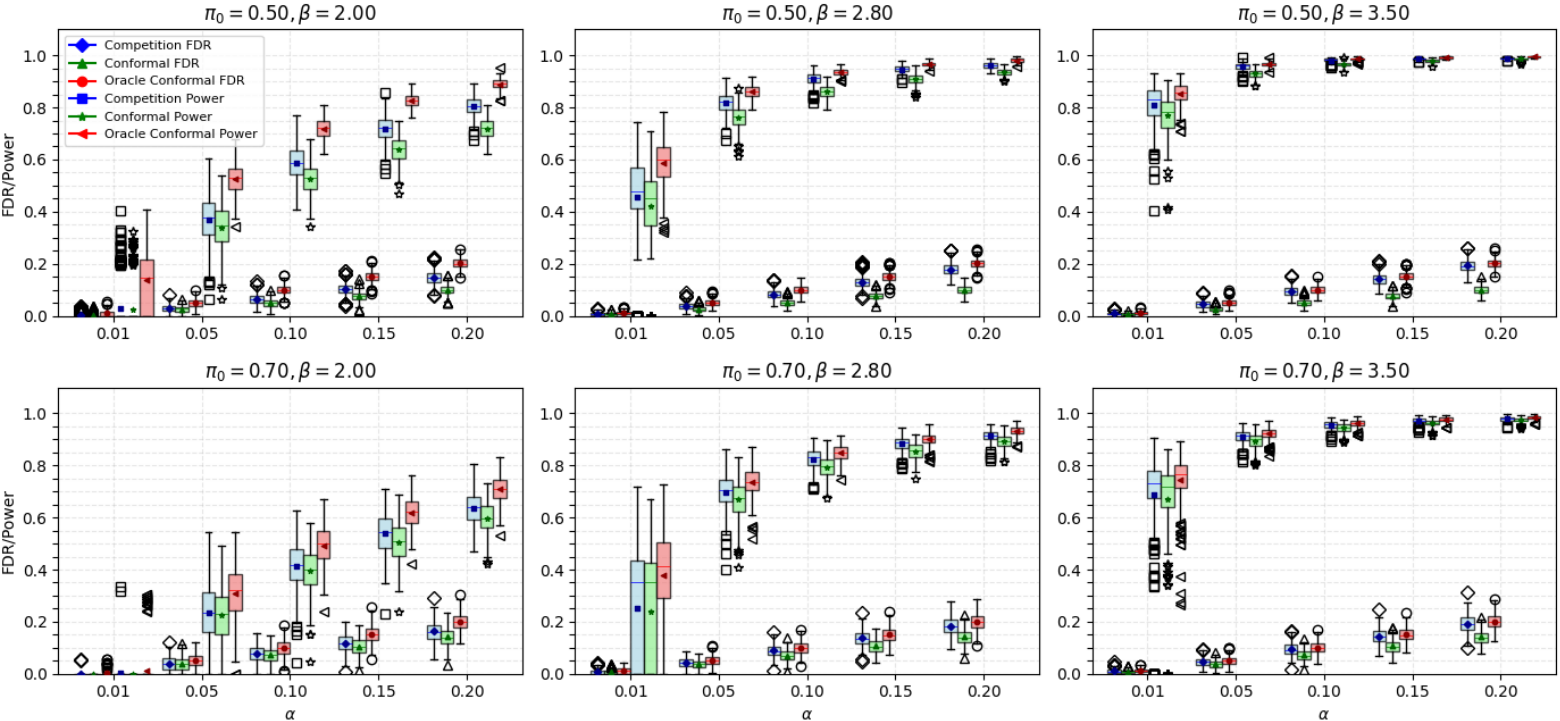}
\caption{The box plots of FDR and power under different target control levels $\alpha$ and parameter setting $\beta,\pi_0$.}
\label{fig:cp_vs_cf_vs_ocf_box}
\end{figure}

\begin{figure}[htbp]
\centering
\includegraphics[width=336.09pt]{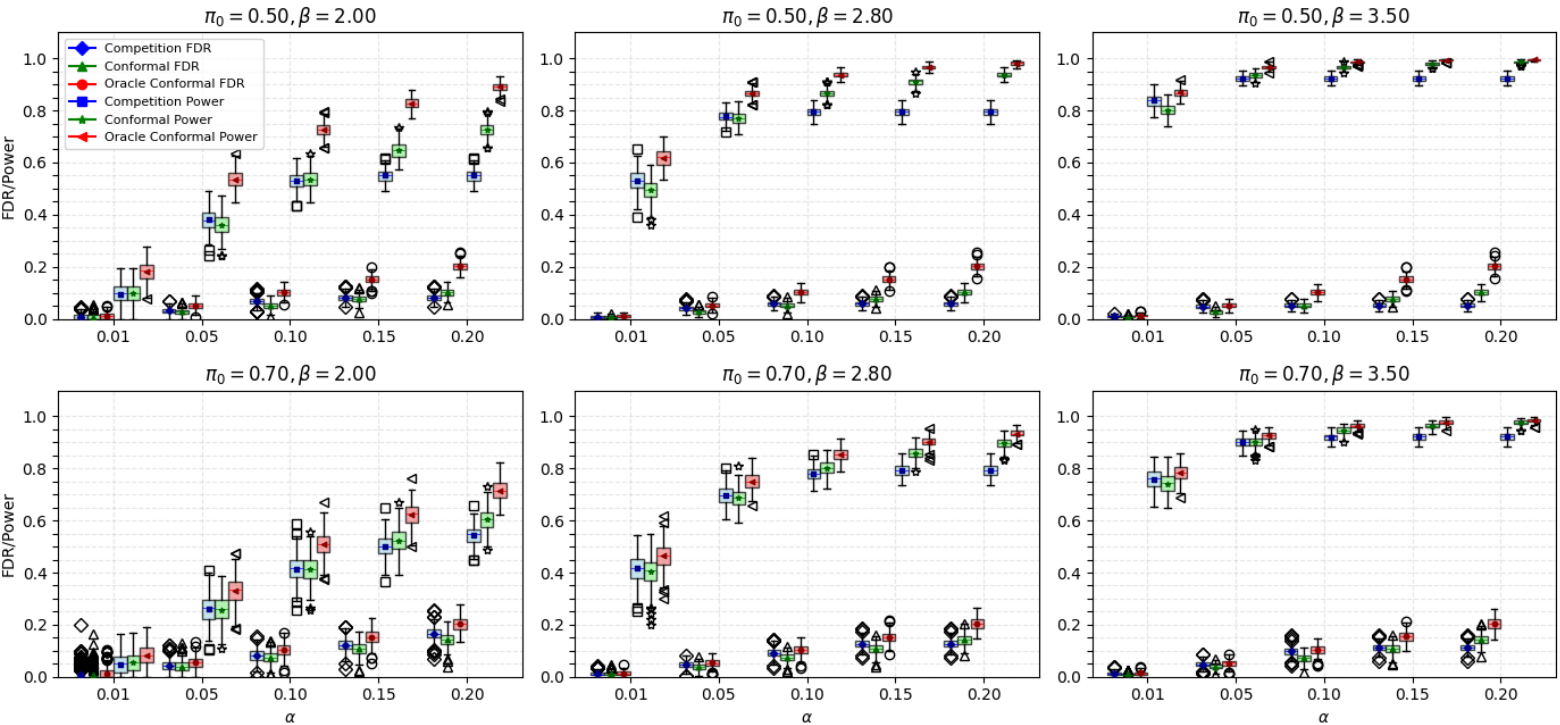}
\caption{The box plots of FDR and power under different target control levels $\alpha$ and parameter setting $\beta,\pi_0$.}
\label{fig:multi_cp_vs_cf_vs_ocf_box}
\end{figure}

\begin{figure}[htbp]
\centering
\includegraphics[width=336.36pt]{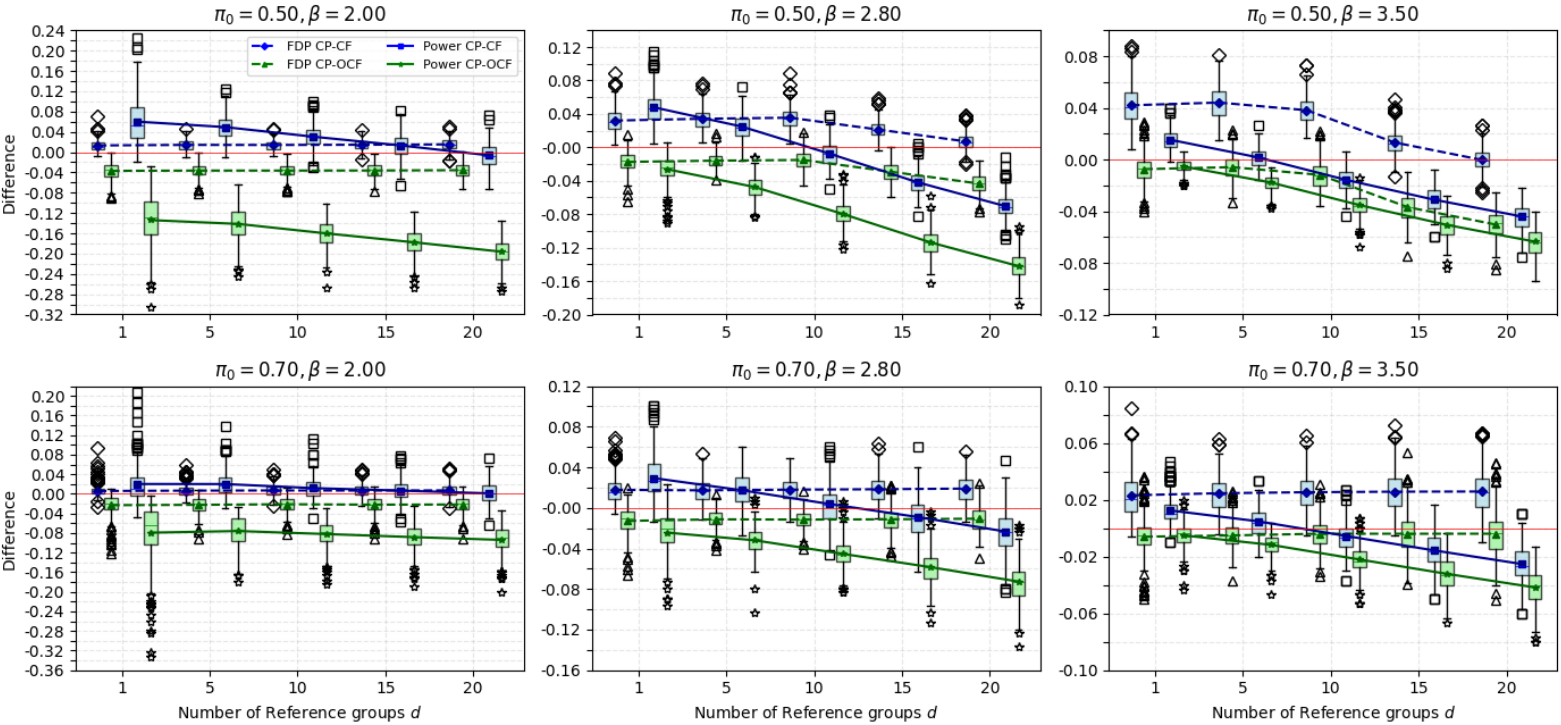}
\caption{The box plots of the differences in FDR and power under different reference set sizes $d$ and parameter setting $\beta,\pi_0$.}
\label{fig:multi_cp_cf_ocf_gap_box}
\end{figure}

We design the following experiment. Let $X_i\sim N(\mu_i,1)$, $i=1,\dots,m$, be the test statistics, with hypothesis set $\cH=[m]$ and null set $\cH_0=[m_0]$. If $i\in\cH_0$, set $\mu_i=0$; if $i\in\cH_1$, set $\mu_i=\beta$, where $\beta$ is the signal strength. In addition, generate reference variables $X_i^{\mathrm{ref}}\sim N(0,1)$, $i=1,\dots,m$, independently of the $X_i$'s. Consider three tests:
\begin{enumerate}
    \item \textbf{Competition test:} construct competition statistics $W_i=\max\{|X_i|,|X_i^{\mathrm{ref}}|\}$ and $L_i=\ind\{|X_i|\geq|X_i^{\mathrm{ref}}|\}$, yielding the rejection set $\cR^{\cp}_{1,\alpha}(\bW,\bL)$.
    \item \textbf{Conformal test:} construct conformal p-values $\bP(\bX,\bX^{\mathrm{ref}})$, yielding the rejection set $\cR_{\alpha}(\bP)$.
    \item \textbf{Oracle conformal test:} use conformal p-values as well, but at the oracle level $\alpha\pi_0^{-1}$, yielding the rejection set $\cR_{\alpha\pi_0^{-1}}(\bP)$, where $\pi_0=m^{-1}m_0$.
\end{enumerate}
We take $m=1000$ and run $300$ replications. The results are shown in Figure~\ref{fig:cp_vs_cf_vs_ocf_box}.
We next consider a larger reference set. Let $\bX^{\mathrm{ref},j}$, $j=1,\dots,d$, denote $d$ reference vectors, and write $\bX^{\mathrm{ref}}=(\bX^{\mathrm{ref},1},\dots,\bX^{\mathrm{ref},d})$. Among the three tests above, the conformal and oracle conformal tests remain unchanged, whereas the competition test now uses the asymmetry coefficient $r=d^{-1}$, yielding the rejection set $\cR^{\cp}_{d^{-1},\alpha}(\bW,\bL)$, where
\[
W_i=\max\bigl\{|X_i|,|X_i^{\mathrm{ref},j}|:j=1,\dots,d\bigr\}, \qquad L_i=\prod_{j=1}^{d}\ind\{|X_i|\geq|X_i^{\mathrm{ref},j}|\}.
\]
Taking $d=20$, the results are shown in Figure~\ref{fig:multi_cp_vs_cf_vs_ocf_box}.
With a larger reference set, the test outcomes change: the stability of both the conformal and the competition tests improves markedly. To examine the effect of the reference set size, we take the competition test as the baseline and compute its differences from the conformal test and from the oracle conformal test. With target level $\alpha=0.1$, the results are shown in Figure~\ref{fig:multi_cp_cf_ocf_gap_box}.
As the reference set size increases, the powers of the conformal and oracle conformal tests increase relative to that of the competition test, while the FDR also rises, though only mildly, and remains controlled throughout.

\subsection{Integrating More Groups of Hypotheses}

The algorithms and proofs presented in our work can be directly generalized from the integration of two groups to the integration of multiple groups. Here we only provide an example of the integration algorithm for multiple groups of randomized p-values.

\begin{algorithm}[htbp]
\caption{Integrated test of groups of randomized p-values}
\LinesNotNumbered
\label{salgo:multi_int_cp_cp}
\emph{Input:} $\bP^g\in\RR_*^{m^g}$, $g=1,\dots,d$, $\alpha\in(0,1)$. \emph{Output:} Rejection sets $\cR_g$, $g=1,\dots,d$.
\begin{enumerate}
  \item Set $m=\sum_{g=1}^{d}m^g$, and find the largest $K\in\{1,\dots,m\}$ such that
    \[
    \sum_{g=1}^{d}\sum_{j=1}^{m^g}\ind\{P^g_j\leq q^{m,\alpha}_K\}\geq K.
    \]
  \item Return $\cR_g=\{j\in[m^g]:P^g_j\leq q^{m,\alpha}_K\}$, $g=1,\dots,d$.
\end{enumerate}
\end{algorithm}
We state without proof that if $\bP^g$, $g=1,\dots,d$, are mutually independent randomized p-variables and each satisfies $\PRDS_{\cH^g}$, then this algorithm controls the FDR. Replacing $m$ by $m_0=\sum_{g=1}^d m_0^g$ yields the corresponding oracle algorithm.

\begin{figure}[htbp]
\centering
\includegraphics[width=336.1pt]{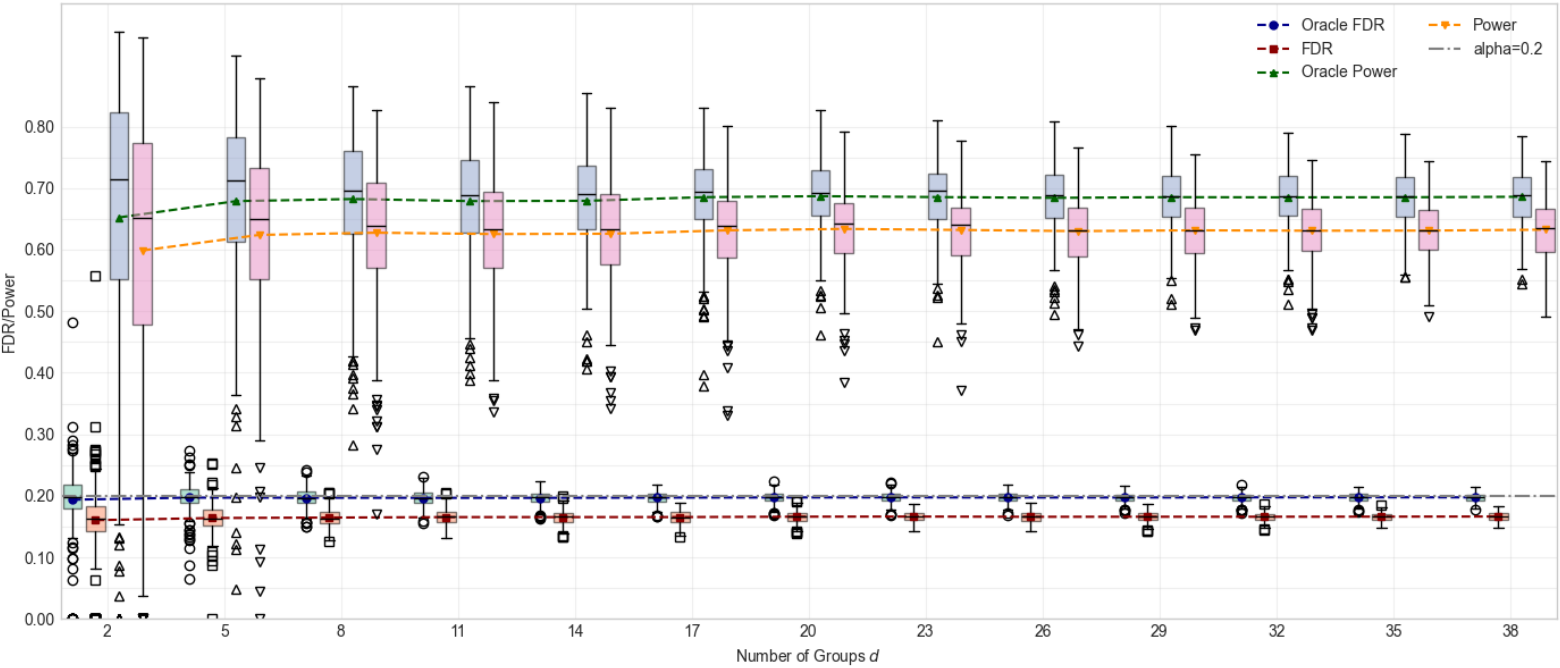}
\caption{The box plots of FDR and power under different numbers of groups $d$.}
\label{fig:multi_integration_weights_focus_multi}
\end{figure}

We conduct a simulation study to illustrate the performance of our algorithm when integrating multiple groups of hypotheses. For group $g$, we generate $m^g\sim\mathrm{Unif}(2000,6000)$, $\pi_0^g\sim\mathrm{Unif}(0.7,0.99)$, $m^{g,\rcal}=1000$, $\boldsymbol{w}^g\sim\mathrm{Unif}(0.5,1.5,b)$ with $b=3$, and $\beta^g\sim\mathrm{Unif}(2,3)$. These yield $d$ groups of independent randomized p-variables, where $m^g$ is the size of the test set, $m^{g,\rcal}$ the size of the calibration set, $\pi_0^g$ the proportion of true null hypotheses, $\boldsymbol{w}^g\in\RR_*^b$ the weight vector, and $\beta^g$ the signal strength. At target level $\alpha=0.2$, we run $300$ replications and estimate the FDR and power by the average FDP and TDP, respectively; Figure~\ref{fig:multi_integration_weights_focus_multi} reports the results for different numbers of groups $d$.
The simulation results show that our integration algorithm not only preserves FDR control and power but also improves the stability of the tests.

\subsection{The Influence of Weights in Integration}

When introducing the integration algorithm, we mentioned that different weights can be assigned to different groups in order to specifically reject more hypotheses from a particular group. Similarly, we only provide an example of the weighted integration algorithm for two groups of randomized p-variables.

\begin{figure}[htbp]
\centering
\includegraphics[width=334.21pt]{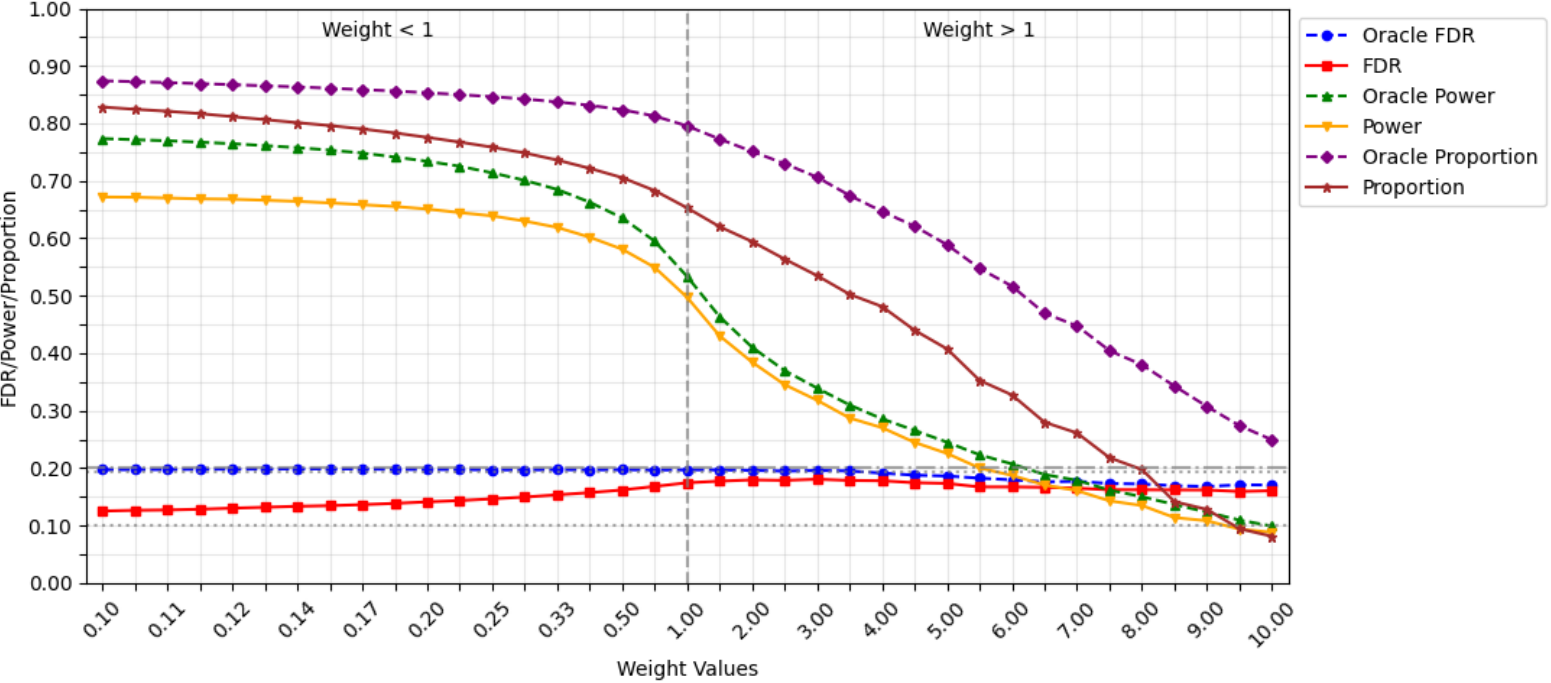}
\caption{The plots of FDR and power under different weights $\rho$.}
\label{fig:rp_rp_integration_weights_exps}
\end{figure}

\begin{algorithm}[htbp]
\caption{Weighted integrated test of randomized p-values}
\LinesNotNumbered
\label{salgo:weight_int_cp_cp}
\emph{Input:} $\bP^1\in\RR_*^{m^1}$, $\bP^2\in\RR_*^{m^2}$, $\alpha\in(0,1)$, $\rho\in\RR_+$. \emph{Output:} Rejection sets $\cR_1,\cR_2$.
\begin{enumerate}
  \item Find the largest $K\in\{0,1,\dots,m^1+m^2\}$ such that
    \[
    \sum_{j=1}^{m^1}\ind\{P^1_j\leq q^{m^1+\rho m^2,\alpha}_K\}
    +\sum_{j=1}^{m^2}\ind\{P^2_j/\rho\leq q^{m^1+\rho m^2,\alpha}_K\}\geq K.
    \]
  \item Return $\cR_1=\{j\in[m^1]:P^1_j\leq q^{m^1+\rho m^2,\alpha}_K\}$ and $\cR_2=\{j\in[m^2]:P^2_j/\rho\leq q^{m^1+\rho m^2,\alpha}_K\}$.
\end{enumerate}
\end{algorithm}
We also construct two independent sets of p-variables for Algorithm~\ref{salgo:weight_int_cp_cp}. Setting the weight vector $\boldsymbol{w}=1$ is degenerate, meaning that we generate standard p-variables rather than randomized p-variables. Table~\ref{tab:data_config} reports the parameter configurations.
\begin{table}[htbp]
\centering
\caption{The parameter setting of weighted integration.}
\label{tab:data_config}
\begin{tabular}{ccccc}
\toprule
&$m^g$ & $\pi_0^g$ & $m^{g,\rcal}$ & $\beta^g$ \\
\midrule
$g=1$ & 2000 & 1/2 & 1000 & 2.2 \\
$g=2$ & 9000 & 29/30 & 1000 & 2.2 \\
\bottomrule
\end{tabular}
\end{table}
Under these settings, at target level $\alpha=0.2$, we run $300$ replications and estimate the FDR and power by the average FDP and TDP, respectively, for different weights $\rho$; Figure~\ref{fig:rp_rp_integration_weights_exps} reports the results.

\bibliography{references}

\end{document}